\documentclass[11pt]{amsart}	
\usepackage[utf8]{inputenc}
\usepackage{amsmath,amsfonts,amssymb,dsfont}
\usepackage[foot]{amsaddr}
\usepackage{lmodern}
\usepackage{graphicx}
\usepackage[english]{babel}
\usepackage{listings}
\usepackage[a4paper]{geometry}
\usepackage[T1]{fontenc}
\usepackage{xspace}
\usepackage[colorinlistoftodos]{todonotes}
\usepackage{pgf,pgfkeys,pgfmath}
\usepackage{upgreek}
\usepackage{fullpage}
 
\usepackage[backref=page]{hyperref}
\hypersetup{
	plainpages=false,
	unicode=false,          
	pdftoolbar=true,        
	pdfmenubar=true,        
	pdffitwindow=true,      
	pdftitle={My title},    
	pdfauthor={Author},     
	pdfsubject={Subject},   
	pdfcreator={Creator},   
	pdfproducer={Producer}, 
	pdfkeywords={keywords}, 
	pdfnewwindow=true,      
	colorlinks=true,       
	linkcolor=blue,          
	citecolor=blue,        
	filecolor=blue,      
	urlcolor=magenta,          
	urlbordercolor=0 1 1
}
\newcommand*{\Scale}[2][4]{\scalebox{#1}{\ensuremath{#2}}}%
\newcommand{\bigmu}{\Scale[1.3]{\upmu}}
\def\div{\operatorname{div}}
\def\dist{\operatorname{dist}}

\newtheorem{thm}{Theorem}[section]
\newtheorem{rem}[thm]{Remark}
\newtheorem{prop}[thm]{Proposition}

\usepackage{todonotes}

\newcommand\numeq[1]%
  {\stackrel{\scriptscriptstyle(\mkern-1.5mu#1\mkern-1.5mu)}{=}}

\newcommand{\myqed}{\thinspace\null\nobreak\hfill\hbox{\vbox{\kern-.2pt\hrule
height.2pt depth.2pt\kern-.2pt\kern-.2pt \hbox to2.5mm{\kern-.2pt\vrule
width.4pt \kern-.2pt\raise2.5mm\vbox to.2pt{}\lower0pt\vtop to.2pt{}\hfil
\kern-.2pt \vrule width.4pt\kern-.2pt}\kern-.2pt\kern-.2pt\hrule
height.2pt depth.2pt \kern-.2pt}}\par\medbreak}

\definecolor{lbcolor}{rgb}{0.95,0.95,0.95}
\definecolor{cblue}{rgb}{0.,0.0,0.6}

\newcommand{\ordres}[2]{\ensuremath{(\mathcal{O}(#1), \mathcal{O}(#2) )}}

\definecolor{Lightgray}{rgb}{0.85, 0.85, 0.85}

\newcommand{\MCH}{{\bf M-CH}}
\newcommand{\NMNCH}{{\bf NMN-CH}}
\def\R{\mathbb{R}}

\def\one#1{\mathds{1}_{#1}}

\title{A multiphase Cahn-Hilliard system with mobilities\\ and the numerical simulation of dewetting}
\author{Elie Bretin$^1$}
\author{Roland Denis$^2$}
\author{Simon Masnou$^2$}
\author{Arnaud Sengers$^2$}
\author{Garry Terii$^2$}
\address{$^1$Univ Lyon, INSA de Lyon, CNRS UMR 5208, Institut Camille Jordan\\ 20 avenue Albert Einstein, F-69621 Villeurbanne, France}
\address{$^2$Univ Lyon, Universit\'e Claude Bernard Lyon 1, CNRS UMR 5208, Institut Camille Jordan, 43 boulevard du 11 novembre
1918, F-69622 Villeurbanne, France}
\email{elie.bretin@insa-lyon.fr}
\email{denis@math.univ-lyon1.fr}
\email{masnou@math.univ-lyon1.fr}
\email{sengers@math.univ-lyon1.fr}
\email{terii@math.univ-lyon1.fr}

\makeatletter
\@namedef{subjclassname@2020}{\textup{2020} Mathematics Subject Classification}
\makeatother
 \subjclass{74N20, 35A35, 53E10, 53E40, 65M32, 35A15}
\keywords{Phase field approximation, multiphase Cahn-Hilliard system, surface diffusion, degenerate mobilities, numerical approximation of dewetting.}

\begin{document}
\maketitle
\begin{abstract}
We propose in this paper a new multiphase Cahn-Hilliard model with doubly degenerate mobilities. 
We prove by a formal asymptotic analysis that it approximates with second order accuracy the multiphase surface diffusion flow with mobility coefficients and surface tensions. To illustrate that it lends itself well to numerical approximation, we propose a simple and effective numerical scheme together with a very compact Matlab implementation. We provide the results of various numerical experiments to show the influence of mobility and surface tension coefficients. Thanks to  its second order accuracy and its good suitability for numerical implementation, our model is very handy for tackling notably difficult surface diffusion problems. In particular, we show that it can be used very effectively to simulate numerically the {dewetting} of thin liquid tubes on arbitrary solid supports without requiring nonlinear boundary conditions.
\end{abstract}

\section{Introduction}
{This paper is devoted to the phase field approximation of multiphase surface diffusion with surface tensions and mobility coefficients. Surface diffusion is a natural process that makes interfaces evolve toward certain equilibrium configurations. Recall that, in the case of a liquid film covering a solid surface in ambient air and considering capillary effects only}, Young~\cite{young1805iii} identified in 1805 the optimal shape at rest of the liquid phase and 
proposed the following law for the contact angle $\theta$ between the liquid and the solid on the contact line where the three phases meet:
\begin{equation*}
\cos(\theta) = \frac{\sigma_{SV}-\sigma_{LS}}{\sigma_{VL}},
\end{equation*}
where $\sigma_{SV},\sigma_{LS},\sigma_{VL}$ represent the surface tensions of the  
solid-vapor  $\Gamma_{SV}$,  liquid-solid $\Gamma_{LS}$, and vapor-liquid $\Gamma_{VL}$ interfaces, respectively. 
Mathematically, Young's law can be derived by minimizing the total energy in the solid-liquid-vapor system.
Ignoring gravity, this total energy reads as
$$
\mathcal{E} = \sigma_{SV}\mathcal{H}^{d-1}(\Gamma_{SV})+\sigma_{LS}\mathcal{H}^{d-1}(\Gamma_{LS})+\sigma_{VL}\mathcal{H}^{d-1}(\Gamma_{VL}).
$$
 
\noindent which is a particular instance of the generic $L$-phase perimeter
\begin{equation} 
  P(\Omega_1,\dots, \Omega_L) = \frac{1}{2} \sum_{i,j=1}^{L}  \sigma_{i,j} \mathcal{H}^{d-1}(\Gamma_{i,j}), \label{per_multip_phase_intro}
  \end{equation}
where $\{\Omega_1,\dots,\Omega_L\}$ is an open partition of an open bounded domain $\Omega \subset \mathbb{R}^d$ and, for all $i,j\in \{1,\dots, L\}$, $\Gamma_{i,j}=\partial^*\Omega_i\cap\partial^*\Omega_j\cap \Omega$ is the interface between phases $i,j$ (with $\partial^*\Omega_i$ the reduced boundary of $\Omega_i$, see~\cite{AmbrosioFuscoPallara2000})
 and  $\sigma_{i,j}$ is the surface tension along this interface. To ensure the lower semicontinuity of the $L$-phase perimeter, see \cite{morgan1997lowersemicontinuity,maggi2012sets,caraballo2009triangle}, we assume that the surface tensions are positive, i.e. $\sigma_{i,j} >0$, and satisfy the triangle inequality  
$$ \sigma_{i,j} + \sigma_{j,k} \geq \sigma_{i,k}\quad \text{ for any } i,j,k,$$

The evolution of the liquid-vapor-solid system toward equilibrium can be approximated by a multiphasic surface diffusion flow.
This motion  can be viewed as the $H^{-1}$ gradient flow of the energy \eqref{per_multip_phase_intro} 
which ensures its decay while maintaining locally the volume of each phase. In particular, the normal velocity $V_{ij}$ 
at the interface  $\Gamma_{ij}$ reads as 
\begin{equation*}
\frac{1}{\nu_{ij}}V_{ij} = \sigma_{ij} \Delta_{\Gamma_{ij}(t)}H_{ij}(t),
\end{equation*}
where $H_{ij}(t)$ denotes the scalar mean curvature on $\Gamma_{ij}(t)$,  
$\Delta_{\Gamma_{ij}(t)}$ is the Laplace-Beltrami operator on the surface, and $\nu_{ij} > 0$ is the surface mobility coefficient. The above expression is the classical form of the velocity in this context, but it can obviously be rewritten to incorporate the degenerate no-motion case where $\nu_{ij}=0$:
\begin{equation*}
V_{ij} = \nu_{ij}\sigma_{ij} \Delta_{\Gamma_{ij}(t)}H_{ij}(t).
\end{equation*}

{The dewetting phenomenon is closely related, in several typical situations, to the above model. Recall that dewetting is the process by which a continuous film forced to cover a surface retracts and breaks down into islands or droplets. This phenomenon occurs not only for liquid films, but also for solid films when heated (solid-state dewetting), see the references in~\cite{Dziwnik_2017}. In general, capillary effects have a prominent role~\cite{Srolovitz1986CapillaryII}}.

 {Classical liquid/solid dewetting} involves $L = 3$ phases: the liquid phase $\Omega_L$, the solid phase $\Omega_S$ and the vapor phase $\Omega_V$.
Moreover, as the surface tension coefficients $( \sigma_{LV}, \sigma_{SV} ,  \sigma_{SL})$ satisfy the
 triangle inequality, they form an additive set of coefficients, i.e. there exist three nonnegative 
 coefficients $\sigma_{L}, \sigma_{S}, \sigma_{V}$ such that
 $$ \sigma_{LV} = \sigma_{L} +  \sigma_{V}, \quad  \sigma_{SV} = \sigma_{S} +  \sigma_{V},\quad \text{ and } \quad {\sigma_{SL} = \sigma_{S} +  \sigma_{L}}.$$
 These coefficients are given by:
 $${\sigma_{S}=\frac{\sigma_{SV}+\sigma_{SL}-\sigma_{LV}}2,\quad \sigma_{V}=\frac{\sigma_{SV}+\sigma_{LV}-\sigma_{SL}}2,\quad \sigma_{L}=\frac{\sigma_{LV}+\sigma_{SL}-\sigma_{SV}}2}$$
 The surface mobilities can be set to 
 $$(\nu_{LV},\nu_{SV},\nu_{SL}) = (1,0^+,0^+),$$
 in order to fix the solid phase.   This set of coefficients is  harmonically additive in the sense 
 that, with the convention $\frac 1{0^+}=+\infty$, there exist three non negative coefficients $\nu_S$, $\nu_L $,  and $\nu_V$ such that
   $$ \nu_{LV}^{-1} = \nu_{L}^{-1} +  \nu_{V}^{-1}, \quad \nu_{SV}^{-1} = \nu_{S}^{-1} +  \nu_{V}^{-1} \text{ and }  \quad 
  \nu_{SL}^{-1} = \nu_{S}^{-1} +  \nu_{L}^{-1}.$$
   Indeed, we can just consider $\nu_S = 0^+$ and $\nu_L = \nu_V = 2$.

~\\
Having in mind the application to dewetting, we assume in the rest of the paper that:
\begin{itemize}
\item  the surface tensions are additive,
i.e. there exist coefficients $\sigma_i \geq 0$,  $i\in\{1,\dots,L\}$, such that
$\sigma_{ij} = \sigma_i + \sigma_j$, $\forall i,j\in\{1,\dots,L\}$;
\item the mobility coefficients  are harmonically
additive, i.e. there exist nonnegative coefficients $\nu_i$ satisfying $\nu_{ij}^{-1} = \nu_{i}^{-1} +  \nu_{j}^{-1}$ (with the convention that $\frac 1{0^+}=+\infty$).
\end{itemize}

With such assumptions, it is easy to reformulate the expression of the $L$-phase perimeter in the more convenient following form
\begin{equation} \label{Intro:MultiphasePerimeter}
P(\Omega_1,\dots,\Omega_L)= \sum_{i=1}^L \sigma_{i} P(\Omega_{i})  = { \sum_{i=1}^L \sigma_{i} \mathcal{H}^{d-1}(\partial^* \Omega_i)} ,
\end{equation}
where $\partial^* \Omega_i$ denotes the reduced boundary of $\Omega_i$. In this form, the  $L$-phase perimeter can be approximated in the sense of $\Gamma$-convergence by a sum of scalar Cahn-Hilliard energies~\cite{Modica1977} defined for every smooth ${\bf u}=(u_1,\dots, u_L)$ by
$$
 P_{\varepsilon}({\bf u}) =  
\begin{cases}
\displaystyle\sum_{k=1}^{L} \sigma_k  \int_\Omega \left( \frac{\varepsilon}{2}  |\nabla u_k|^2 + \frac{1}{\epsilon} {W}(u_k) \right) dx  & \text{ if }   \sum_{k=1}^{L} u_k = 1, \\
 +\infty & \text{otherwise,}
\end{cases}
$$
In this definition each $u_i$ represents a smooth approximation of the characteristic function $\one{\Omega_i}$, ${W}(s) = \frac{s^2(1-s)^2}{2}$
is a double-well potential, and the parameter $\varepsilon$ characterizes the width 
of the diffuse interface, i.e. how much each $\nabla u_idx$ is concentrated around the Hausdorff measure supported on the reduced boundary of $\Omega_i$.

\begin{rem}
{
For $L=3$ phases, surface tensions satisfying the triangle inequality are always additive. It is not always the case as soon as $L\geq 4$, but under weak conditions on the surface tensions and a suitable choice of multi-well potentials, it is again possible~\cite{bretin2017new} to design phase field approximations that are convenient for numerical simulation but whose analysis is rather difficult.}
\end{rem}

Since the multiphase surface diffusion flow is the $H^{-1}$-gradient flow of~\eqref{Intro:MultiphasePerimeter},  a natural idea to approximate it is to consider the $H^{-1}$-gradient flow of 
$P_{\varepsilon}$  which yields the Cahn-Hilliard system
\begin{equation}\label{global-CH}
\left\lbrace
\begin{aligned}
& \varepsilon^2 \partial_t u_k =  \nu_k \Delta \left( \sigma_k \mu_k +  \lambda \right) \\
& \mu_k = W'(u_k) - \varepsilon^2 \Delta u_k,
\end{aligned}
\right.
\end{equation}
where $\lambda$ is the Lagrange multiplier associated with the partition constraint $\sum u_k = 1$. 
Here, we follow \cite{bretin2018multiphase} to handle the set $\{\nu_{i,j}\}$ of mobilities  and we use explicitly  
 its harmonically additive decomposition. 
 
\begin{rem}
{
There are  physical situations, e.g. total wetting, where triangle inequality fails. In the solid-liquid-vapor configuration, the total wetting corresponding to a liquid film (no contact line between vapor and solid) is associated with  $ \sigma_{SV}\geq\sigma_{SL} +\sigma_{LV}$, and the total wetting due to a gaz film (no contact line between liquid and solid) is associated with $\sigma_{SL}\geq \sigma_{SV}+\sigma_{LV}$.  The strict inequalities are not consistent with the lower semicontinuity required for energy minimization, yet the limit cases $ \sigma_{SV}=\sigma_{SL} +\sigma_{LV}$ or $\sigma_{SL}= \sigma_{SV}+\sigma_{LV}$ have to be considered. In the first case, $\sigma_L=\sigma_S=0$, and in the latter case $\sigma_V=\sigma_S=0$. Our approach can actually easily handle the situations where $\sigma_S=0$, and even  $\sigma_S\leq 0$, using a simple coupling with a null mobility $\nu_S=0$. It allows in practice to preserve a well-posed phase field system.}
\end{rem}

The asymptotic expansion of the phase field system~\eqref{global-CH} is delicate and, to the best of our
knowledge, no rigorous analysis of its convergence has been made so far. The main obstacle
to overcome is the non local nature of the system, which is particularly significant in the multiphase case. \\

In \cite{bretin2020approximation}, we reviewed various two-phase Cahn-Hilliard systems and we
proposed a new one. It basically involves degenerate mobilities that vanish in 
pure phase regions, therefore localize the system and allow to prove asymptotic results. 
In the next paragraph we sum up the properties and choices of parameters in the biphasic case (see  \cite{bretin2020approximation} for details) before the extension to the multiphase case.
 
Recall that \cite{pego1989front,alikakos1994convergence} proved that the classical Cahn-Hilliard equation
\begin{equation*}
\left\lbrace
\begin{aligned}
& \varepsilon^2 \partial_t u = \Delta \mu, \\
& \mu = W'(u) - \varepsilon^2 \Delta u,
\end{aligned}
\right.
\end{equation*}
does not converge to surface diffusion flow but rather to the Hele-Shaw model which is non local.
Cahn et al. \cite{cahn1996cahn} introduced a new system involving a concentration-dependent mobility $M$. 
It is often referred  as a degenerate mobility in the sense that no motion occurs in the pure state regions. 
The model proposed by Cahn et al. is  the following equation that we will refer to as \MCH:
\begin{equation*}
\left\lbrace
\begin{aligned}
& \varepsilon^2 \partial_t u = \div\left(M(u)\mu\right),\\
& \mu = W'(u) - \varepsilon^2 \Delta u.
\end{aligned}
\right.
\end{equation*}
A formal convergence to the correct motion is shown in  \cite{cahn1996cahn}. However, the particular model studied by Cahn et al invoves a logarithmic potential $W$, which raises numerical issues. Instead, the potential commonly chosen in the literature and the one that we will use for the remainder of this paper is the smooth potential
$$
W(s)= \frac{1}{2}s^2(1-s)^2.
$$
The choice of the mobility $M$ has been discussed theoretically in \cite{gugenberger2008comparison,lee2015degenerate,MR3466205}. It is proven by a formal asymptotic method that the choice $M(u)=u(1-u)$ does not lead to the correct velocity as an additional bulk diffusion term appears. These conclusions have been  corroborated numerically in \cite{dai2012motion,dai2014coarsening,MR3457961} where undesired coarsening effects are observed. Actually a quartic mobility $M(u)=u^2(1-u)^2$ is necessary to recover the correct velocity. These conclusions have been extended to the anisotropic case in~\cite{Dziwnik_2017}. From now on, we fix
$$
M(s) = s^2(1-s)^2.
$$
While the \MCH ~model has the correct sharp interface limit and produces satisfactory numerical results, 
it has a well identified drawback: in the asymptotic, the leading error term is of order $1$ and becomes 
relevant when reaching the pure states $0$ or $1$, causing oscillations and an imprecise profile for the solution. 
The problem is twofold. Firstly, the solution does not remain within the physical range of $[0,1]$,
which means that in the multiphase context, some phases might be negative in some areas and larger than $1$ in others (in other words, what is called {\it positivity property} by some authors is not fulfilled). Secondly, as illustrated in \cite{bretin2020approximation}, the approximation being of order $1$ only, it induces numerical volume losses despite the natural volume preservative nature of the Cahn-Hilliard equation.\\

{
The failure to meet the positivity condition and its numerical illustration with Fourier spectral approaches in~\cite{MR3457961,bretin2020approximation} seem to be in contradiction with the analytical result of \cite{MCH_Garcke_Elliott} regarding the existence of weak solutions contained in $[0,1]$ of a Cahn-Hilliard model with degenerate mobility. We believe it is not a contradiction, it rather illustrates the nonuniqueness of the solutions to the Cahn-Hilliard equations due to possible bifurcations when the pure states $0$ or $1$ are attained. The nonuniqueness is for example illustrated in  \cite{MR1742748} where a finite elements approach is introduced that captures confined solutions only, but varying the mesh or the time step gives various solutions with different behavior. Recent numerical analyses of these solutions have shown their singular behavior~\cite{MR3466205,MR4284406}, and although the finite elements approach used certainly allows to reduce oscillations, the solutions' singularities and asymptotic behavior impact significantly the pointwise and integral approximation errors. In particular, a very fine resolution is necessary to capture accurately the solutions and the correct domain of values if the constraint to take values in $[0,1]$ is not forced with a potential. In addition, the accuracy of these numerical solutions to approximate the continuous solution is no better than $O(\varepsilon)$ and, in particular, they remain at distance greater than a positive multiple of $\varepsilon$ from the pure states $0$ and $1$.\\}

Regarding numerical accuracy, the authors of \cite{ratz2006surface} managed to improve it by introducing 
another degeneracy in the model. It has been successfully adapted in various applications, 
see for example \cite{albani2016dynamics,naffouti2017complex,salvalaglio2017morphological,salvalaglio2015faceting}. 
However, the proposed model does not derive from an energy, it is thus more difficult to prove rigorously theoretical properties and to extend to complex multiphase applications. A variational adaptation has been proposed in
\cite{salvalaglio2019doubly} where the second degeneracy is injected in the energy.
But because it relies on modifying the energy, the approach is hard to extend to complex multiphase or anisotropic 
applications. \\

In \cite{bretin2020approximation}, we proposed a different approach where an additional mobility $N$
is incorporated in the metric of the gradient flow instead of plugging it into the energy, and thus the geometry of the evolution problem.
The so-called \NMNCH ~model proposed in~\cite{bretin2020approximation} reads as
\begin{equation*}
\begin{cases}
   \varepsilon^2 \partial_t u &= N(u) \div\left(M(u) \nabla (N(u)\mu) \right) \\
 \mu &= W'(u) - \varepsilon^2 \Delta u,\\
  \end{cases}
\end{equation*}
The presence of two supplementary terms $N(u)$ is needed to ensure the variational nature of the model.
Using formal asymptotic expansion, we showed in~\cite{bretin2020approximation} that a good choice for $N$ is 
$$
N(s) = \frac{1}{\sqrt{M(s)}}= \frac{1}{s(1-s)},
$$
Indeed, it allows to nullify the error term of order $1$ in the solution, making the \NMNCH ~model of order~$2$.
The profile obtained for the solution $u$ is very accurate and the volume conservation is ensured up to an error of order $2$, to be compared with the order $1$ for \MCH. As observed in~\cite{bretin2020approximation}, another choice for $N$ which avoids issues with the pure phases $s=0,1$ without changing the conclusions of the asymptotic expansion is $N(s)=\frac 1{\sqrt{s^2(1-s)^2+\gamma\varepsilon^2}}$, with $\gamma>0$.

\begin{rem}
 {
 Regarding the positivity property, as mentioned above, we actually believe that the key point is not so much the numerical solution being confined in $[0,1]$, but rather the quality of the numerical approximation. As illustrated in~\cite{bretin2020approximation}, the numerical solution obtained with the \MCH ~model for approximating the evolution by surface diffusion of a thin structure is well contained in $[0,1]$, but the approximation error in $O(\varepsilon)$ prevents it from representing correctly the continuous solution. In contrast, with the  \NMNCH ~model in $O(\varepsilon^2)$ that we propose and using the same type of numerical method, we obtain a much more realistic numerical solution. It may not be valued in $[0,1]$, but it is a more accurate approximation.}
\end{rem}

In this paper, we extend the \MCH ~and \NMNCH ~models to the case of $L$ phases.
From the modeling viewpoint, this amounts to integrating in the model the influence of surface tensions $\sigma_{ij}$ and phase 
mobilities $\nu_{ij}$. To this end, we adapt to the Cahn-Hilliard system the work of \cite{bretin2018multiphase} done for the Allen-Cahn 
system. In particular, we propose to analyze the two following phase field models, where in both cases $\lambda$ is the Lagrangian multiplier which encodes the partition constraint $\sum_{k=1}^L u_k = 1$:

\begin{itemize}
 \item The \MCH ~multiphase field model defined for $k\in\{1,\dots,L\}$ by:
 \begin{equation}
\label{cahnHilliardMulti}
\left\lbrace
\begin{aligned}
& \varepsilon^2\partial_t u_k = \nu_k\div\left(M(u_k)\nabla( \sigma_k \mu_k + \lambda) \right), \\
& \mu_k = W'(u_k) - \varepsilon^2 \Delta u_k, \\
\end{aligned}
\right.
\end{equation}
with mobility $M(s) = 2 W(s)$.
 \item The \NMNCH ~multiphase field model defined for $k\in\{1,\dots,L\}$ by:
 \begin{equation}
\label{Eq:NMNmulti}
\left\lbrace
\begin{aligned}
& \varepsilon^2 \partial_t u_k = \nu_kN(u_k)\div\left(M(u_k)\nabla( \sigma_k N(u_k)\mu_k + \lambda) \right), \\
& \mu_k = W'(u_k) - \varepsilon^2 \Delta u_k.\\
\end{aligned}
\right.
\end{equation}
\end{itemize}
with mobilities $M(s) = 2 W(s)$ and $N(s) = 1/\sqrt{M(s)}$. This model is well defined whenever $u\not=0,1$, which is the case near the interface $\{u=\frac 1 2\}$. To give sense to the model in the whole domain, it can be rewritten in two different ways:
\begin{itemize}
\item either by transferring $N$ to the left-hand side to obtain the alternative model
\begin{equation}
\label{Eq:OuterSystemNMN-init}
\text{(\NMNCH ~reformulation I)}\qquad \left\lbrace
\begin{aligned}
& \varepsilon^2 g(u_k) \partial_t u_k = \nu_k\div\left(M(u_k)\nabla( \sigma_k \mu_k + \lambda) \right),\\
& g(u_k)\mu_k = W'(u_k) - \varepsilon^2 \Delta u_k.\\
\end{aligned}
\right.
\end{equation} 
where $g(u_k) = \sqrt{M(u_k)}$ is always well-defined. Such a reformulation (strictly equivalent where $M(u_k)$ does not vanish) will be used for the asymptotic expansion.
\item or by modifying the definitions of $M, N$ to prevent them from vanishing while preserving the conclusions of the asymptotic expansion. This is the case with the following model:
 $$
\text{(\NMNCH ~reformulation II)}\qquad  \begin{cases}
  \partial_t u_k &= \nu_k \tilde N(u_k) \div\left(\tilde M(u_k)\nabla \tilde N(u_k) ( \sigma_k \mu_k + \lambda) \right) \\
   \mu_k & = \frac{W'(u_k)}{\varepsilon^2} - \Delta u_k \\ 
 \end{cases}
 $$
 \end{itemize}
 where the mobilities $\tilde M$ and $\tilde N$ are defined by $\tilde M(s) = 2W(s) + \gamma \varepsilon^2$ and $\tilde N(s) = \frac{1}{\sqrt{\tilde M(s)}}$, with $\gamma>0$. Obviously, $\tilde M$ never vanishes and $\tilde N$ is well-defined everywhere. We will explain in the first lines of Section~\ref{sec:asymptotic-NMNCH} why this reformulation has the same asymptotic properties as the original model~\eqref{Eq:NMNmulti}. The \NMNCH ~reformulation II model will be used for numerical approximation (because numerical errors require a choice for $M$ that prevents cancellations).

\subsection{Outline of the paper}

We first proceed to a formal asymptotic analysis of the \MCH ~and \NMNCH ~multiphase models.
In particular, we show that the limit law of each model is indeed the multiphase surface diffusion flow with the advantage that \NMNCH ~guarantees an approximation error of order $2$ in $\varepsilon$. In a second section, devoted to numerical approximation, we first introduce a numerical scheme suitable for both models. This scheme is based on a Fourier-spectral
convex-concave semi implicit approach in the spirit of~\cite{MR1676409,bretin2020approximation}.
We provide numerical experiments which illustrate the stability of our scheme and the asymptotic properties
of both phase field models.  In the last section, we consider the special case of  the wetting / dewetting phenomenon for which
we derive a simplified, yet equivalent model using the liquid phase only. We illustrate this model with 3D numerical experiments using either 
smooth or rough surfaces, and choosing various set of parameters to get different Young angle conditions.

%
%
%
%
%
%
\section{Formal matched asymptotic expansions}
\label{Sec:AsymptoticSection}
In this section, we give a formal proof of Propositions~\ref{result:modelM} and~\ref{result:modelNMN} below
using the method of matched asymptotic expansions. These results involve the so-called \emph{optimal profile} $q$ associated with the potential $W$ and defined by the equation $q'(z) = - \sqrt{2W(q(z))}$ with a suitable constraint on $q(0)$. In the case where $W(s)=\frac{1}{2}s^2(1-s)^2$ and $q(0)=\frac 1 2$, one gets 
$$q(z) = \frac{1-\tanh\left(\frac{z}{2}\right)}{2}.$$
The following constants are also used in both propositions:
$$ c_W = \int_\mathbb{R} \left( q'(z) \right)^2 dz, \quad   c_M = \int_\mathbb{R} M(q(z))dz \quad \text{ and } \quad  c_N = \int_\mathbb{R} \frac{q'(z)}{N(q(z))} dz.$$
Remark that with our particular choices for $N$ and $q$, one has $c_N=-c_W$.

\begin{prop}
\label{result:modelM}
For $i,j\in\{1,\dots,L\}$ with $i\not=j$, let $\Omega_i^\varepsilon=\{x,\;u_i(x)\geq \frac 1 2\}$ and 
$$\Gamma_{ij}^{\varepsilon}=\partial \Omega_i^\varepsilon\cap\{x,\;u_j\geq u_k,\;k\in\{1,\dots,L\}\setminus\{i\}\}.$$
The solution ${\bf u}^{\varepsilon}$ to the \MCH ~model defined for $k\in\{1,\dots,L\}$ by
\begin{equation*}
\left\lbrace
\begin{aligned}
& \varepsilon^2 \partial_t u_k = \nu_k \div\left(M(u_k)\nabla( \sigma_k  \mu_k + \lambda) \right), \\
& \mu_k = W'(u_k) - \varepsilon^2 \Delta u_k, \\
\end{aligned}
\right.
\end{equation*}
 satisfies (formally) near the interface $\Gamma_{ij}^{\varepsilon}$ the following asymptotic expansions:
\begin{equation*}
\left\lbrace
\begin{aligned}
& u^{\varepsilon}_i = q\left( \frac{\operatorname{dist}(x,\Omega_i^{\varepsilon})}{\varepsilon} \right) + \mathcal{O}(\varepsilon), \\
& u^{\varepsilon}_j = 1- q\left( \frac{\operatorname{dist}(x,\Omega_i^{\varepsilon})}{\varepsilon} \right) + \mathcal{O}(\varepsilon), \\
& u^{\varepsilon}_k = \mathcal{O}(\varepsilon). \\
\end{aligned}
\right.
\end{equation*}
where $\dist(\cdot,\Omega_i^{\varepsilon})$ denotes the signed distance function to $\Omega_i^{\varepsilon}$.\\
Moreover, the normal velocity $V^{\varepsilon}_{ij}$ at the interface $\Gamma_{ij}^{\varepsilon}$ satisfies (formally):
\begin{equation*}
\frac{1}{\nu_{ij}} V_{ij}^{\varepsilon} = \sigma_{ij} c_M c_W \Delta_{\Gamma_{ij}^{\varepsilon}} H_{ij} + \mathcal{O}(\varepsilon).
\end{equation*}
\end{prop}
\begin{prop}
\label{result:modelNMN}
With the notations of Proposition~\ref{result:modelM}, the solution ${\bf u}^{\varepsilon}$ to the \NMNCH ~model defined for $k\in\{1,\dots,L\}$ by
\begin{equation*}
\left\lbrace
\begin{aligned}
& \varepsilon^2 \partial_t u_k = \nu_kN(u_k)\div\left(M(u_k)\nabla( \sigma_k N(u_k)\mu_k + \lambda) \right), \\
& \mu_k = W'(u_k) - \varepsilon^2 \Delta u_k, \\
\end{aligned}
\right.
\end{equation*}
 satisfies (formally) near the interface $\Gamma_{ij}^{\varepsilon}$ the following asymptotic expansions:
\begin{equation*}
\left\lbrace
\begin{aligned}
& u_i^{\varepsilon} = q\left( \frac{\operatorname{dist}(x,\Omega_i^{\varepsilon})}{\varepsilon} \right) + \mathcal{O}(\varepsilon^2), \\
& u_j^{\varepsilon} = 1- q\left( \frac{\operatorname{dist}(x,\Omega_i^{\varepsilon})}{\varepsilon} \right) + \mathcal{O}(\varepsilon^2), \\
& u_k^{\varepsilon} = \mathcal{O}(\varepsilon^2). \\
\end{aligned}
\right.
\end{equation*}
Moreover, the normal velocity $V_{ij}^{\varepsilon}$ at the interface $\Gamma_{ij}^{\varepsilon}$ satisfies (formally):
\begin{equation*}
\frac{1}{\nu_{ij}} V_{ij}^{\varepsilon} = \sigma_{ij} \frac{c_W c_M}{(c_N)^2} \Delta_{\Gamma_{ij}^{\varepsilon}} H_{ij} + \mathcal{O}(\varepsilon).
\end{equation*}
\end{prop}

\begin{figure}[htbp]
\centering
	\includegraphics[width=12cm]{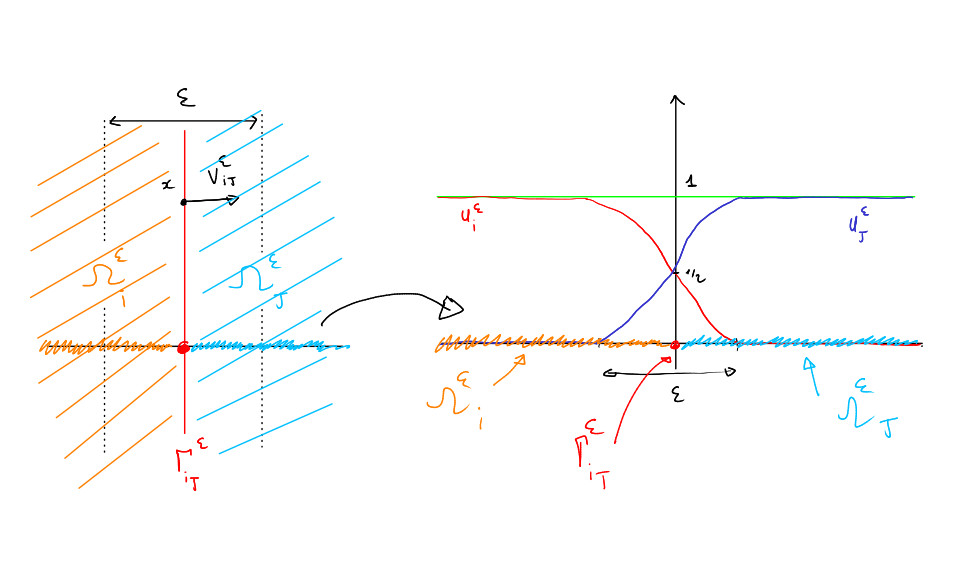} 
\caption{Illustration of the notations used and schematics profiles of the solutions $u_i^{\varepsilon}$ and $u_j^{\varepsilon}$ near an interface $\Gamma^\varepsilon_{ij}$.}
\label{fig_illustration}
\end{figure}

The results stated in both propositions are  {illustrated in Figure~\ref{fig_illustration}}. To prove these results, we first recall the necessary tools following the notations of 
\cite{alfaro2013convergence,chen2011mass,MR3383330} and the presentation in~\cite{bretin2020approximation}.  
Then we proceed to the asymptotic expansion for the \MCH ~model. The proof is shown in 
dimension 2 only  for the sake of simplicity of notations and readability, but it can be readily extended to higher dimensions. We end up with the \NMNCH ~model, which we have to 
rewrite to avoid indeterminate forms. Some calculations remain the same as for the biphasic case presented in~\cite{bretin2020approximation},  {but the presence of the Lagragian multiplier which does not exist for two phases only requires additional calculations}.

\subsection{Formal asymptotic analysis toolbox} 

In this multiphase context, we study the behavior of the system in two regions:
near the interface $\Gamma:=\Gamma^\varepsilon_{ij}$ separating two given phases $i\not=j$, and far from it. 
We denote $u_k$ the solution for an arbitrary phase $k$. Whether $k$ can designate $i$ or $j$ in an equation
will be clear from the context.

To derive the method we require that the interface $\Gamma=\Gamma^\varepsilon_{ij}$ remains smooth enough 
so that there exist $\delta >0$ and  a neighborhood  
$$\mathcal{N}=\mathcal{N}_{ij}^{\delta}(\Gamma) = \lbrace x\in\Omega / |d(x,t)| < 3\delta \rbrace,$$
in which the signed distance function $d:=d_{ij}$ to $\Gamma$ is well-defined. $\mathcal{N}$ is 
called the \emph{inner region} near the interface and its complement the \emph{outer region}.\\

\noindent {\it Outer variables:}~\\
Far from the interface, we consider the \emph{outer functions} $(u_k,\mu_k)$ depending 
on the standard \emph{outer variable} $x$. The systems remain the same, e.g., for \MCH: 
\begin{equation}
\label{Eq:OuterSystemMCH}
\left\lbrace
\begin{aligned}
& \varepsilon^2 \partial_t u_k = \nu_k\div(M(u_k)\nabla(\sigma_k\mu_k+ \lambda)), \\
& \mu_k = -\varepsilon^2 \Delta u_k + W'(u_k).
\end{aligned}
\right.
\end{equation}

\begin{figure}[htbp]
\centering
	\includegraphics[width=10cm]{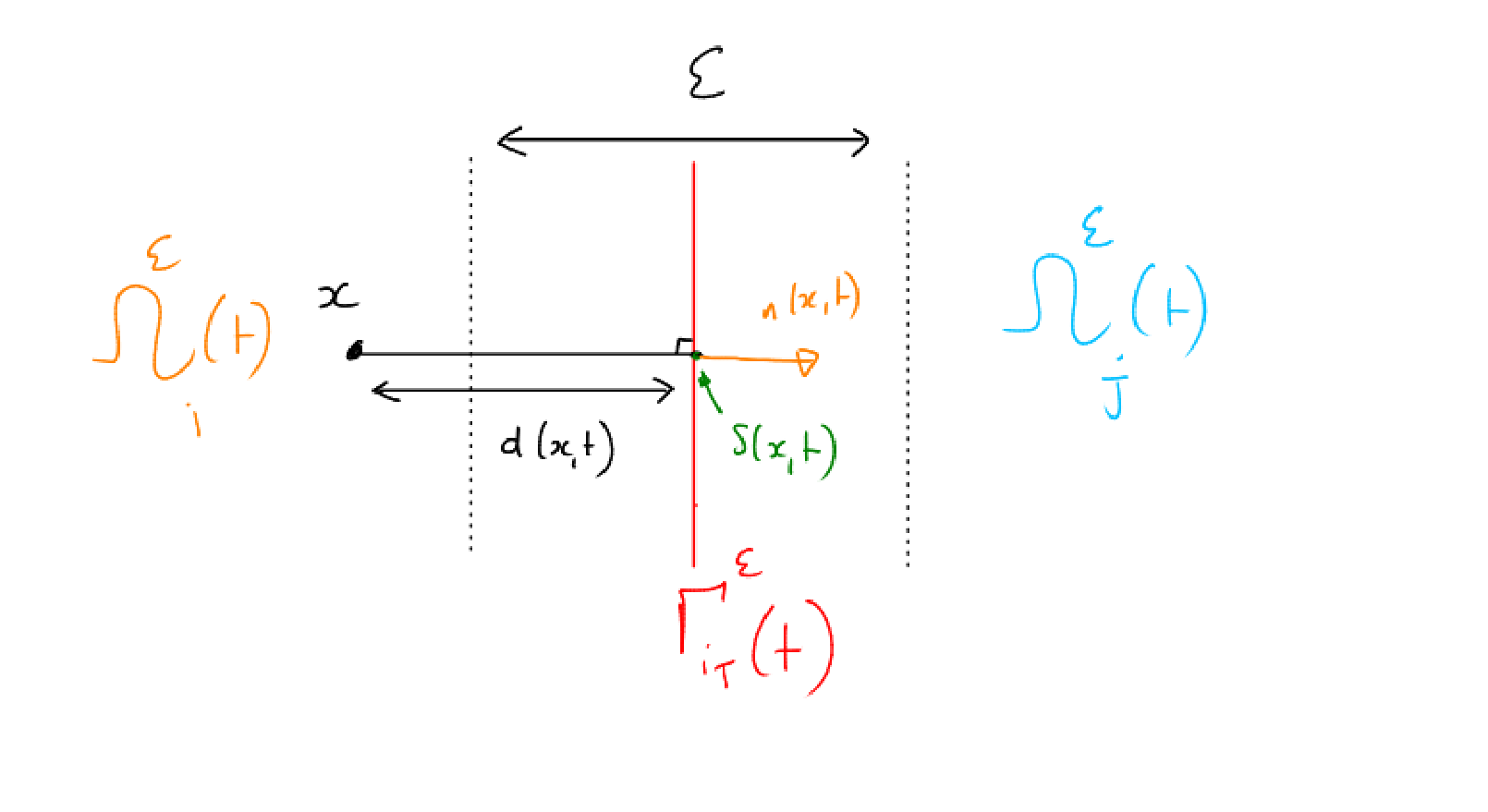} 
\caption{The inner variables are $s$ and $z=\frac d\varepsilon$.}
\label{fig_inner_variable}
\end{figure}

\noindent {\it Inner variables:}~\\
Inside $\mathcal{N}$ we consider the \emph{inner variables} $(z,s)$ associated with the original variables $(x,t)$ in the following way (see Figure~\ref{fig_inner_variable}): $z=\frac{d(x,t)}{\varepsilon}$ is a variable along the normal direction $n(x,t)$ to the interface $\Gamma$ and $s=S(x,t)$ is associated with a parameterization $X_0(s,t)$ of $\Gamma$. We define the \emph{inner functions} $U_k,\bigmu_k$ depending on $(z,s)$ as follows:
\begin{equation*}
\left\lbrace
\begin{aligned}
& U_k(z,s,t):= U_k\left( \frac{d(x,t)}{\varepsilon}, S(x,t), t \right) = u_k(x,t) \\
& \bigmu_k(z,s,t): = \bigmu_k\left( \frac{d(x,t)}{\varepsilon}, S(x,t), t \right) = \mu_k(x,t) \\
\end{aligned}
\right.
\end{equation*}
In order to express the derivatives of $U_k$, we first need to calculate the gradient and the Laplacian 
of $d$ and $S$.  The properties of $d$ are well-known, see for instance~\cite{AmbrosioGeom}:
\begin{equation*}
\left\lbrace
\begin{aligned}
& \nabla d(x,t) = n(x,t), \\
& \Delta d(x,t) = \sum_{l=1}^{d-1} \frac{\kappa_l(\pi(x))}{1+\kappa_l(\pi(x))d(x,t)} =\frac{H}{1+\varepsilon z H} \text{ in dimension 2}. 
\end{aligned}
\right.
\end{equation*}
where $\pi$ is the orthogonal projection onto $\Gamma$ and $\kappa_1,\dots, \kappa_{d-1}$ are the principal curvatures on $\Gamma$.

Given a point $X_0(s,t)$ on $\Gamma$, let 
\begin{equation*}
X(z,s,t) = X_0(s,t) + \epsilon z n(s,t).
\end{equation*}
whose orthogonal projection onto $\Gamma$ is $X_0(s,t)$. 
The equation connecting the variable $s$ and the function $S$ is:
\begin{equation*}
s = S(X_0(s,t) + \varepsilon z n(s,t),t).
\end{equation*}
Deriving this equation with respect to $z$ leads to
 \begin{equation*}
0  = \varepsilon n \cdot \nabla S = \varepsilon \nabla d \cdot \nabla S{,}
\end{equation*}
which implies there is no cross derivative term. The derivation of the same equation with 
respect to $s$ gives
\begin{equation*}
1  = \left( \partial_s X_0 + \varepsilon z H \partial_s n \right) \cdot \nabla S = (1+\varepsilon z H) \tau \cdot \nabla S.\\
\end{equation*}
Since $\nabla S$ is orthogonal to $n$, therefore collinear with the tangent $\tau$, we have that:
\begin{equation*}
\nabla S = \frac{1}{1+\varepsilon z H} \tau
\end{equation*}
Taking the divergence, we find $\Delta S$:
\begin{equation*}
\left. \begin{aligned}
\Delta S & = \div\left( \frac{\tau}{1+\varepsilon z H}  \right) = \nabla \left( \frac{1}{1+\varepsilon z H} \right) \cdot \tau + \frac{1}{1+\varepsilon z H} \div(\tau), \\
& = \frac{1}{1+\varepsilon z H} \partial_s \left( \frac{1}{1+\varepsilon z H}\right) + \frac{1}{1+\varepsilon z H} \tau \cdot  \partial_s \tau, \\
&= - \frac{\varepsilon z \partial_s H}{(1+\varepsilon z H)^3}.
\end{aligned} \right.
\end{equation*}
To express the connection between the derivatives of $U_k,\bigmu_k$ and $u_k,\mu_k$,
we come back to the definition of the inner functions:
\begin{equation*}
u_k(x,t) = U_k\left(\frac{d(x,t)}{\varepsilon}, S(x,t),t) \right).
\end{equation*}
Successive derivations with respect to $x$ give the following equations
\begin{equation*}
\left\lbrace \begin{aligned}
& \nabla u_k = \nabla d \frac{1}{\varepsilon}\partial_z U_k + \nabla S \partial_s U_k, \\
& \Delta u_k = \Delta d \frac{1}{\varepsilon}\partial_z U_k + \frac{1}{\varepsilon^2} \partial_{zz} U_k + \Delta S \partial_s U_k + |\nabla S|^2 \partial_{ss} U_k, \\
& \div\left( M(u_k) \nabla (N(u_k)\mu_k)\right) = \frac{1}{\varepsilon^2} \partial_{z} {\left( M_k \partial_z(N_k\bigmu_k)\right)} + \frac{M_k}{\varepsilon} \Delta d \partial_z (N_k\bigmu_k),  \\
& \qquad \qquad \qquad \qquad \qquad \quad + |\nabla S|^2 \partial_s \left( M_k \partial_s (N_k\bigmu_k) \right) + \Delta S M_k \partial_s (N_k\bigmu_k).
\end{aligned}\right.
\end{equation*}
The \emph{inner system} of the \MCH ~model near the interface $\Gamma$ finally reads as:
\begin{equation*}
\left\lbrace
\begin{aligned}
& \frac{\varepsilon^2}{\nu_k}\left(\partial_t U_k + \partial_t S \partial_s U_k  \right) - \frac{\varepsilon}{\nu_k} V_{ij} \partial_z U_k = \frac{1}{\varepsilon^2} \partial_z \left( M(U_k)\partial_z \left( \sigma_k \bigmu_k + \Lambda \right)\right), \\
& \qquad \qquad \qquad \qquad \qquad \qquad \qquad \qquad + \frac{M(U_k)}{\varepsilon} \partial_z \left( \sigma_k \bigmu_k + \Lambda \right) \Delta d_{ij} + T_1(s),\\
& \bigmu_k = W'(U_k) - \partial_{zz} U_k {-} \varepsilon \Delta d_{ij} \partial_z U_k {-}  \varepsilon^2 T_2(s), \\
& \Delta d_{ij} = \frac{H_{ij}}{1+\varepsilon z H_{ij}} = H_{ij} - \varepsilon z H_{ij}^2 + \mathcal{O}(\varepsilon^2),\\
& T_1(s) = \frac{\partial_s(M(U_k)\partial_s (\sigma_k\bigmu_k+\Lambda))}{(1+\varepsilon z H_{ij})^2} - \frac{M(U_k) \varepsilon z \partial_s H_{ij}}{(1+\varepsilon z H_{ij})^3}  \partial_s (\sigma_k\bigmu_k+\Lambda),\\
& T_2(s) = \frac{1}{(1+\varepsilon z H)^2}\partial_{ss} U_k - \frac{\varepsilon z \partial_s H_{ij}}{(1+\varepsilon z H_{ij})^3} \partial_s U_k.\\
\end{aligned}
\right.
\end{equation*}
Note that the terms in $T_1$ and $T_2$ are high-order tangential terms that play a role only at the fourth order in the asymptotic expansion.\\

\noindent {\it Independence in $z$ of the normal velocity $V_{ij}$:}~\\
The normal velocity $V_{ij}(s,t)$ of the interface is defined by:
\begin{equation*}
\left. \begin{aligned}
V_{ij}(s,t) & = \partial_t X_{0}(s,t) \cdot n(s,t). \\
\end{aligned}\right.
\end{equation*}
In the neighborhood $\mathcal{N}$, we have the following property 
(which is a direct consequence of the definition of the signed distance function):
\begin{equation*}
d(X_{0}(s,t) + \varepsilon z n(s,t),t) = \varepsilon z.
\end{equation*}
Deriving this with respect to $t$ yields:
\begin{equation*}
V_{ij}(s,t) = \partial_t X_{0}(s,t) \cdot \nabla d(X_{0}(s,t)+ \varepsilon z n(s,t),t) = - \partial_t d \left( X(z,s,t),t\right).
\end{equation*}
Thus, the function $\partial_t d(x,t)$ is independent of $z$ and can be extended in the whole neighborhood
by choosing
\begin{equation*}
V_{ij}(X_{0}(s,t) + \varepsilon z n,t) := -\partial_t d(X_{0}(s,t) +\varepsilon z n,t) = V_{ij}(s,t).
\end{equation*}
This property of independence is crucial to be able to extract the velocity from integrals in $z$ in the following derivations.\\ 

\noindent {\it Taylor expansions:}~\\
We assume the following Taylor expansions for our functions:
\begin{equation*}
\left.
\begin{aligned}
& u_k(x,t) = u_k^{(0)}(x,t)+\varepsilon u_k^{(1)}(x,t)+ \varepsilon^2 u_k^{(2)}(x,t) + \cdots \\
& U_k(z,s,t) = U_k^{(0)}(z,s,t) + \varepsilon U_k^{(1)}(z,s,t) + \varepsilon^2 U_k^{(2)}(z,s,t) + \cdots \\
& \mu_k(x,t) = \mu_k^{(0)}(x,t) + \varepsilon \mu_k^{(1)}(x,t) + \varepsilon^2 \mu_k^{(2)}(x,t) + \cdots \\
& \bigmu_k(z,s,t) = \bigmu_k^{(0)}(z,s,t) + \varepsilon\bigmu_k^{(1)}(z,s,t) + \varepsilon^2\bigmu_k^{(2)}(z,s,t) + \cdots \\
& \lambda(x,t) = \varepsilon \lambda^{(1)}(x,t)+\varepsilon^2 \lambda^{(2)}(x,t) + \cdots \\
& \Lambda(z,s,t) = \varepsilon\Lambda^{(1)}(z,s,t) + \varepsilon^2 \Lambda^{(2)}(z,s,t) + \cdots
\end{aligned}
\right.
\end{equation*}

Since the numbering of the phase is present as a subscript, we indicate the order in the Taylor expansion
as a superscript in brackets. We can then compose these expansions with a regular function $F$:
\begin{equation*}
\left.
\begin{aligned}
F(U_k) = & \ F(U_k^{(0)}) + \varepsilon F'(U_k^{(0)})U_k^{(1)}  + \varepsilon^2 \left[ F'(U_k^{(0)})U_k^{(2)}+\frac{F''(U_k^{(0)})}{2} (U_k^{(1)})^2\right] \\
	   & \quad \quad  \ + \varepsilon^3\left[F'(U_k^{(0)})U_k^{(3)} + F''(U_k^{(0)}) U_k^{(1)}U_k^{(2)} + \frac{F'''(U_k^{(0)})}{6}(U_k^{(1)})^3 \right] + \cdots \\
\end{aligned}
\right.
\end{equation*}
To simplify the notations within the asymptotic expansion, we adopt the following notations for $M(u_k)$
\begin{equation*}
M(u_k) = m_k^{(0)} + \varepsilon m_k^{(1)} + \varepsilon^2 m_k^{(2)} + \cdots,
\end{equation*}
where 
\begin{equation*}
\left\lbrace \begin{aligned}
& m_k^{(0)} = M(u_k^{(0)}),\\
& m_k^{(1)} = M'(u_k^{(0)})u_k^{(1)}, \\
& m_k^{(2)} = M'(u_k^{(0)})u_k^{(2)} + \frac{M''(u_k^{(0)})}{2}(u_k^{(1)})^2.\\
\end{aligned} \right.
\end{equation*}
We adopt the same convention for any generic \emph{outer function} $F(u_k)$ or \emph{inner function} $F(U_k)$:
\begin{equation*}
\left.
\begin{aligned}
& F(u_k) = f_k^{(0)} + \varepsilon f_k^{(1)} + \varepsilon^2 f_k^{(2)} + \varepsilon^3 f_k^{(3)} +\cdots \\
& F(U_k) = F_k^{(0)} + \varepsilon F_k^{(1)} + \varepsilon^2 F_k^{(2)} + \varepsilon^3 F_k^{(3)} +\cdots \end{aligned}
\right.
\end{equation*}

\noindent {\it Flux matching condition between inner and outer equations:}~\\
Instead of using the matching conditions directly between the first equations of the inner and outer systems, 
it is more convenient to do the matching for the flux 
$$j_k = M(u_k) \nabla (\sigma_k\mu_k+\lambda).$$
$j_k$  has the following Taylor expansion for the \MCH ~model 
\begin{equation}
\label{Eq:MCHOuterFlux}
\left.
\begin{aligned}
j_k & = \left[ m_k^{(0)}\nabla(\sigma_k \mu_k^{(0)}  )\right] + \varepsilon\left[m_k^{(1)} \nabla (\sigma_k\mu_k^{(0)}) +m_k^{(0)}\nabla (\sigma_k\mu_k^{(1)})\lambda^{(1)}) \right] \\
		 & \quad + \varepsilon^2 \left[ m_k^{(2)}\nabla(\sigma_k\mu_k^{(0)}) + m_k^{(1)}\nabla(\sigma_k\mu_k^{(1)} +\lambda^{(1)})+ m_k^{(0)} \nabla (\sigma_k\mu_k^{(2)}+\lambda^{(2)}) \right] + \mathcal{O}(\varepsilon^3).
\end{aligned}
\right.
\end{equation}
In inner coordinates, we only need to express the normal part 
$$J_{k,n}:= J_k \cdot n = \frac{ M(U_k)}{\varepsilon} \partial_z (\sigma_k\bigmu_k+\Lambda),$$
as the tangential terms are of higher order. The normal part expands as
\begin{equation}
\label{Eq:MCHInnerFlux}
\left.
\begin{aligned}
J_{k,n} & = \ \ \frac{1}{\varepsilon} \left[ M_k^{(0)} \partial_z(\sigma_k\bigmu_k^{(0)}) \right],\\
 	& \quad \quad + \ \ \left[M_k^{(1)} \partial_z (\sigma_k\bigmu_k^{(0)}) + M_k^{(0)} \partial_z (\sigma_k\bigmu_k^{(1)}+\Lambda^{(1)}) \right],\\
	&\quad \quad + \varepsilon \left[M_k^{(2)}\partial_z(\sigma_k\bigmu_k^{(0)}) + M_k^{(1)}\partial_z (\sigma_k\bigmu_k^{(1)}+\Lambda^{(1)})+ M_k^{(0)} \partial_z (\sigma_k\bigmu_k^{(2)}+\Lambda^{(2)}) \right], \\
	& \quad \quad + \varepsilon^2 \left[M_k^{(3)}\partial_z (\sigma_k\bigmu_k^{(0)}) + M_k^{(2)}\partial_z (\sigma_k\bigmu_k^{(1)}+\Lambda^{(1)}), \right. \\
	& \qquad \qquad \qquad \left.  + M_k^{(1)} \partial_z (\sigma_k\bigmu_k^{(2)}+\Lambda^{(2)}) +M_k^{(0)}\partial_z (\sigma_k\bigmu_k^{(3)}+\Lambda^{(3)}) \right] + \mathcal{O}(\varepsilon^3).
\end{aligned}
\right.
\end{equation}
The flux matching conditions allow to match the limit as $z\rightarrow\pm\infty$ of terms of \eqref{Eq:MCHOuterFlux} 
with the corresponding order terms of \eqref{Eq:MCHInnerFlux}. 

We can now investigate order by order the behavior of the  \MCH ~model. We have to study up to
the fourth order term where the leading order of the velocity will appear in the first equation of
\eqref{Eq:OuterSystemMCH}.  After that, we adapt the argument to the \NMNCH ~model, 
where a reformulation of the problem will be necessary
 to avoid indeterminate forms in the asymptotic expansion.

\subsection{Formal matched asymptotic expansion for the multiphasic \MCH ~model} 
We first establish Proposition~\ref{result:modelM} regarding the properties of the \MCH ~model. 
We recall that we study the behavior of the different terms of the system near the interface $\Gamma^\varepsilon_{ij}$ 
separating phases $i$ and $j$. We assume the following matching conditions for the two phases:
\begin{equation*}
\left. \begin{aligned}
& \lim_{z\rightarrow + \infty} U_i^{(0)} = 0, \quad \lim_{z\rightarrow -\infty} U_i^{(0)} = 1, \\
& \lim_{z\rightarrow + \infty} U_j^{(0)} = 1, \quad \lim_{z\rightarrow -\infty} U_j^{(0)} = 0. \\
\end{aligned} \right.
\end{equation*}
For the other phases, we require the following matching conditions
\begin{equation*}
\lim_{z\rightarrow \pm \infty} U_k^{(0)} = 0, \quad \lim_{z\rightarrow \pm\infty} U_k^{(1)} = 0.
\end{equation*}
{ The definition of $\Gamma_{ij}^{\varepsilon}$ also shows that 
$$  U_i^{(0)}(0,s,t) = \frac 1 2,\;  U_i^{(1)}(0,s,t) = 0,\; \text{ and }  U_i^{(2)}(0,s,t) = 0.$$
}
~\\
{\it First order:}\\
At order \ordres{\varepsilon^{-2}}{1} the inner system reads as
\begin{equation*}
\left\lbrace
\begin{aligned}
& 0 = \partial_z \left( M_k^{(0)} \partial_z(\sigma_k\bigmu_k^{(0)})\right),\\
& \bigmu_k^{(0)} = W'(U_k^{(0)}) - \partial_{zz} U_k^{(0)}.\\
\end{aligned}
\right.
\end{equation*}
From the first equation of the system, we deduce that $M_k^{(0)}\partial_z(\sigma_k\bigmu_k^{(0)})$ is constant. 
The matching conditions for the outer flux \eqref{Eq:MCHOuterFlux} and the inner flux \eqref{Eq:MCHInnerFlux} at
order $\varepsilon^{-1}$ impose this constant to be zero. Then there is a constant $A_k^{(0)}$ such that
\begin{equation*}
\sigma_k\bigmu_k^{(0)} = A_k^{(0)}.
\end{equation*} 
Collecting all this information, we obtain that:
\begin{equation*}
\partial_{zz} U_k^{(0)} - W'(U_k^{(0)}) = {-} \frac{A_k^{(0)}}{\sigma_k} \quad \forall k \in \{1,\dots,L \}.
\end{equation*}
{Finally, the  matching conditions on $U_k^{(0)}$ and the initial condition $U_i^{(0)}(0,s,t) = \frac 1 2$ imply that}
\begin{equation*}
\left\lbrace
\begin{aligned}
& U_i^{(0)} = q(z), \\
& U_j^{(0)} = q(-z) = 1-q(z),\\
& U_k^{(0)} = 0 \quad \forall k \in \{1,\dots,L\}\setminus\lbrace i,j \rbrace, \\
& \bigmu_k^{(0)} =0 \quad \forall k\in \{1,\dots,L\}, \\
\end{aligned}
\right.
\end{equation*}
where  $q$ is the optimal phase field profile.\\

\noindent {\it Second order:}\\
At order \ordres{\varepsilon^{-1}}{\varepsilon} the outer system reads as
\begin{equation}
\label{Eq:Inner2ClassicalMultiphase}
\left\lbrace
\begin{aligned}
& 0 = \partial_z \left(M_k^{(0)} \partial_z (\sigma_k\bigmu_k^{(1)} + \Lambda^{(1)} )\right),\\
& \bigmu_k^{(1)} = W''(U_k^{(0)})U_k^{(1)} - \partial_{zz} U_k^{(1)} - H_{ij} \partial_z U_k^{(0)}. \\
\end{aligned}
\right.
\end{equation}
It follows that there exists a function $B_k^{(1)}$ constant in $z$ such that 
\begin{equation*}
M_k^{(0)}\partial_z(\sigma_k\bigmu_k^{(1)} +\Lambda^{(1)}) = B_k^{(1)}.
\end{equation*}
By the matching condition between the outer flux \eqref{Eq:MCHOuterFlux} and the inner flux \eqref{Eq:MCHInnerFlux}
at order $1$, it holds that
\begin{equation*}
\lim_{z\rightarrow \pm \infty} M_k^{(0)} \partial_z (\sigma_k\bigmu_k^{(1)}+\Lambda^{(1)}) = 0.
\end{equation*}
We deduce  that $B_k^{(1)}=0$ and that there exists a function $A_k^{(1)}$ constant in $z$ such that
\begin{equation*}
\sigma_k\bigmu_k^{(1)} +\Lambda^{(1)} = A_k^{(1)}.
\end{equation*}
Subtracting the case  $k=i$ from the case $k=j$ gives
\begin{equation*}
\sigma_j\bigmu_j^{(1)} - \sigma_i\bigmu_i^{(1)} = A_j^{(1)} -A_i^{(1)}.
\end{equation*}
The term $A_{ij}^{(1)}:=  A_j^{(1)} -A_i^{(1)}$ can be determined using the second equation of 
\eqref{Eq:Inner2ClassicalMultiphase}. Indeed, recall that 
\begin{equation*}
\left\lbrace
\begin{aligned}
&  \sigma_i \bigmu_i^{(1)} = \sigma_i  W''(U_i^{(0)})U_i^{(1)} - \sigma_i  \partial_{zz} U_i^{(1)} - \sigma_i  H_{ij} q', \\
& \sigma_j \bigmu_j^{(1)} = \sigma_j W''(U_{{j}}^{(0)})U_j^{(1)} - \sigma_j  \partial_{zz} U_j^{(1)} +\sigma_j H_{ij} q'.
\end{aligned}
\right.
\end{equation*}
We multiply both equations by $q'$ and integrate the difference. We can eliminate the terms in $U^{(1)}$ through integration by parts:
\begin{equation*}
\left. 
\begin{aligned}
{ \int_\mathbb{R} W''(U_i^{(0)}) U_i^{(1)} \partial_z U_0^{(0)} - \partial_{zz} U_i^{(1)} \partial_z U_0^{(0)}  dz} & =
\int_\mathbb{R} { \partial_z(W'(U_i^{(0)})) U_i^{(1)} - \partial_{zz} U_i^{(1)} \partial_z U_i^{(0)}} dz,\\
& =\left[ W'(U_i^{(0)}) U_i^{(1)} - \partial_z U_i^{(1)}  \partial_z U_i^{(0)}  \right]_{-\infty}^{+\infty},\\
		& \qquad - \int_\mathbb{R} \partial_z U_i^{(1)} \left(\underbrace{ W'(U_i^{(0)}) - \partial_{zz}U_i^{(0)}}_{=0} \right) dz, \\
		& = 0.
\end{aligned}
\right.
\end{equation*}
It follows that
\begin{equation}
\label{Eq:MultiphaseAIJM}
\left.
\begin{aligned}
A_{ij}^{(1)} & = -\int_\mathbb{R} (\sigma_j\bigmu_j^{(1)}-\sigma_i\bigmu_i^{(1)}) q'\,dz = -(\sigma_j+\sigma_i) H_{ij}\int_\mathbb{R}(q')^2\,dz = -\sigma_{ij} c_W H_{ij}.\\
\end{aligned}
\right.
\end{equation}
On the other hand, summing the second equation of system \eqref{Eq:Inner2ClassicalMultiphase} for the phases $i$ and $j$ gives
\begin{equation*}
\bigmu_i^{(1)}+\bigmu_j^{(1)} = W''(q) \left[U_i^{(1)}+U_j^{(1)} \right] - \partial_{zz}\left[U_i^{(1)} + U_j^{(1)} \right].
\end{equation*}
Multiplying by $q'$ and integrating by parts gives:
\begin{equation*}
\int_\mathbb{R} \left( \bigmu_i^{(1)}+\bigmu_j^{(1)} \right) q' dz= 0.
\end{equation*}
{
Thus, for all $s$,  there exists a profile $\zeta_s$ satisfying $\int q' \zeta_s dz = 0 $ and such that
$\bigmu_i^{(1)}+\bigmu_j^{(1)} = \zeta_s(z)$.}
Combining this with equation \eqref{Eq:MultiphaseAIJM}, we obtain 
\begin{equation*}
\left\lbrace
\begin{aligned}
& \bigmu_i^{(1)} = c_W H_{ij} + \zeta_s(z) \frac{\sigma_i}{\sigma_{ij}}, \\
& \bigmu_j^{(1)} = {-}c_W H_{ij} + \zeta_s(z)  \frac{\sigma_j}{\sigma_{ij}}. \\
\end{aligned}
\right.
\end{equation*}
and then,
$$
\begin{cases}
W''(U_i^{(0)})U_i^{(1)}-\partial_{zz}U_i^{(1)} = H_{ij}(c_W+q')+\zeta_s(z)\frac{\sigma_j}{\sigma_{ij}},\\
W''(U_j^{(0)})U_j^{(1)}-\partial_{zz}U_j^{(1)} = -H_{ij}(c_W+q')+\zeta_s(z)\frac{\sigma_i}{\sigma_{ij}}.
\end{cases}
$$
which leads to 

\begin{equation*}
\begin{cases}
U_i^{(1)}(z,s)=H_{ij}\eta(z)+c(s)\frac{\sigma_j}{\sigma_{ij}} \omega_s(z),\\
U_j^{(1)}(z,s)=-H_{ij}\eta(z)+c(s)\frac{\sigma_i}{\sigma_{ij}} \omega_s(z).
\end{cases}
\end{equation*}

{
Here $\eta$ and $\omega_s$ are two profiles defined as the solutions to { $W''(q)y-y''=q'+c_W$ and $W''(q)y-y''=\zeta_s$}, respectively, with appropriate initial conditions.  Note that such profiles exist because $ \int(q'+c_W)q' dz = 0 $ and 
$ \int \zeta_s q' dz = 0 $.}\\

We deduce from the above system that if $H_{ij}\neq 0$, then $U_i^{(1)}$ and $U_j^{(1)}$ cannot vanish both together which yields the important conclusion that
the leading order error term for the solution of the system is no better than $\varepsilon$: the \MCH ~model is {\it always of order $1$} when the mean curvature is non zero.
It justifies the interest of the \NMNCH ~model which is of second order. \\

\noindent {\it Third order:}\\
At order \ordres{1}{\varepsilon^2} the inner system reads as 
\begin{equation}
\label{Eq:Inner3Multiphase}
\left\lbrace
\begin{aligned}
& 0 = \partial_z \left( M_k^{(0)} \partial_z (\sigma_k\bigmu_k^{(2)}+ \Lambda^{(2)}) \right), \\
& \bigmu_k^{(2)} = \frac{W'''(U_k^{(0)})}{2} (U_k^{(1)})^2 + W''(U_k^{(0)}) U_k^{(2)} - \partial_{zz} U_k^{(2)} {-} H_{ij} \partial_z U_k^{(1)} {+} z H_{ij}^{(2)} \partial_z U_k^{(0)}.\\
\end{aligned}
\right.
\end{equation} 
In the first equation, we used the results from the first two orders and left out the term that vanishes. 
From the first equality of \eqref{Eq:Inner3Multiphase}, we find that:
\begin{equation*}
M_k^{(0)} \partial_z(\sigma_k\bigmu_k^{(2)} + \Lambda^{(2)}) = B_k^{(2)}.
\end{equation*}
The matching conditions between the flux \eqref{Eq:MCHOuterFlux} and \eqref{Eq:MCHInnerFlux} 
at order $\varepsilon$ yield (by removing all the null terms):
\begin{equation*}
B_k^{(2)} = \lim_{z\rightarrow \pm \infty} M_k^{(0)} \partial_z ( \sigma_k\bigmu_k^{(2)} + \Lambda^{(2)} ) = 0.
\end{equation*}
This means that the term $\sigma_k\bigmu_k^{(2)}+\Lambda^{(2)}$ is constant in $z$ and will not intervene in the flux term 
of order $\varepsilon^2$. \\

\noindent {\it Fourth order:}\\
Collecting the previous results, the first equation of the inner system at order $\varepsilon$ for the phases $i$ and $j$ simplifies to 
\begin{equation*}
\left\lbrace
\begin{aligned}
-& \frac{1}{\nu_i} V_{ij} q' = \partial_z \left[ M_i^{(0)} \partial_z (\sigma_i\bigmu_i^{(3)} + \Lambda^{(3)})\right] + \partial_s\left[ M_i^{(0)} \partial_s(\sigma_i\bigmu_i^{(1)} +\Lambda^{(1)}) \right],\\
& \frac{1}{\nu_j} V_{ij} q' = \partial_z \left[ M_j^{(0)} \partial_z (\sigma_j\bigmu_j^{(3)}+\Lambda^{(3)}) \right] + \partial_s\left[ M_j^{(0)} \partial_s(\sigma_j\bigmu_j^{(1)} +\Lambda^{(1)}) \right], \\
\end{aligned}
\right.
\end{equation*}
We subtract the two equations, and integrate. We divide the computation in three:
\begin{itemize}
\item The left hand side gives:
\begin{equation*}
\left(\frac{1}{\nu_i} + \frac{1}{\nu_j} \right) V_{ij} = \frac{1}{\nu_{ij}} V_{ij}.
\end{equation*}
\item Collecting the result from the previous paragraphs, we find the following matching between the outer flux \eqref{Eq:MCHOuterFlux} 
and the inner flux \eqref{Eq:MCHInnerFlux} at order $\varepsilon^2$:
\begin{equation*}
\lim_{z\rightarrow \pm \infty}  M_i^{(0)} \partial_z (\sigma_i\bigmu_i^{(3)} + \Lambda^{(3)}) = m_i^{(1)} \nabla(\sigma_i\mu_i^{(1)}+\lambda^{(1)}).
\end{equation*}
Because $M'(0) = M'(1) = 0$, the limit term is zero and then
\begin{equation*}
\int_\mathbb{R} \partial_z \left( M_i^{(0)} \partial_z (\sigma_i\bigmu_i^{(3)} + \Lambda^{(3)}) \right)dz= 0.
\end{equation*}
The corresponding term for the $j$-th phase is treated similarly.
\item Using \eqref{Eq:MultiphaseAIJM}, the second term of the right hand side is (noting that $M_i^{(0)} = M_j^{(0)}$):
\begin{equation*}
\int_\mathbb{R} M_i^{(0)} \partial_{ss}(\sigma_i\bigmu_i^{(1)} - \sigma_j\bigmu_j^{(1)})dz = \sigma_{ij} \left( \int_\mathbb{R} M(q(z)) dz \right) c_W \partial_{ss}H_{ij}.
\end{equation*}
\end{itemize}
Finally, we obtain that 
\begin{equation*}
\frac{1}{\nu_{ij}}V_{ij} = \sigma_{ij} c_W c_M \partial_{ss}H_{ij}.
\end{equation*}
which concludes the proof of Proposition~\ref{result:modelM}.

\subsection{Formal matched asymptotic expansion for the multiphase \NMNCH ~model}\label{sec:asymptotic-NMNCH}
We now give a proof of Proposition \ref{result:modelNMN} concerning the properties of \NMNCH. 
{We assume the same matching conditions as for the \MCH ~model}. \\

\noindent {\it Reformulation of the model:}\\
 It is more convenient to rewrite the  \NMNCH ~model by transferring $N(u_k)$ to the left hand side of the system, which yields the \NMNCH ~reformulation I model we already mentioned:
\begin{equation}
\label{Eq:OuterSystemNMN}
\left\lbrace
\begin{aligned}
& \varepsilon^2 g(u_k) \partial_t u_k = \nu_k\div\left(M(u_k)\nabla( \sigma_k \mu_k + \lambda) \right),\\
& g(u_k)\mu_k = W'(u_k) - \varepsilon^2 \Delta u_k.\\
\end{aligned}
\right.
\end{equation} 
where $g(u_k) = \sqrt{M(u_k)}= \frac{1}{N(u_k)}$ when $N(u_k)$ is well-defined, and, as before, $\lambda$ is the Lagrangian multiplier which encodes the partition constraint {$\sum_{k=1}^{L} u_k=0$}. As already said, the advantage of such a formulation is that $g(u_k)$ is always well defined even if $u_k = 0$, 
which is not the case for $N(u_k)$. Note also that the definition of $\mu_k$ has been changed but we keep the same notation for simplicity.\\

Remark that similar calculations as those shown below can be done for the \NMNCH ~reformulation II model which is used for numerical approximation, and the same conclusions of Proposition~\ref{result:modelNMN} hold. Actually, using the additional term $\gamma \varepsilon^2$ in the definition of $M$  does not change the asymptotic results for at least the first four orders of interest. Indeed, this term appears to be associated with $\bigmu_0$ (see below) which is zero, and $\bigmu_1$ whose derivative in $z$ vanishes. \\

\medskip
\noindent The inner system for \NMNCH ~reformulation I now reads (for simplicity, we drop the expression "reformulation I" in the calculations below):
\begin{equation*}
\left\lbrace
\begin{aligned}
& \frac{\varepsilon^2 G_k}{\nu_k} \partial_t U_k +\frac{\varepsilon^2 G_k}{\nu_k} \partial_tS \partial_s U_k - \frac{\varepsilon G_k}{\nu_k} V_{ij} \partial_zU_k = \frac{1}{\varepsilon^2} \partial_{z}\left(M_k \partial_z(\sigma_k\bigmu_k+\Lambda)\right) \\
&\qquad \qquad \qquad \qquad \qquad \qquad \qquad \qquad \qquad \qquad + \frac{M_k}{\varepsilon} \Delta d_{ij} \partial_z(\sigma_k\bigmu_k+\Lambda) + T_1(s), \\
& G_k\bigmu_k = W'(U_k) - \partial_{zz} U_k - \varepsilon \Delta d_{ij} \partial_z U_k + \varepsilon^2 T_2(s), \\
& \Delta d_{ij} = \frac{H_{ij}}{1+\varepsilon z H_{ij}} = H_{ij} - \varepsilon z H_{ij}^2 + \mathcal{O}(\varepsilon^2),\\
& T_1(s) = \frac{1}{(1+\varepsilon z H_{ij})^2} \partial_{s}(M_k\partial_s(\sigma_k\bigmu_k+\Lambda)) - \frac{\varepsilon z M_k\partial_s H_{ij}}{(1+\varepsilon z H_{ij})^3} \partial_s(\sigma_k\bigmu_k+ \Lambda), \\
& T_2(s) = \frac{1}{(1+\varepsilon z H_{ij})^2} \partial_{ss}U_k-\frac{\varepsilon z \partial_s H_{ij}}{(1+\varepsilon z H_{ij})^3} \partial_s U_k. \\
\end{aligned}
\right.
\end{equation*}
Because $N(u_k)$ is now on the left hand side of the system in the form of $G_k$, the flux term $j = M(u_k)\nabla\left(\sigma_k\mu_k+\lambda_k\right)$ is the same as the one for \MCH. The flux matching condition is then also equal to the one given by \eqref{Eq:MCHOuterFlux} and \eqref{Eq:MCHInnerFlux}. \\

\noindent {\it First order:}\\
At order \ordres{\varepsilon^{-2}}{1} the inner system reads:
\begin{equation*}
\left\lbrace
\begin{aligned}
& 0 = \partial_z \left( M_k^{(0)} \partial_z(\sigma_k\bigmu_k^{(0)})\right),\\
& G_k^{(0)}\bigmu_k^{(0)} = W'(U_k^{(0)}) - \partial_{zz} U_k^{(0)}. \\
\end{aligned}
\right.
\end{equation*}
From the first equation of the system, we deduce that $M_k^{(0)}\partial_z(\sigma_k\bigmu_k^{(0)})$ is constant.
The matching conditions on the outer \eqref{Eq:MCHOuterFlux} and inner fluxes \eqref{Eq:MCHInnerFlux} at order $\varepsilon^{-1}$ 
impose this constant to be zero. Then there exists a constant $A_k^{(0)}$ in $z$ such that
\begin{equation*}
\sigma_k\bigmu_k^{(0)} = A_k^{(0)}.
\end{equation*} 
Collecting all this information, we have
\begin{equation*}
\partial_{zz} U_k^{(0)} - W'(U_k^{(0)}) = \frac{G\left(U_k^{(0)}\right)A_k^{(0)}}{\sigma_k}, \quad \forall k\in \{1,\ldots,L\}.
\end{equation*}
Then, using the matching conditions and {the initial conditions $U_i^{(0)}(0,s,t)  = \frac{1}{2}$} leads to
\begin{equation*}
\left\lbrace
\begin{aligned}
& U_i^{(0)} = q(z), \\
& \bigmu_i^{(0)} = 0, \\
& U_j^{(0)} = q(-z) = 1-q(z),\\
& \bigmu_j^{(0)} = 0, \\
& U_k^{(0)} = 0 \quad \forall k \in \{1,\ldots,L\}\setminus\lbrace i,j \rbrace . \\
\end{aligned}
\right.
\end{equation*}
Notice that $\bigmu_k^{(0)}$ is a constant in $z$ that can be nonzero.\\

\noindent {\it Second order:}\\
 At order \ordres{\varepsilon^{-1}}{\varepsilon} the inner system reads
\begin{equation}
\label{Eq:Inner2NMNmulti}
\left\lbrace
\begin{aligned}
& 0 = \partial_z \left(M_k^{(0)} \partial_z (\sigma_k\bigmu_k^{(1)} + \Lambda^{(1)} )\right),\\
& G_k^{(0)} \bigmu_k^{(1)} + G_k^{(1)} \bigmu_k^{(0)} = W''(U_k^{(0)})U_k^{(1)} - \partial_{zz} U_k^{(1)} - H_{ij} \partial_z U_k^{(0)}. \\
\end{aligned}
\right.
\end{equation}
For $k\neq i,j$, using the fact that $U_k^{(0)}=0$, the second equation can be rewritten as
\begin{equation*}
0 = \left(1-\bigmu_k^{(0)}\right)U_k^{(1)}-\partial_{zz}U_k^{(1)}.
\end{equation*}
As the matching conditions show that $\lim_{z \pm \infty} U_k^{(1)} = 0$, it follows that $U_k^{(1)}=0$.\\

Now turning to the $i$-th phase (resp. $j$-th), there exists a function $B_i^{(1)}$ constant in $z$ such that
\begin{equation*}
M_i^{(0)}\partial_z(\sigma_i\bigmu_i^{(1)} +\Lambda^{(1)}) = B_i^{(1)}.
\end{equation*}
From the matching condition between outer \eqref{Eq:MCHOuterFlux} and inner flux \eqref{Eq:MCHInnerFlux} at order $1$, we deduce that
\begin{equation*}
\lim_{z\rightarrow \pm \infty} M_i^{(0)} \partial_z (\sigma_i\bigmu_i^{(1)}+\Lambda^{(1)}) = 0,
\end{equation*}
and $B_i^{(1)}=0$ (resp $B_j^{(1)}=0$). Then there exist functions $A_i^{(1)}, A_j^{(1)}$ constant in $z$ such that
\begin{equation*}
\left\lbrace
\begin{aligned}
& \sigma_i\bigmu_i^{(1)} +\Lambda^{(1)} = A_i^{(1)},\\
& \sigma_j\bigmu_j^{(1)} + \Lambda^{(1)} = A_j^{(1)}. \\
\end{aligned}
\right.
\end{equation*}
Subtracting the $i$-th term to the $j$-th term leads to
\begin{equation*}
\sigma_j\bigmu_j^{(1)} - \sigma_i\bigmu_i^{(1)} = A_j^{(1)} -A_i^{(1)} := A_{ij}^{(1)}.
\end{equation*}
Moreover, recall that $U_k^{(1)}=0$ for $k\neq i,j$, which implies that $U_i^{(1)}=-U_j^{(1)}$ as  $\sum_{k=1}^N U_k^{(1)} = 0$.
Using the symmetry properties $g(1-s) = g(s)$ and $W''(1-s)=W''(s)$, it follows that:
\begin{equation*}
\left.
\begin{aligned}
 g(U_j^{(0)})\bigmu_j^{(1)} & = W''(U_j^{(0)}) U_j^{(1)} - \partial_{zz} U_j^{(1)}- H_{ij}\partial_zU_j^{(0)},\\
 & = -W''(U_i^{(0)})U_i^{(1)}+ \partial_{zz}U_i^{(1)} + H_{ij}\partial_z U_i^{(0)},\\
 & = -g(U_i^{(0)})\bigmu_i^{(1)}, \\
 & = -g(U_j^{(0)})\bigmu_i^{(1)}.\\
\end{aligned}
\right.
\end{equation*}
thus $ g(U_j^{(0)})\left( \bigmu_i^{(1)}+\bigmu_j^{(1)}\right)=0$. Finally, as g$(q)\neq 0$, $\bigmu_i^{(1)}$ 
and $\bigmu_j^{(1)}$ are necessarily constant in $z$ and
\begin{equation*}
\bigmu_j^{(1)} = -\bigmu_i^{(1)}.
\end{equation*}
It shows that we can  express $A_{ij}^{(1)}$ as
\begin{equation*}
A_{ij}^{(1)} =\left(\sigma_j+\sigma_i\right)\bigmu_j^{(1)} = -\sigma_{ij}\bigmu_i^{(1)}.
\end{equation*}
Now, multiplying the second equation of \eqref{Eq:Inner2NMNmulti} for phase $i$ 
\begin{equation*}
g(U_i^{(0)}) \bigmu_i^{(1)} = W''(U_i^{(0)})U_i^{(1)} - \partial_{zz} U_i^{(1)} - H_{ij} q', 
\end{equation*}
by the profile $q'$ and integrating over $\R$ shows that 
\begin{equation}
\label{Eq:MultiphaseAIJNMN}
  \bigmu_i^{(1)} = - \frac{c_W}{c_N} H_{ij} \quad \text{ and } \quad  A_{ij}^{(1)} =  \frac{c_W}{c_N} \sigma_{ij} H_{ij}.
\end{equation}
Indeed, on the one hand we have
\begin{equation*}
\left.
\begin{aligned}
\int_\mathbb{R} (\partial_z(W'(U_i^{(0)})) U_i^{(1)} - \partial_{zz} U_i^{(1)} \partial_z U_i^{(0)})dz & = \left[ W'(U_i^{(0)}) U_i^{(1)} - \partial_z U_i^{(1)}  \partial_z U_i^{(0)}  \right]_{-\infty}^{+\infty}, \\
		& \qquad - \int_\mathbb{R} \partial_z U_i^{(1)} \left(\underbrace{ W'(U_i^{(0)}) - \partial_{zz}U_i^{(0)}}_{=0}  \right)dz, \\
		& = 0,
\end{aligned}
\right.
\end{equation*} 
and on the other hand
\begin{equation*}
\int_\mathbb{R} g(q) \bigmu_i^{(1)} q'dz = \bigmu_i^{(1)}\int_\mathbb{R} \frac{q'}{N(q)}dz = c_N \bigmu_i^{(1)}, \quad\text{ and } \int_\mathbb{R} H_{ij} (q')^2dz = c_W H_{ij}. 
\end{equation*}

Finally, it follows that
$$ \partial_{zz} U_i^{(1)} - W''(q) U_i^{(1)} = - H_{i,j} (q' - \frac{c_W}{c_N}g(q)).$$
Now, recall that the choice $M(s) = \sqrt{2 W(s)}$ and  $N(s) = 1/\sqrt{M(s)}  = 1/\sqrt{2 W(s)}$ for the mobilities
implies that  
$$g(q) = 1/N(q) =  \sqrt{2 W(q)} = - q' \;\;\text{ and }\;\;  c_N = -c_W.$$
This is the key point to understand why in this case the term $U_i^{(1)}$ is null as a solution of
$$ \partial_{zz} U_i^{(1)} - W''(q) U_i^{(1)} = 0.$$
The same argument gives $U_j^{(1)}=0$.  In summary, we have $U_k^{(1)} = U_i^{(1)} = U_j^{(1)}= 0$.
It means that the leading error order term in the solutions $U_i$ and $U_j$ is of magnitude $\varepsilon^2$ while
the other phases are absent.\\

\noindent {\it Third order:}\\
At order \ordres{1}{\varepsilon^2}, using both previous orders, the inner system simplifies to 
\begin{equation*}
\left\lbrace
\begin{aligned}
& 0 = \partial_z \left( M_k^{(0)} \partial_z (\sigma_k\bigmu_k^{(2)}+ \Lambda^{(2)}) \right),\\
& G_k^{(0)} \bigmu_k^{(2)} = \frac{W'''(U_k^{(0)})}{2} (U_k^{(1)})^2 + W''(U_k^{(0)}) U_k^{(2)} - \partial_{zz} U_k^{(2)} + H_{ij} \partial_z U_k^{(1)} - z H_{ij}^2 \partial_z U_k^{(2)}. \\
\end{aligned}
\right.
\end{equation*}
From the first equality, we find that
\begin{equation*}
M_k^{(0)} \partial_z(\sigma_k\bigmu_k^2 + \Lambda^2) = B_2.
\end{equation*}
The matching conditions at order $\varepsilon$ for the fluxes given by \eqref{Eq:MCHOuterFlux} and \eqref{Eq:MCHInnerFlux}  
yield also (by removing all the null terms): 
\begin{equation*}
B_2 = \lim_{z\rightarrow \pm \infty} M_k^{(0)} \partial_z ( \sigma_k\bigmu_k^{(2)} + \Lambda^{(2)} ) = 0.
\end{equation*}
Therefore, $\sigma_k\bigmu_k^{(2)}+\Lambda^{(2)}$ is constant in $z$ and will not intervene in the flux term of order $\varepsilon^2$.

\begin{rem}
It is possible to show that  $U_i^{(2)}$ and $U_j^{(2)}$ are of the form  
$$U_j^{(2)} = -U_i^{(2)} = H_{ij}^2 \zeta(z),$$
where $\zeta$ is the profile defined by 
\begin{equation*}
\left\lbrace
\begin{aligned}
& y''(z) - W''(q) y(z) = z q' \\
& y(0) = 0 \\
\end{aligned}
\right.
\end{equation*}
and decreasing to zero at infinity.
\end{rem}
\noindent{\it Fourth order:}\\
Eliminating all vanishing terms, the first equations for the $i$-th and $j$-th phases of the inner system read
\begin{equation}
\label{Eq:NMNFourthorder}
\left.
\begin{aligned}
& - \frac{1}{\nu_i} V_{ij} g(U_i^{(0)}) \partial_z U_i^{(0)} = \partial_z\left[M_i^{(0)}\partial_z\left(\sigma_i \bigmu_i^{(3)}+ \Lambda^{(3)} \right) \right]dz + \partial_s \left[M(U_i^{(0)}) \partial_s \left( \sigma_i \bigmu_i^{(1)} + \Lambda^{(1)} \right) \right]dz,\\
& - \frac{1}{\nu_j} V_{ij} g(U_j^{(0)}) \partial_z U_j^{(0)} = \partial_z\left[M_j^{(0)}\partial_z\left(\sigma_j \bigmu_j^{(3)}+ \Lambda^{(3)} \right) \right]dz + \partial_s \left[M(U_j^{(0)}) \partial_s \left( \sigma_j \bigmu_j^{(1)} + \Lambda^{(1)} \right) \right]dz.
\end{aligned}
\right.
\end{equation}
Integrating over $\mathbb{R}$ yields to
\begin{equation*}
\left.
\begin{aligned}
-& \frac{c_N}{\nu_i} V_{ij} = \int_\mathbb{R} \partial_z \left[M(U_i^{(0)})\partial_z\left( \sigma_i \bigmu_i^{(3)}+ \Lambda^{(3)}\right) \right]dz + \int_\mathbb{R} \partial_s \left[ M(U_i^{(0)})q'\partial_s \left(\sigma_i \bigmu_i^{(1)}+\Lambda^{(1)} \right)\right]dz,\\
+& \frac{c_N}{\nu_j} V_{ij} = \int_\mathbb{R} \partial_z \left[M(U_j^{(0)})\partial_z\left( \sigma_j \bigmu_j^{(3)}+ \Lambda^{(3)}\right) \right]dz + \int_\mathbb{R} \partial_s \left[ M(U_j^{(0)})q'\partial_s \left(\sigma_j \bigmu_j^{(1)}+\Lambda^{(1)} \right)\right]dz. \\
\end{aligned}
\right.
\end{equation*}
The matching conditions for the fluxes at order $\mathcal{O}(\varepsilon^2)$ show that the first integral is zero.
Note that most of the terms in the fluxes have been proven to be zero  in the previous orders. 

On the other hand, the second integral can be expressed with the terms from the second order calculations and the properties of the profile $q$
give that
\begin{equation*}
\left.
\begin{aligned}
\int_\mathbb{R} \partial_s\left[ M(U_i^{(0)})\partial_s\left( \sigma_i \bigmu_i^{(1)} + \Lambda^{(1)}\right) \right]dz & = \partial_{ss} \left(\sigma_i \bigmu_i^{(1)}+\Lambda^{(1)}  \right) \int_\mathbb{R} M(q)dz, \\
& = c_M \partial_{ss}\left( \sigma_i \bigmu_i^{(1)}+\Lambda^{(1)} \right).\\
\end{aligned}
\right.
\end{equation*}
The same result can be obtained for the integral in $j$. Subtracting the first equation of \eqref{Eq:NMNFourthorder} to the second, we get
$$ \frac{1}{\nu_{ij}} c_N V_{ij}  = \left(\frac{1}{\nu_i}+\frac{1}{\nu_j} \right) c_N V_{ij}
= c_M \partial_{ss} \left( - \sigma_i^{(1)} \bigmu_i^{(1)} - \Lambda^{(1)} + \sigma_j \bigmu_j^{(1)} + \Lambda^{(1)} \right) =  c_M \partial_{ss}A_{ij}^{(1)}.$$
Using \eqref{Eq:MultiphaseAIJNMN}, it follows that 
\begin{equation*}
\frac{1}{\nu_{ij}}V_{ij} = \frac{c_W c_M}{(c_N)^2} \sigma_{ij} \partial_{ss} H_{ij},
\end{equation*}
which concludes the proof of Proposition~\ref{result:modelNMN}.

\section{Numerical approximation}
\label{Sec:NumericalSection}

In this section, we show how to compute effectively numerical approximations of the solutions to phase field models \MCH ~and \NMNCH,  
and we provide various numerical illustrations of the performances and properties of both models 
in dimensions $2$ and $3$. The numerical approximation is performed with the original \MCH ~model and with the \NMNCH ~reformulation II model (see the introduction), whose definitions are recalled:

\begin{itemize}
 \item {\MCH} \\
 $$
 \begin{cases}
  \partial_t u_k &= \nu_k\div\left(M(u_k)\nabla( \sigma_k \mu_k + \lambda) \right) \\
   \mu_k & = \frac{W'(u_k)}{\varepsilon^2} - \Delta u_k \\
   1  &= \sum_k u_k
 \end{cases}
 $$
 where the mobility $M$ is defined by $M(u) = \frac{1}{c_N^2} 2 W(u)$. Here, the constant $|c_N| = \frac{1}{6}$ is added
to get the same limit law as with  our new Cahn--Hilliard model.
 \item {\NMNCH\text{-reformulation II}} \\
$$
 \begin{cases}
  \partial_t u_k &= \nu_k N(u_k) \div\left(M(u_k)\nabla N(u_k) ( \sigma_k \mu_k + \lambda) \right) \\
   \mu_k & = \frac{W'(u_k)}{\varepsilon^2} - \Delta u_k \\ 
   1 &= \sum_k u_k
 \end{cases}
 $$
 \end{itemize}
 where the mobilities $M, N$ are defined by $M(s) = 2W(s) + \gamma \varepsilon^2$, with $\gamma>0$, and $N(s) = \frac{1}{\sqrt{M(s)}}$. For simplicity we keep the original notations $M,N$ although the definitions are different, and still for the sake of simplicity we drop the expression "reformulation II". We set $\gamma = 1$ for all numerical experiments presented below.\\

Various schemes have already been proposed in the literature, see \cite{YANG2020105276,Li2016,LEE20121009,LEE20084787}, 
to deal with multiphase Cahn-Hilliard type equations, especially when the number of phases is $L=3$  \cite{barrett_parametric_2007,boyer_study_2006,kang_conservative_2004,kim_phase_2007,kim_numerical_2009,copetti_numerical_2000,bhattacharyya_study_2003} or $L=4$~\cite{kim_generalized_2009,kitashima_phase-field_2008,zhou_multimaterial_2006,LEE20084787}. \\

Recall that the Cahn-Hilliard system is of fourth-order in space, 
which introduces severe restrictions
on the time step for most classical methods due to numerical instability.
To overcome these difficulties,  a natural idea is to adapt the strategy of convex splitting of the Cahn--Hilliard energy
which was first proposed by Eyre~\cite{MR1676409}. This technique has become very popular for it provides simple, efficient, and stable schemes to approximate
various evolution problems with a gradient flow structure~\cite{MR2418360,MR2519603,MR2799512,MR3100769,MR3564350,MR3682074,DU2020425}. 
For instance,  a first- and second-order splitting scheme was proposed in \cite{backofen2017convexity,salvalaglio2019doubly,Doubly_anisotropic}
to address the case of the Cahn--Hilliard equation with mobility. However, these approaches are based 
on finite elements and require the resolution of linear systems at each step, which can be ill-conditioned in the case of 
degenerate mobilities. As an alternative, we proposed recently in \cite{bretin2020approximation} a semi-implicit Fourier spectral method
in the spirit of
\cite{Chen_fourier,Bretin_brassel,bretin2017new,bretin_largephases,bretin2018multiphase}. The idea is to exploit the 
variational structure of the mobility  by using an additional convex splitting of the associated metric. It gives a very simple, effective, and stable scheme even in the case of degenerate mobilities. \\

An accurate non linear multigrid method was proposed in~\cite{LEE20084787} to approximate the solution to the Cahn-Hilliard equation.
However, this approach requires  the resolution of a $2L \times 2L$ system of equations which can be problematic when $L$ is large.
Based on the first-order convex splitting method, Lee et al~\cite{LEE20121009} developed a practically unconditionally gradient-stable conservative nonlinear numerical scheme  for converting the $L$-phase Cahn--Hilliard system into a system of 
$L$ Cahn--Hilliard equations. This reduces significantly the computational cost. More recently, Yang and Kim~\cite{YANG2020105276} 
proposed an unconditionally stable with second-order accuracy based on the Crank-Nicolson scheme and adopted the idea of stabilized method~\cite{PhysRevE.60.3564}. \\

In this paper, we extend to multiphase the approach we proposed in~\cite{bretin2020approximation}.
The novelty is to split the treatment of the Lagrange multiplier via the splitting of the metric so that 
\MCH ~and \NMNCH ~can 
be solved in
a decoupled way. This means that we only need to solve L biphasic Cahn--Hilliard equations at each iteration, 
as in~\cite{LEE20121009}.  \\
   
In the following, we first recall the schemes we have introduced in \cite{bretin2020approximation} 
when only two phases (i.e. one single function $u$) are considered for both \MCH ~and \NMNCH ~models. 
Then we extend to the multiphase case by using a semi-implicit treatment 
of the Lagrange multiplier which is explicitly given in Fourier space. For each model, a {\bf Matlab} script is provided
to give an example of implementation. Next, we provide a numerical comparison of phase field models in space dimension 2. 
In addition, some illustrations are provided to show the influence of mobilities and surface tensions using 
the \NMNCH ~model.  These illustrations show also that our models can handle Cahn--Hilliard problems in complex domain without imposing any boundary  condition or additional surface energy, but rather by simply imposing a null mobility at the appropriate interfaces.  
Then we conclude the section with an application to the dewetting problem using a simplified model that involves the liquid phase only.and that is equivalent in this particular context to \NMNCH.
 
\subsection{Spatial and time discretization: a Fourier-spectral approach}
All equations are solved on a square-box $Q = [0,L_1]\times \cdots \times [0,L_d]$ with periodic boundary conditions.
We recall that the Fourier $\boldsymbol K$-approximation of a function $u$ defined in a box 
$Q = [0,L_1]\times \cdots \times [0,L_d]$ is given by
$$u^{\boldsymbol K}(x) = \sum_{{\boldsymbol k}\in K_N  } c_{\boldsymbol k} e^{2i\pi{\boldsymbol \xi}_k\cdot x},$$
where  $K_N =  [ -\frac{N_1}{2},\frac{N_1}{2}-1 ]\times [ -\frac{N_2}{2},\frac{N_2}{2}-1] \cdots \times   [ -\frac{N_d}{2},\frac{N_d}{2}-1] $,   ${\boldsymbol k} = (k_1,\dots,k_d)$ and ${\boldsymbol \xi_k} = (k_1/L_1,\dots,k_d/L_d)$. In this formula, the $c_{\boldsymbol k}$'s denote the $K^d$ first discrete Fourier coefficients of $u$. 
The inverse discrete Fourier transform leads to 
$u^{K}_{\boldsymbol k} =   \textrm{IFFT}[c_{\boldsymbol k}]$ 
where $u^{K}_{\boldsymbol k}$ denotes the value of $u$ at the points 
$x_{\boldsymbol k} = (k_1 h_1,\dots, k_d h_d)$ and where $h_{i} = L_{i}/N_{i}$ for $i\in\{1,\dots,d\}$. Conversely,
$c_{\boldsymbol k}$ can be computed as the discrete Fourier transform of $u^K_{\boldsymbol k},$ {\em i.e.}, $c_{\boldsymbol k} =
\textrm{FFT}[u^K_{\boldsymbol k}].$

Given a time discretization parameter $\delta_t > 0$, we construct a sequence $(u^n)_{n \geq 0}$ 
of approximations of ${u}$ at times $n \delta_t$.

\subsection{Numerical scheme for the \MCH ~model}
  
We first recall  the numerical approach introduced in \cite{bretin2020approximation} to 
compute numerical solutions of the \MCH ~model in a biphasic context. 
In such a case, the Cahn--Hilliard equation reads as 
  $$
 \begin{cases}
  \partial_t u &=  \div\left(M(u)\nabla( \mu) \right),\\
   \mu & = \frac{W'(u)}{\varepsilon^2} - \Delta u. \\
 \end{cases}
 $$
 
Our approach can be viewed as  a Fourier semi-implicit scheme which reads as
 $$
\begin{cases}
 (u^{n+1} - u^{n})/\delta_t &=  m \Delta \mu^{n+1} +  \div( (M(u^n) - m) \nabla \mu^n ) \\
 \mu^{n+1} &=  \left( - \Delta u^{n+1} + \frac{\alpha}{\varepsilon^2} u^{n+1} \right) + \left( \frac{1}{\varepsilon^2} (W'(u^{n}) - \alpha u^{n}) \right).   
\end{cases}
$$
where $m$ and $\alpha$ are two stabilization parameters. More precisely, this scheme derives from a convex-concave splitting of 
the Cahn--Hilliard energy 
$$  \int_Q (\frac{|\nabla u|^2}{2} +  \frac{1}{\varepsilon^2} W(u)) dx =  \frac{1}{2}\int_Q(|\nabla u|^2  + \frac{\alpha}{\varepsilon^2} u^2) dx + \int_{Q} \frac{1}{\varepsilon^2} (W(u) - \alpha \frac{u^2}{2})dx,  $$
but also of  the associated metric
$$  \frac{1}{2} \int_Q M(u^{n})|\nabla \mu|^2 dx   = \frac{1}{2} \int_Q m |\nabla \mu|^2 dx +  \frac{1}{2} \int (M(u^{n})-m) |\nabla \mu|^2 dx.$$
As we explained in \cite{bretin2020approximation}, the scheme seems to decrease the Cahn--Hilliard energy 
as soon as each  explicit term is concave, which is true when setting $m  = \max_{s \in [0,1]} M(s)$ 
and $\alpha \geq \max_{s \in [0,1]} \left| W''(s) \right|$. \\

Alternatively this scheme reads in a matrix form as
$$ \begin{pmatrix}
    I_d &  - \delta_t m \Delta \\
    \Delta - \alpha/\varepsilon^2 & I_d
   \end{pmatrix}  
   \begin{pmatrix}
    u^{n+1} \\
    \mu^{n+1}
   \end{pmatrix}
   =    \begin{pmatrix}
    B^1_{u^{n},\mu^n}\\
    B^2_{u^{n},\mu^n}
   \end{pmatrix},   
   $$
   where 
   $$
    \begin{pmatrix}
    B^1_{u^{n},\mu^n}\\
    B^2_{u^{n},\mu^n}
   \end{pmatrix} =
    \begin{pmatrix}
    u^{n} +  \delta_t \div( (M(u^n) - m) \nabla \mu^n ) \ \\
     \frac{1}{\varepsilon^2} (W'(u^{n}) - \alpha u^{n})
   \end{pmatrix}.
  $$
Finally, the couple $(u^{n+1},\mu^{n+1})$ can be computed using the system
$$
\begin{cases}
 u^{n+1} &= L_M \left[   B^1_{u^{n},\mu^n} +  \delta_t m \Delta B^2_{u^{n},\mu^n} \right],\\
 \mu^{n+1}&= L_M\left[ ( - \Delta  B^1_{u^{n},\mu^n} + \alpha/\varepsilon^2 B^1_{u^{n},\mu^n} ) +   B^2_{u^{n},\mu^n} \right],
\end{cases}
$$
where the operator
$$L_M = \left( I_d +  \delta_t m \Delta ( \Delta - \alpha/\varepsilon^2 I_d) \right)^{-1},$$
can be computed very efficiently in Fourier space.\\

\begin{rem}
This scheme is very efficient as it does not require any resolution of a linear system.
Moreover, this scheme seems to be stable without assumption on $\delta_t$ in the sense that it decreases
the Cahn-Hilliard energy. It is also not difficult to show  that the mass of $u$ is conserved along the iterations, i.e.,
 $$ \int_Q u^{n+1} dx = \int_Q u^{n} dx.$$
\end{rem}

 Following this method, we now propose a similar scheme for the multiphase \MCH ~model
  $$
 \begin{cases}
  \partial_t u_k &= \nu_k\div\left(M(u_k)\nabla( \sigma_k \mu_k + \lambda) \right), \\
   \mu_k & = \frac{W'(u_k)}{\varepsilon^2} - \Delta u_k, \\
   1  &= \sum_k u_k.
 \end{cases}
 $$
 
 The scheme we propose is based  on the same convex-concave splitting of
 the Cahn-Hilliard equation and its associated metric. In the multiphase context, we obtain
 $$
\begin{cases}
 (u^{n+1}_k - u^{n}_k)/\delta_t &=  \nu_k \left( m  \Delta \left[ \sigma_k \mu^{n+1}_k + \lambda^{n+1} \right] +  \div \left[ (M(u^n_k) - m) \nabla \left[\sigma_k \mu^n_k + \lambda^n \right] \right] \right),\\
 \mu^{n+1}_k &=  \left( - \Delta u^{n+1}_k + \frac{\alpha}{\varepsilon^2} u^{n+1}_k \right) + \left( \frac{1}{\varepsilon^2} (W'(u^{n}_k) - \alpha u^{n}_k) \right),   
\end{cases}
$$
where the Lagrange multiplier $\lambda^{n+1}$ is associated to the partition constraint  $\sum_k u^{n+1}_k = 1$. \\

More precisely, the couple  $(u^{n+1}_k, \mu^{n+1}_k)$ can be expressed as 
  $$
  \begin{cases}
  u^{n+1}_k &= u^{n+1/2}_k + \delta_t \nu_k m L_{M_k} \left[ \Delta \lambda^{n+1}   \right], \\
  \mu^{n+1}_k &= \mu^{n+1/2}_k + \delta_t  \nu_k m L_{M_k} \left[ (- \Delta + \frac{\alpha}{\varepsilon^2}) \Delta \lambda^{n+1} \right], 
  \end{cases}
  $$
 where
 \begin{itemize}
  \item the operator $L_{M_k}$ is given by 
  $$L_{M_k} = \left( I_d +  \delta_t m \sigma_k \nu_k \Delta ( \Delta - \alpha/\varepsilon^2 I_d) \right)^{-1}.$$
  \item the couple $(u^{n+1/2}_k, \mu^{n+1/2}_k)$ is defined as the solution to the decoupled system
  $$
\begin{cases}
 u^{n+1/2}_k &= L_{M_k} \left[   B^1_{u^{n}_k,\mu_k^n} +  \delta_t m \sigma_k \nu_k \Delta B^2_{u^{n}_k,\mu^n_k} \right],\\
 \mu^{n+1/2}_k &= L_{M_k}\left[ ( - \Delta  + \alpha/\varepsilon^2) B^1_{u^{n}_k,\mu^n_k}  +   B^2_{u^{n}_k,\mu^n_k} \right],
\end{cases}
 $$
 where
 $$ 
   B^1_{u^{n}_k,\mu^n_k} =  u^{n}_k +  \delta_t \nu_k \div \left[ (M(u^n_k) - m) \nabla \left[  \sigma_k \mu^n_k + \lambda^k \right] \right],$$
  and 
  $$B^2_{u^{n}_k,\mu^n_k} =  \frac{1}{\varepsilon^2} (W'(u^{n}_k) - \alpha u^{n}_k).
$$
 \end{itemize}

In particular, $\lambda^{n+1}$ satisfies the equation  
$$ \delta_t m \left( \sum_k \nu_k L_{M_k} \right) \Delta  \lambda^{n+1} = 1 - \sum_k u^{n+1/2}_k,$$
therefore,
$$ \lambda^{n+1} = \frac{1}{\delta_t m} \left[ \sum_k \nu_k L_{M_k} \Delta   \right]^{-1}( 1 - \sum_k u^{n+1/2}_k).$$
Here the operator   $\left[ \sum_k \nu_k  L_{M_k}  \Delta \right]^{-1}$ is still homogeneous and can be computed
easily  in Fourier space. \\


From the previous equations, we can implement the scheme within the {\bf Matlab} framework almost 
immediately, see the $54$-lines {\bf Matlab} 
script of Table~\ref{fig:matlab_code_model_M} which approximates the solution to the  \MCH ~model.
In particular :

\begin{itemize}
 \item We consider here a discretized  computation box $Q = [-1/2,1/2]^2$  using $N = 2^8$ nodes in each direction. 
 The initial condition of $u$ is a uniform noise and the numerical parameters are set to $\varepsilon = 1/N$, $\delta_t = \varepsilon^4$,  
 $\alpha = 2$, and $m= \max_{s\in [0,1]} M(s)$. 
 
 \item First we define the terms $u_k^{n+1/2}$ and $\mu_k^{n+1/2}$ (lines 29-39) as in \cite{bretin2020approximation}.
 Then, we determine $\lambda^{n+1}$ (lines 42-45) which allows us to correct and obtain $u_k^{n+1}$ and $\mu_k^{n+1}$ (lines 48-52).

 \item Line $24$ corresponds to the definition of the Fourier-symbol associated with the operator $L_{M_k}$. The application of $L_{M_k}$ can then be performed by using a simple multiplication in Fourier space with the array $M_L$.
 
  \item Each computation of a gradient or a divergence is made in Fourier space. For instance the 
  divergence $\div \left[ (M(u^n_k) - m) \nabla \left[  \sigma_k \mu^n_k + \lambda^k \right] \right]$ is computed on line $31$.
   
 \item The computation of $\lambda^{n+1}$ is illustrated on lines 42-45. $\lambda^{n+1}$ is first computed in Fourier space 
 using the Fourier-symbol of the operator $\left[ \sum_k \nu_k  L_{M_k}  \Delta \right]^{-1}$. Then $\lambda^{n+1}$
 is obtained by applying the discrete inverse Fourier transform.
 \end{itemize}


 \begin{table}[htbp]
\centering
\begin{lstlisting}
clear all; colormap('jet');

 %%%%%%%%%%%%%%%%% parameters %%%%%%%%%%%%%%%%%%%%%%%%
N = 2^8; epsilon =1/N; dt = epsilon^4; T = 10^(-4);
W_prim = @(U) (U.*(U-1).*(2*U-1));
MobMM = @(U) 2*36/2*((((U).*(1-U)).^2) );
alpha = 2;x = linspace(0,1,N); c=max(MobMM(x));
%%%%%%%%%%%%%%%% initial conditions  %%%%%%%%%%%%%%
U(:,:,1) = 2*rand(N,N)/3; U(:,:,2) = rand(N,N).*(1 -U(:,:,1) );
U(:,:,3) = 1-(U(:,:,1) + U(:,:,2) );
Mu = 0*U;lambda = 0; lambda_fourier = 0;Mu_fourier = zeros(N,N,3);
for k=1:3, U_fourier(:,:,k) = fft2(U(:,:,k)); end

%%%%%%%%%%%%%%%%%%% surface tension coefficients - mobilities %%%%%%%%%%
sigma12 =1; sigma13 =1; sigma23 =1;
sigma(1) = (sigma12 + sigma13 - sigma23)/2;
sigma(2) = (sigma12 + sigma23 - sigma13)/2;
sigma(3) = (sigma23 + sigma13 - sigma12)/2;
mob(1) = 1; mob(2)=  1; mob(3) = 1;
%%%%%%%%%%%%%%%%% Diffusion Fourier %%%%%%%%%%%%%%%
k =  [0:N/2,-N/2+1:-1];
[K1,K2] = meshgrid(k,k);
Delta = -4*pi^2*((K1.^2 + (K2).^2));
for k=1:3, M_L(:,:,k) = 1./(1 + dt*sigma(k)*mob(k)*(c*Delta).*(Delta - alpha/epsilon^2)); end
%%%%%%%%%%%%%%%%%%%%%%%%%%%%%%%%%%%%%%%%%%%%%%%%%%
k=1;
for i=1:T/dt,
%%%%%%%%%%%% computation of u^{n+1/2}_k, \mu^{n+1/2}_k
for k=1:3,    
mobUk = MobMM(U(:,:,k));
div_mob_laplacien_fourier =  2*1i*pi*K1.*fft2((mobUk-c).*ifft2(2*1i*pi*K1.*(sigma(k)*Mu_fourier(:,:,k) + lambda_fourier)))...
    + 2*1i*pi*K2.*fft2((mobUk-c).*ifft2(2*1i*pi*K2.*(sigma(k)*Mu_fourier(:,:,k) + lambda_fourier)));
B1 = U_fourier(:,:,k)  +  dt*mob(k)*div_mob_laplacien_fourier;
B2 = fft2(W_prim(U(:,:,k))/epsilon^2 - alpha/epsilon^2*U(:,:,k));
U_fourier(:,:,k) = M_L(:,:,k).*(B1 + dt*mob(k)*sigma(k)*c*Delta.*B2);
U(:,:,k) = real(ifft2(U_fourier(:,:,k)));
Mu_fourier(:,:,k) = M_L(:,:,k).*((alpha/epsilon^2 - Delta).*B1 + B2);  
Mu(:,:,k) = ifft2(Mu_fourier(:,:,k));
end

%%%%%%%%%%%%%%%%%%%% computation of lambda and correction %%%%%%%%%%
Err_sum = fft2(1 - sum(U,3));
lambda_fourier = (1./((c*(M_L(:,:,1)*mob(1) +M_L(:,:,2)*mob(2) + M_L(:,:,3)*mob(3))).*Delta)).*Err_sum/dt;
lambda_fourier(1,1) = 0;
lambda = real(ifft2(lambda_fourier));

for k=1:3,
term = dt*(mob(k)*(c*Delta.*lambda_fourier));
U_fourier(:,:,k) = U_fourier(:,:,k) + M_L(:,:,k).*term;
Mu_fourier(:,:,k) = Mu_fourier(:,:,k) - (Delta - alpha/epsilon^2).*M_L(:,:,k).*term; 
U(:,:,k) = real(ifft2(U_fourier(:,:,k)));
Mu(:,:,k) = ifft2(Mu_fourier(:,:,k));    
end
end






\end{lstlisting}
	 \caption{{\bf Matlab} implementation of our scheme to approximate in dimension $2$ the solutions to the \MCH ~model.}
\label{fig:matlab_code_model_M}
\end{table}

\subsection{Numerical scheme for the \NMNCH ~model}
The case of the \NMNCH ~model is slightly more complicated.  We first recall the numerical scheme 
introduced in \cite{bretin2020approximation} for only two phases, then we explain how to generalize it in the multiphase context. 
Recall that the \NMNCH ~model reads in the biphase case as
  $$
 \begin{cases}
  \partial_t u &=  N(u) \div\left(M(u)\nabla( N(u) \mu ) \right) \\
   \mu & = \frac{W'(u)}{\varepsilon^2} - \Delta u \\
 \end{cases}
 $$
 and that the Fourier semi-implicit  scheme we proposed in~\cite{bretin2020approximation} to approximate its solutions is
  $$
\begin{cases}
 (u^{n+1} - u^{n})/\delta_t &=  [m \Delta - \beta ] \mu^{n+1}  +  H(u^n,\mu^n),\\
 \mu^{n+1} &=  \left( - \Delta u^{n+1} + \frac{\alpha}{\varepsilon^2} u^{n+1} \right) + \left( \frac{1}{\varepsilon^2} (W'(u^{n}) - \alpha u^{n}) \right),   
\end{cases}
$$
where 
$$   H(u^n,\mu^n) =  N(u^n)\div( (M(u^n) \nabla (N(u^{n}) \mu^n) ) - m \Delta \mu^{n} + \beta \mu^{n}.$$

\begin{rem}
Recall that  this approach is based on the convex-concave splitting of the associated metric
$$\frac{1}{2} \int_Q  M(u) \left| \nabla ( N(u) \mu) \right|^2 dx =   J_{u,c}(\mu)  + J_{u,e}(\mu),$$
with
 $$ J_{u,c}(\mu) =   \frac{1}{2} \int_Q m |\nabla \mu|^2 dx + \frac{1}{2}  \int_Q \beta \mu^2  dx$$
 and
 $$ J_{u,e}(\mu) =  \int_Q G(u) \cdot \nabla \mu \mu dx + \frac{1}{2}  \int_Q (|G(u)|^2 - \beta) \mu^2  dx + \frac{1}{2} \int_Q (1 -m) |\nabla \mu|^2 dx.$$
 Here, $G(u) = -\frac{1}{2} \nabla(\log(M(u)))$, and as it is bounded is $H^1(Q)$,  a sufficiently large choice for
 $m$ and $\beta$ should  ensure the concavity of $J_{u,e}(\mu)$ and the stability of the scheme. In practice, we set $m=1$ and 
 $\beta = 1/\varepsilon^2$ for our numerical experiments and with these values we did not observe any sign of instability regardless of the choice of the time step $\delta_t$. \\
 \end{rem}
 
Finally, the couple $(u^{n+1},\mu^{n+1})$ is  solution of the system
$$ \begin{pmatrix}
    I_d &  - \delta_t (m \Delta - \beta I_d)  \\
    \Delta - \alpha/\varepsilon^2 & I_d
   \end{pmatrix}  
   \begin{pmatrix}
    u^{n+1} \\
    \mu^{n+1}
   \end{pmatrix}
   = \begin{pmatrix}
    B^1_{u^{n},\mu^n}\\
    B^2_{u^{n},\mu^n}
   \end{pmatrix},
   $$
  with 
   $$
   \begin{pmatrix}
    B^1_{u^{n},\mu^n}\\
    B^2_{u^{n},\mu^n}
   \end{pmatrix}
   =
    \begin{pmatrix}
    u^{n} +  \delta_t H(u^n,\mu^n)\ \\
     \frac{1}{\varepsilon^2} (W'(u^{n}) - \alpha u^{n}),
   \end{pmatrix},
   $$
   and satisfies
$$
\begin{cases}
 u^{n+1} &= L_{NMN} \left[   B^1_{u^{n},\mu^n} +  \delta_t ( m \Delta B^2_{u^{n},\mu^n} - \beta B^2_{u^{n},\mu^n} ) \right] \\
 \mu^{n+1}& = L_{NMN}\left[ ( - \Delta  B^1_{u^{n},\mu^n} + \alpha/\varepsilon^2 B^1_{u^{n},\mu^n} ) +   B^2_{u^{n},\mu^n} \right].
\end{cases}
$$
where the operator $L_{NMN} = \left( I_d +  \delta_t (m \Delta - \beta I_d)( \Delta - \alpha/\varepsilon^2 I_d) \right)^{-1}$ can be computed very efficiently in Fourier space. \\

We now propose to extend this approach to the multiphase case:

  $$
 \begin{cases}
  \partial_t u_k &= \nu_k N(u_k) \div \left[ M(u_k)\nabla[  N(u_k)  \sigma_k \mu_k + \lambda )] \right], \\
   \mu_k & = \frac{W'(u_k)}{\varepsilon^2} - \Delta u_k, \\
   1  &= \sum_k u_k.
 \end{cases}
 $$
 
The scheme reads 
   $$
\begin{cases}
 (u^{n+1}_k - u^{n}_k)/\delta_t &= \nu_k [m \Delta - \beta ] (\sigma_k \mu^{n+1}_k + \lambda^{n+1})   +  H_k(u^n_k,\mu^n_k,\lambda^n) \\
 \mu^{n+1}_k &=  \left( - \Delta u^{n+1}_k + \frac{\alpha}{\varepsilon^2} u^{n+1}_k \right) + \left( \frac{1}{\varepsilon^2} (W'(u^{n}_k) - \alpha u^{n}_k) \right),   
\end{cases}
$$
 where 
 $$   H_k(u^n_k,\mu^n_k,\lambda^n) =  \nu_k \left( N(u_k^n)\div( (M(u_k^n) \nabla (N(u_k^{n}) (\sigma_k \mu_k^n + \lambda^n)) ) - [m \Delta - \beta](\sigma_k \mu_k^{n} + \lambda^n ) \right).$$
and $\lambda^{n+1}$ is associated to the partition constraint $\sum_k u_k^{n+1} = 1$.\\

Let us now introduce the couple $(u^{n+1/2}_k,\mu_k^{n+1/2})$ defined by
 
$$
\begin{cases}
 u_k^{n+1/2} &= L_{NMN,k} \left[   B^1_{u_k^{n},\mu_k^n,\lambda^n} +  \delta_t \nu_k \sigma_k ( [m \Delta - \beta]  B^2_{u_k^{n},\mu_k^n,\lambda^n} ) \right] \\
 \mu_k^{n+1/2}& = L_{NMN,k}\left[ ( - \Delta  B^1_{u_k^{n},\mu_k^n,\lambda^n} + \alpha/\varepsilon^2 B^1_{u_k^{n},\mu_k^n,\lambda^n} ) +   B^2_{u_k^{n},\mu_k^n,\lambda_k^n} \right].
\end{cases}
$$
where 

$$L_{NMN,k} = \left( I_d +  \delta_t \nu_k \sigma_k (m \Delta - \beta I_d)( \Delta - \alpha/\varepsilon^2 I_d) \right)^{-1}$$
and
 $$ 
    B^1_{u^{n},\mu^n,\lambda^n} =  u^{n} +  \delta_t H(u^n,\mu^n,\lambda^n) \text{ and } \quad  B^2_{u^{n},\mu^n,\lambda^n} =    \frac{1}{\varepsilon^2} (W'(u^{n}) - \alpha u^{n}).
$$
It is not difficult to see that 
  $$
  \begin{cases}
  u^{n+1}_k &= u^{n+1/2}_k + \delta_t \nu_k L_{NMN,k} \left[ [m \Delta - \beta] \lambda^{n+1}   \right] \\
  \mu^{n+1}_k &= \mu^{n+1/2}_k + \delta_t  \nu_k L_{NMN,k} \left[ (- \Delta + \frac{\alpha}{\varepsilon^2}) ( m \Delta - \beta) \lambda^{n+1} \right]   
  \end{cases}
  $$
which shows that $\lambda^{n+1}$ satisfies 
$$ \lambda^{n+1} = \frac{1}{\delta_t } \left[ \sum_k \nu_k L_{NMN,k} ( m\Delta - \beta)   \right]^{-1}( 1 - \sum_k u^{n+1/2}_k).$$
where the  operator   $ \left[ \sum_k \nu_k L_{NMN,k} ( m\Delta - \beta)   \right]^{-1}$ is homogeneous and can be, again, computed
easily in Fourier space. 
 
\begin{figure}[htbp]
\centering
	\includegraphics[width=3.5cm]{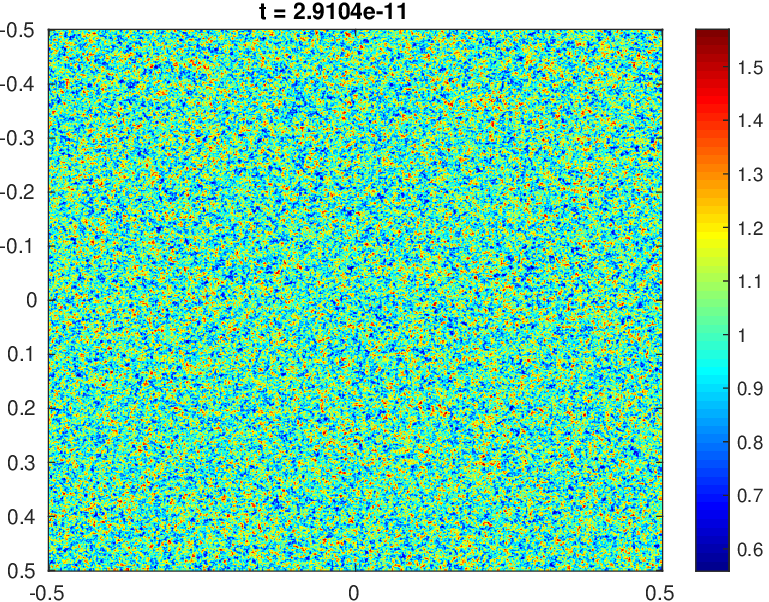}
	\includegraphics[width=3.5cm]{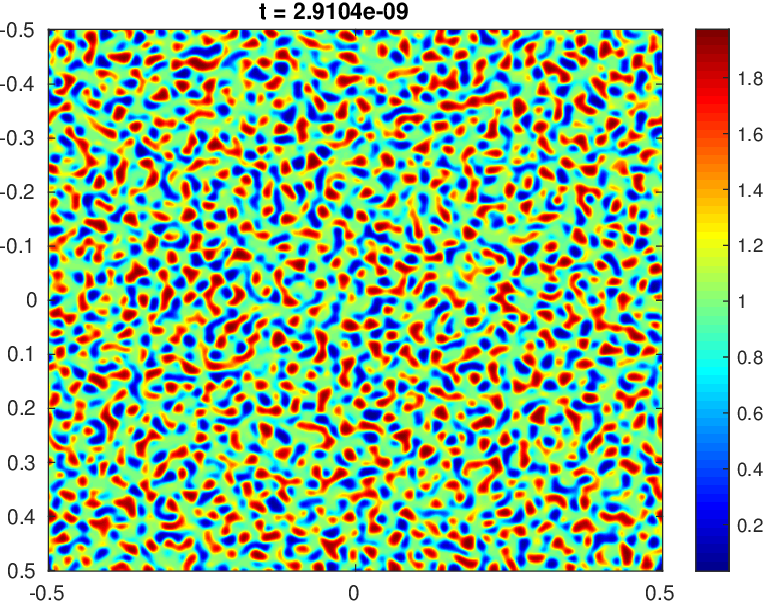}
	\includegraphics[width=3.5cm]{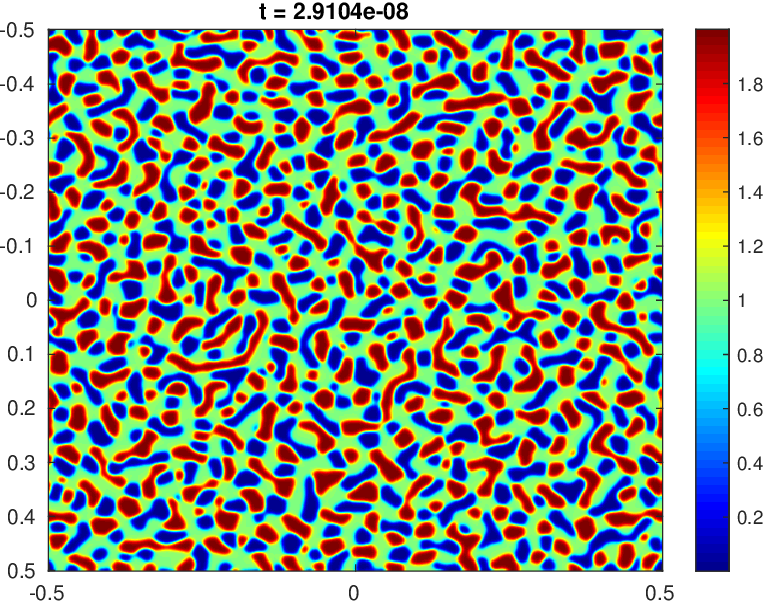}
	\includegraphics[width=3.5cm]{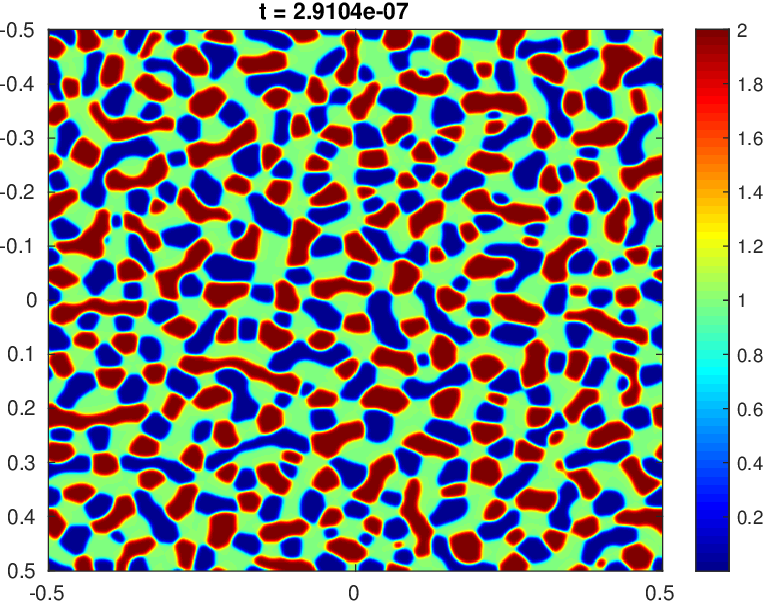}
\caption{First numerical experiment using the \NMNCH ~model; the solutions $u$ are computed with the {\bf Matlab}
script of Table~\ref{fig:matlab_code_model_NMN}. We plot the function $x \mapsto u_2(x) + 2u_3(x)$ on each picture so that the first, second, and third phases appear in blue, green, and red, respectively.}
\label{fig_Rand_NMN}
\end{figure}

Similarly to the \MCH ~model, this scheme approximating the solutions to the \NMNCH ~model can be easily implemented, see in Table~\ref{fig:matlab_code_model_NMN} a {\bf Matlab} script with less than $60$ lines.
In particular : 
\begin{itemize}
 \item We consider here a computation box $Q = [-1/2,1/2]^2$ discretized  with $N = 2^8$ nodes in each direction. 
 The initial condition of $u$ is a uniform noise and the numerical parameters are set to 
$\varepsilon = 1/N$, $\delta_t =  \varepsilon^4$,  
 $\alpha = 2$, $\beta = 2/\varepsilon^2$, and $m=1$.
    
 \item The implementation is almost identical to the previous model. Only the treatment of the divergence term $H_k(u_k^n,\mu_k^n,\lambda^n)$ makes a difference. The computation is done in lines $33$ to $36$ and is based on the
following equality: 

\begin{eqnarray*}
N(u) \div(M(u) \nabla (N(u)\mu)) &=&  \sqrt{M(u)} \Delta \left(N(u)\mu\right) + N(u) \nabla (M(u)) \cdot \nabla (N(u)\mu) \\
                                  &=&  \sqrt{M(u)} \Delta \left(N(u)\mu \right) +   2 \nabla \left[\sqrt{M(u)}\right] \cdot \nabla (N(u)\mu), 
\end{eqnarray*}

as $N(u) = 1/\sqrt{M(u)}$, see \cite{bretin2020approximation} for more details. 
 
\item Figure~\ref{fig_Rand_NMN} shows the function $u_2^n+2u_3^n$ computed at different times $t^n$ using this script.

\end{itemize}

\vspace*{3mm}

We believe that the proposed implementation illustrates well the simplicity, efficiency, and stability of our numerical scheme.

\begin{table}[htbp]
\centering
\begin{lstlisting}
clear all;
%%%%%%%%%%%%%%%%% parameters %%%%%%%%%%%%%%%%%%%%%%%%
N = 2^8; epsilon =1/N; dt = epsilon^4; T = 10^(-4);
alpha = 2;gamma=1; beta = 2/epsilon^2;
W_prim = @(U) (U.*(U-1).*(2*U-1));
MobM = @(U) 1/2*((((U).*(1-U)).^2+epsilon^2) );
MobN = @(U) 1./sqrt(MobM(U) );
%%%%%%%%%%%%%%%%% initial condition  %%%%%%%%%%%%%%
U(:,:,1) = 2*rand(N,N)/3; U(:,:,2) = rand(N,N).*(1 -U(:,:,1) );
U(:,:,3) = 1-(U(:,:,1) + U(:,:,2) );
Mu = 0*U;lambda = 0; lambda_fourier = 0;Mu_fourier = zeros(N,N,3);
for k=1:3, U_fourier(:,:,k) = fft2(U(:,:,k)); end
%%%%%%%%%%%%%%%%%%% surface tension coefficients - mobilities %%%%%%%%%%
sigma12 =1; sigma13 =1; sigma23 =1;
sigma(1) = (sigma12 + sigma13 - sigma23)/2;
sigma(2) = (sigma12 + sigma23 - sigma13)/2;
sigma(3) = (sigma23 + sigma13 - sigma12)/2;
mob(1) = 1; mob(2)=  1; mob(3) = 1;
%%%%%%%%%%%%%%%%% Kernel %%%%%%%%%%%%%%%
k =  [0:N/2,-N/2+1:-1]; [K1,K2] = meshgrid(k,k);
Delta = -4*pi^2*((K1.^2 + (K2).^2)); M_L = zeros(N,N,3);
for k=1:3,
M_L(:,:,k) = 1./(1 + dt*sigma(k)*mob(k)*(gamma*Delta - beta) .*(Delta - alpha/epsilon^2));
end
%%%%%%%%%%loop %%%%%%%%%%%%%%%%%%%
for i=1:T/dt,  
%%%%%%%%%%%% computation of u^{n+1/2}_k, \mu^{n+1/2}_k
for k=1:3,
mobMUk = MobM(U(:,:,k)); mobNUk = MobN(U(:,:,k));
sqrtMk = sqrt(mobMUk); sqrtMk_fourier = fft2(sqrtMk); 
nabla1_sqrtMk=  real(ifft2(2*pi*1i*K1.*sqrtMk_fourier )); nabla2_sqrtMk=  real(ifft2(2*pi*1i*K2.*sqrtMk_fourier ));
muN_fourier = fft2((sigma(k)*Mu(:,:,k) + lambda).*mobNUk);
nabla1_muN =  real(ifft2(2*pi*1i*K1.*muN_fourier ));
nabla2_muN =  real(ifft2(2*pi*1i*K2.*muN_fourier ));
laplacien_muN =  real(ifft2(Delta.*muN_fourier ));
NdivMgradNMu = sqrtMk.*laplacien_muN + 2*(nabla1_sqrtMk.*nabla1_muN +nabla2_sqrtMk.*nabla2_muN);
B1 = U_fourier(:,:,k) +  dt*(mob(k)*fft2(NdivMgradNMu) - mob(k)*(gamma*Delta-beta).*((sigma(k)*Mu_fourier(:,:,k)+lambda_fourier)));
B2 = fft2(W_prim(U(:,:,k))/epsilon^2 - alpha/epsilon^2*U(:,:,k));
U_fourier(:,:,k) = M_L(:,:,k).*(B1 + dt*mob(k)*sigma(k)*(gamma*Delta-beta).*B2);
U(:,:,k) = real(ifft2(U_fourier(:,:,k)));
Mu_fourier(:,:,k) = M_L(:,:,k).*((alpha/epsilon^2 - Delta).*B1 + B2);  
Mu(:,:,k) = real(ifft2(Mu_fourier(:,:,k)));
end
%%%%%%%%%%%%%%%%%%%% computation of lambda and correction %%%%%%%%%%
Err_sum = fft2(1 - sum(U,3));
weight =  (M_L(:,:,1)*mob(1) +M_L(:,:,2)*mob(2) + M_L(:,:,3)*mob(3)).*(gamma*Delta-1*beta);
lambda_fourier = (1./(weight)).*Err_sum/dt;
lambda = real(ifft2(lambda_fourier));
for k=1:3,
term = (mob(k)*M_L(:,:,k)./((M_L(:,:,1)*mob(1) +M_L(:,:,2)*mob(2) + M_L(:,:,3)*mob(3)))).*Err_sum;   
U_fourier(:,:,k) = U_fourier(:,:,k) + term;
Mu_fourier(:,:,k) = Mu_fourier(:,:,k) - (Delta - alpha/epsilon^2).*term; 
U(:,:,k) = real(ifft2(U_fourier(:,:,k)));
Mu(:,:,k) = real(ifft2(Mu_fourier(:,:,k)));
end
end






\end{lstlisting}
	 \caption{{\bf Matlab} implementation of our scheme to approximate in dimension $2$ the solutions to the \NMNCH ~model.}
\label{fig:matlab_code_model_NMN}
\end{table}

\subsection{Numerical validation}

\subsubsection{Asymptotic expansion and flow: a numerical comparison}

\begin{figure}[htbp]
\centering
	\includegraphics[width=3.8cm]{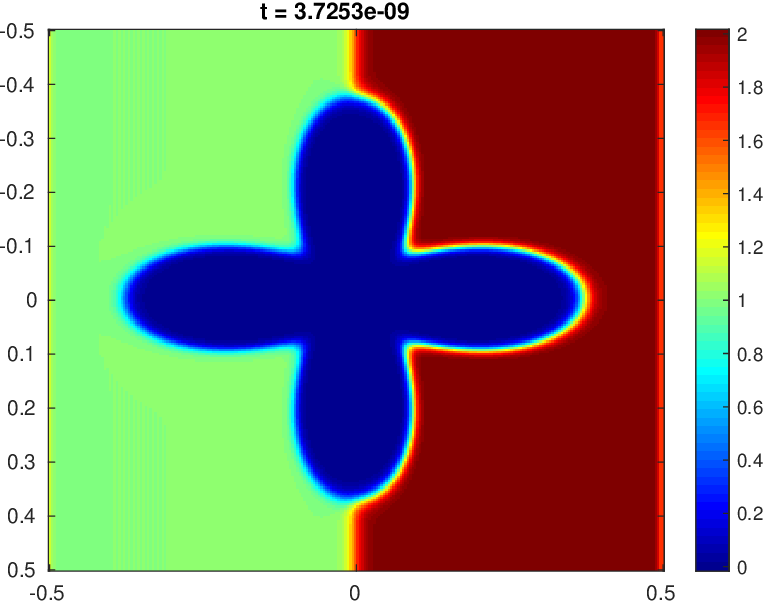}
	\includegraphics[width=3.8cm]{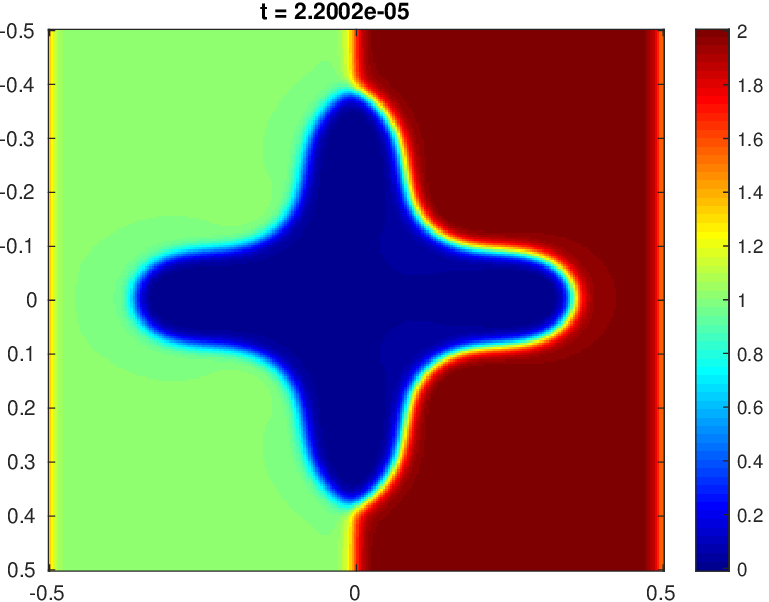}
	\includegraphics[width=3.8cm]{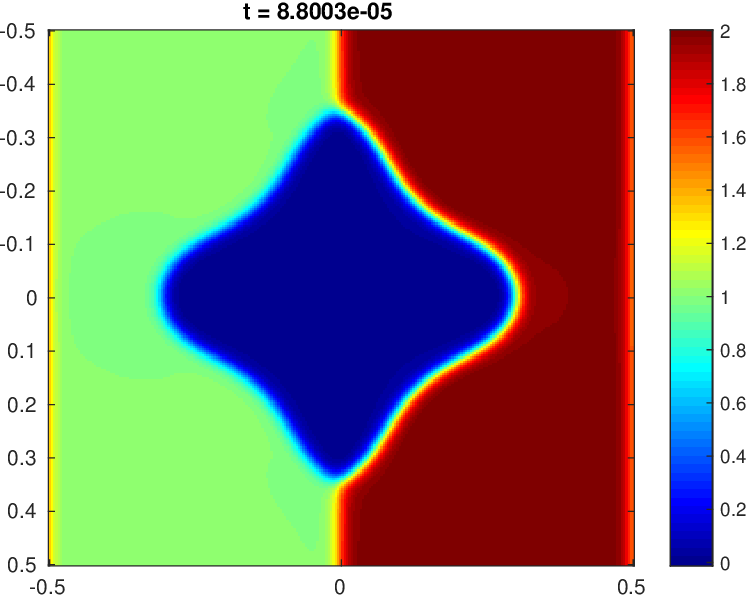}
	\includegraphics[width=3.8cm]{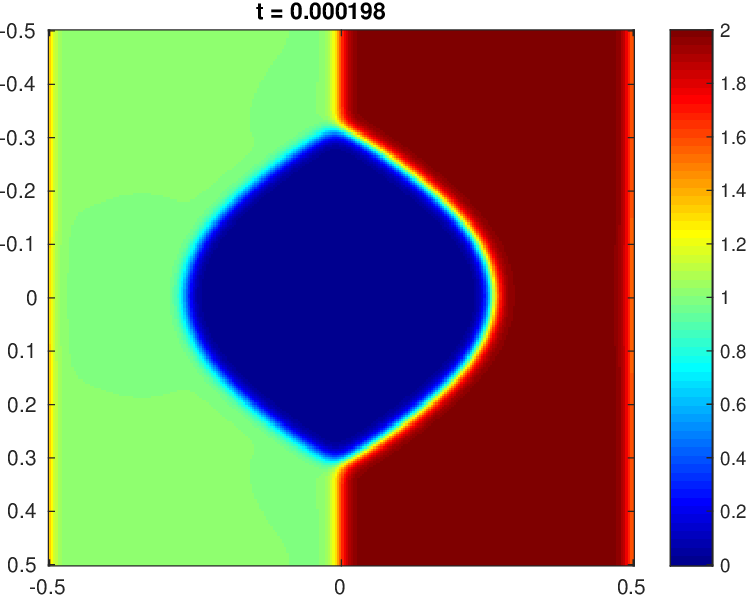} \\
	\includegraphics[width=3.8cm]{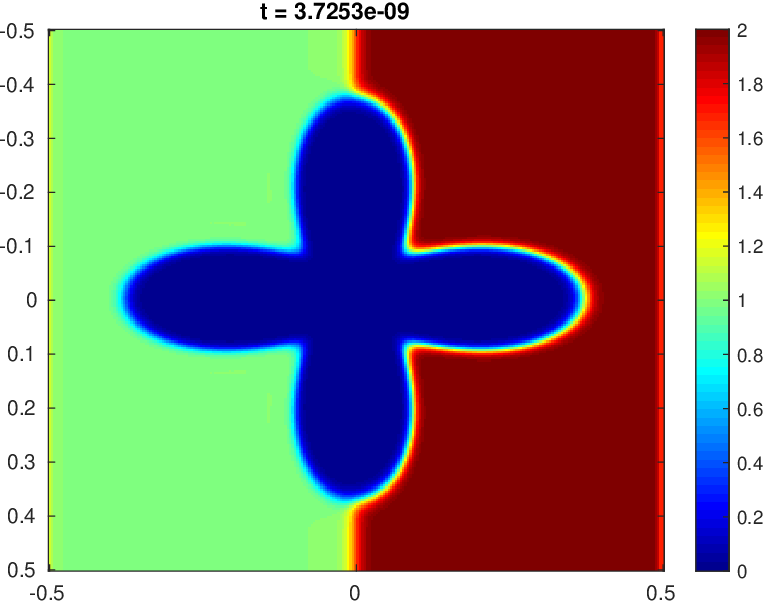}
	\includegraphics[width=3.8cm]{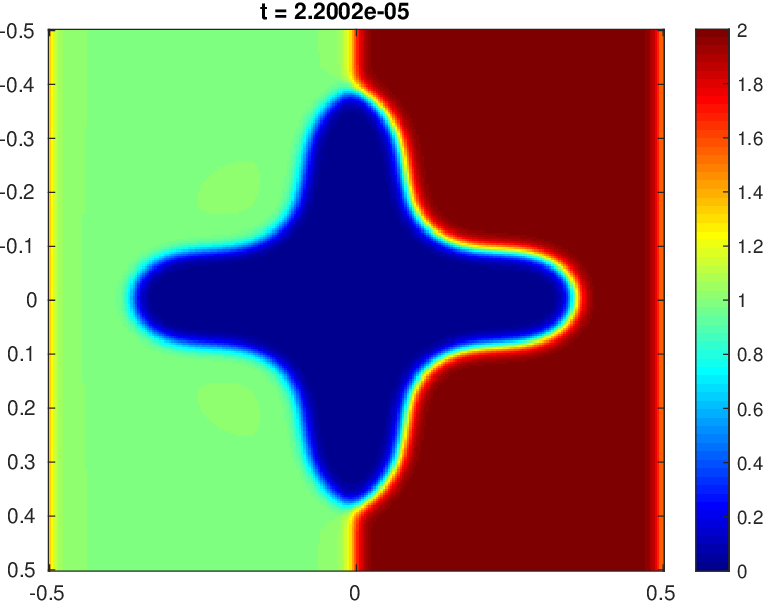}
	\includegraphics[width=3.8cm]{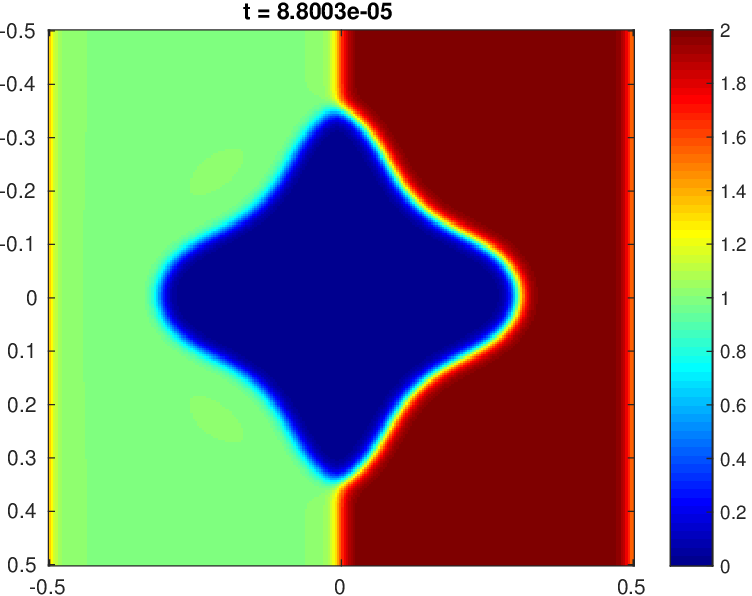}
	\includegraphics[width=3.8cm]{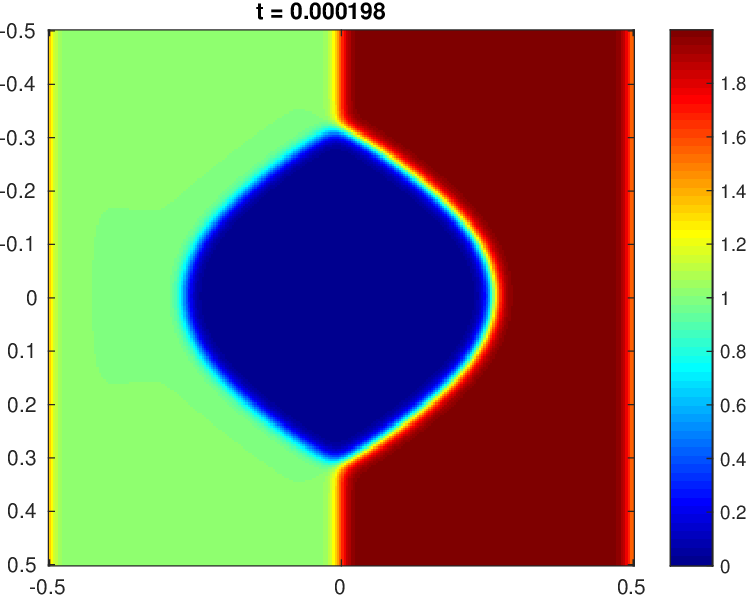} 
\caption{First numerical comparison of the two {\bf CH} models: evolution of the solution ${\bf u}$ along the iterations.
First line using \MCH, second line using \NMNCH.}
\label{fig_test2}
\end{figure}

The first numerical example concerns the evolution of an initial connected set. For each Cahn--Hilliard model,
we plot on Figure~\ref{fig_test2} the phase field function $u_2^{n} + 2 u_3^{n}$ computed at different times $t$. 
Each experiment is performed using the same numerical parameters:
$N = 2^8$, $\varepsilon = \delta_x$, $\delta_t = \varepsilon^4$,
$\alpha = 2$, $m=1$, and $\beta= 2/\varepsilon^2$. 

The first and second lines of Figure~\ref{fig_test2} correspond  to
the solutions given by the \MCH ~and the \NMNCH ~models, respectively.  
Notice that the numerical experiments obtained with both models are very similar and 
should give a good approximation of the surface diffusion flow.
In addition, for each model, the stationary flow limit appears to correspond to a ball of the same mass as that of the initial set. \\

\begin{figure}[htbp]
\centering
	\includegraphics[width=4.5cm]{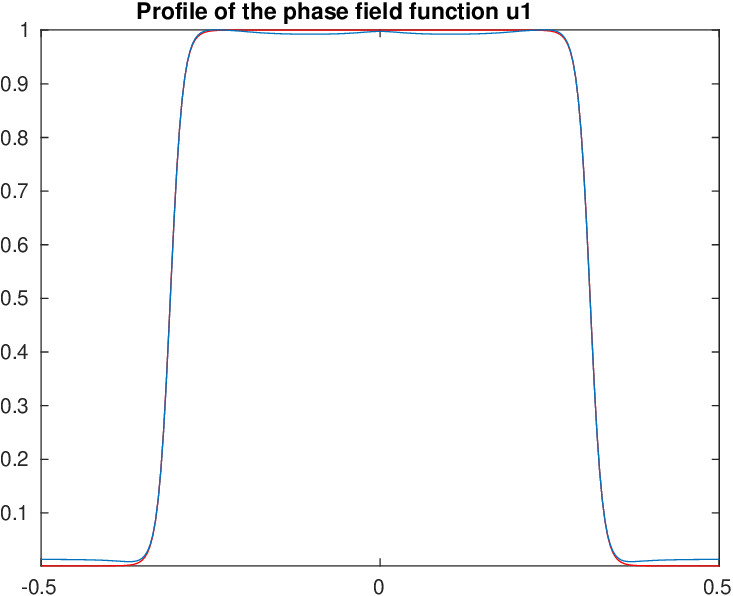}
	\includegraphics[width=4.5cm]{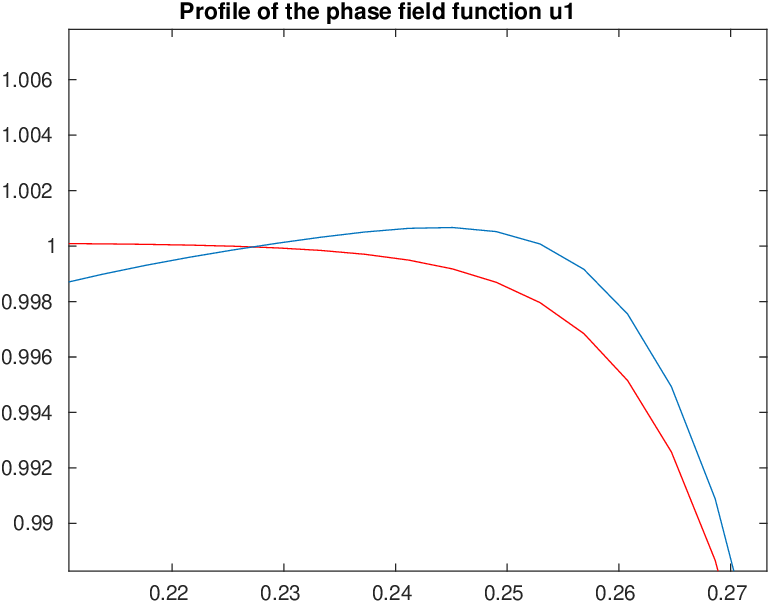}
	\includegraphics[width=4.5cm]{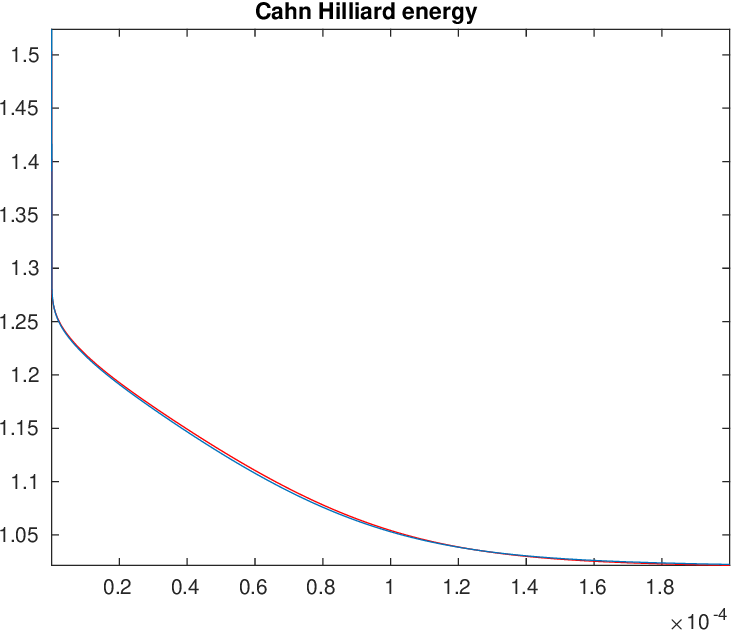}
\caption{Comparison of the two models regarding the solutions' profile and energy (\MCH ~in blue,  \NMNCH ~in red). 
Left figure: slice of $u$: $x_1 \mapsto u_1(x_1,0)$; Middle figure: zoom on the slice of $u_1$;  Right figure:
evolution of the Cahn--Hilliard energy along the flow.}
\label{fig_test2_profile}
\end{figure}

To illustrate  the asymptotic expansion performed in Section \ref{Sec:AsymptoticSection},
we plot on Figure~\ref{fig_test2_profile} (first two pictures) 
the slice $x_1 \mapsto u_1(x_1,0)$  at the final time $T = 10^{-4}$. The profile associated to the \MCH ~model 
is plotted in red and clearly indicates that the solution $u$ does not remain in the interval $[0,1]$  
with an overshoot of order $O(\varepsilon)$.   In contrast, the profile obtained using the \NMNCH ~
model (in green) seems to be very close to $q$ and remains in $[0,1]$ up to an error of order $O(\varepsilon^2)$. 
Finally, we plot the evolution of the Cahn--Hilliard energy 
along the flow for each model on the last picture of Figure~\ref{fig_test2_profile}. 
We can clearly observe a decrease of the energy in each case.\\

In conclusion, this first numerical experiment confirms the asymptotic expansion obtained in the previous section,
and highlights the interest of our \NMNCH ~model to approximate surface diffusion flows. 


\subsubsection{Influence of the mobility coefficients using the \NMNCH ~model }

\begin{figure}[htbp]
\centering
	\includegraphics[width=3.5cm]{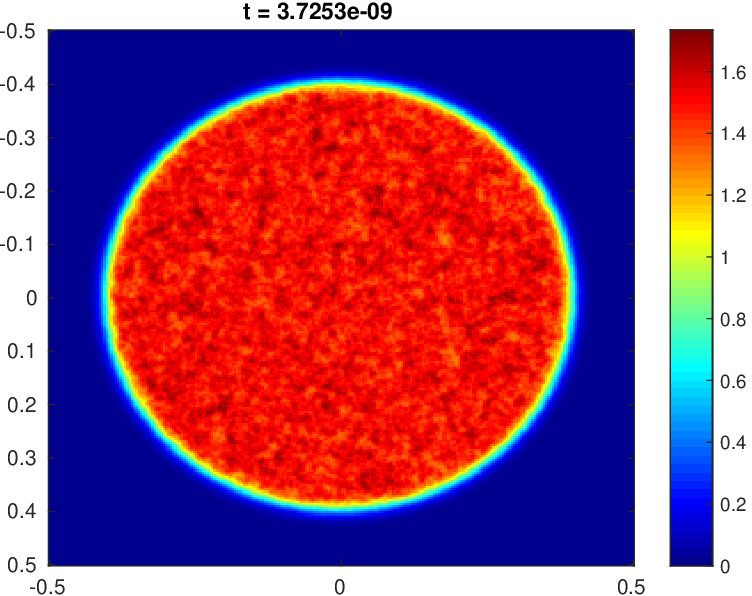}
	\includegraphics[width=3.5cm]{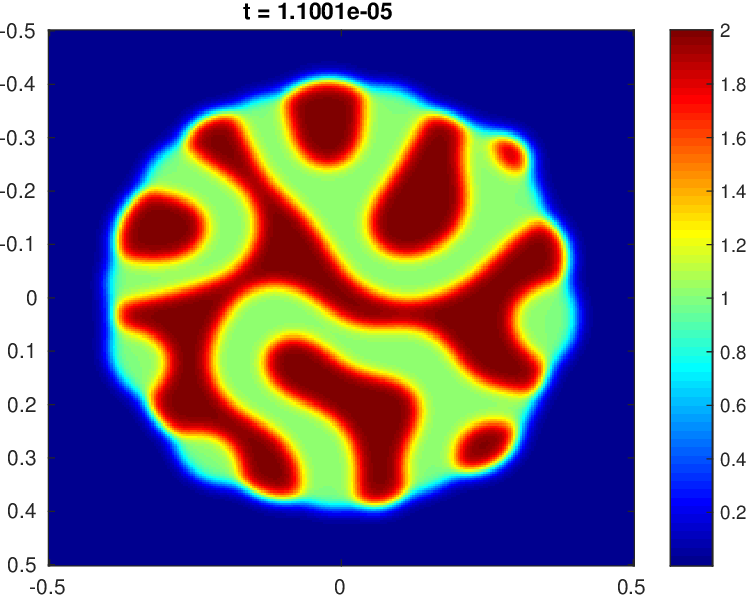}
	\includegraphics[width=3.5cm]{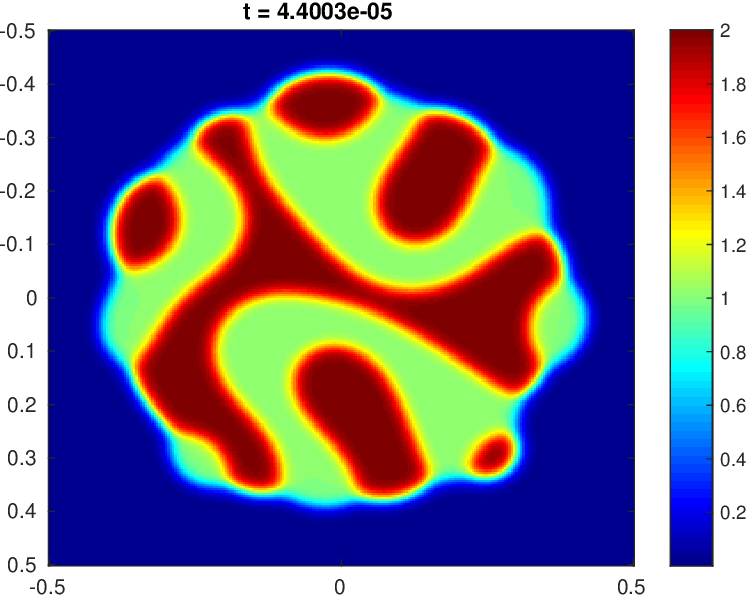}
	\includegraphics[width=3.5cm]{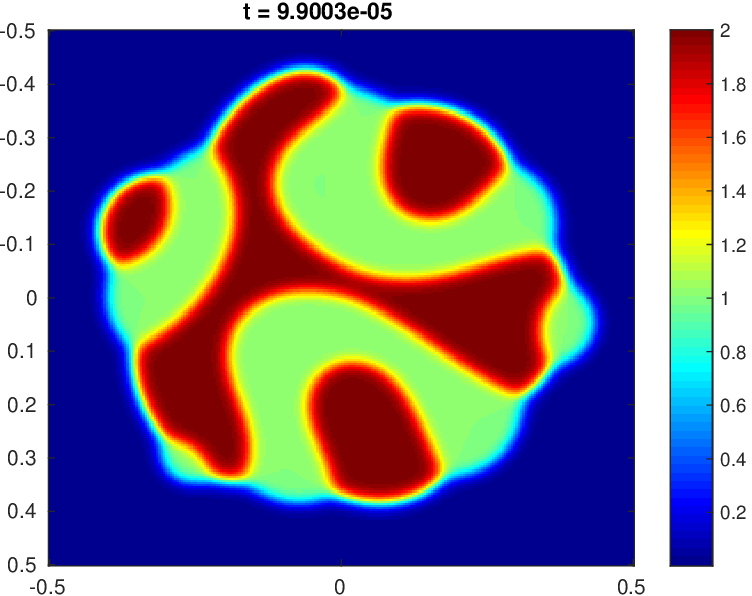} \\
    \includegraphics[width=3.5cm]{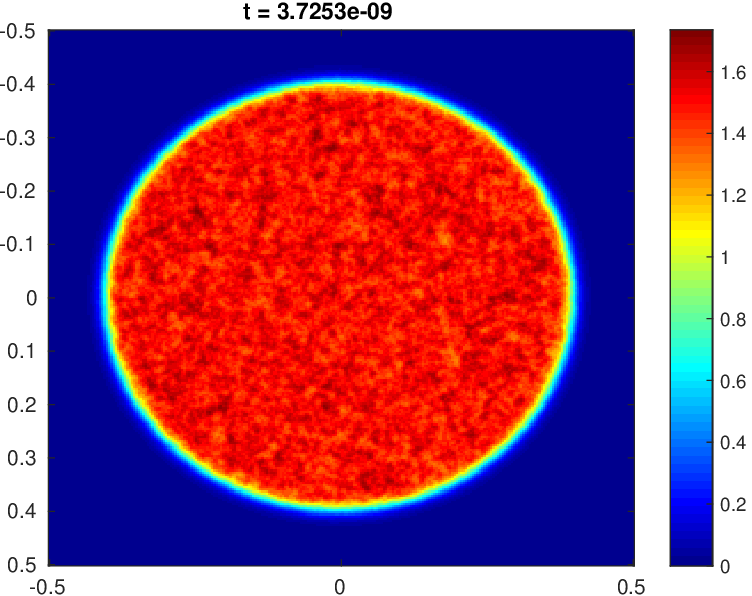}
	\includegraphics[width=3.5cm]{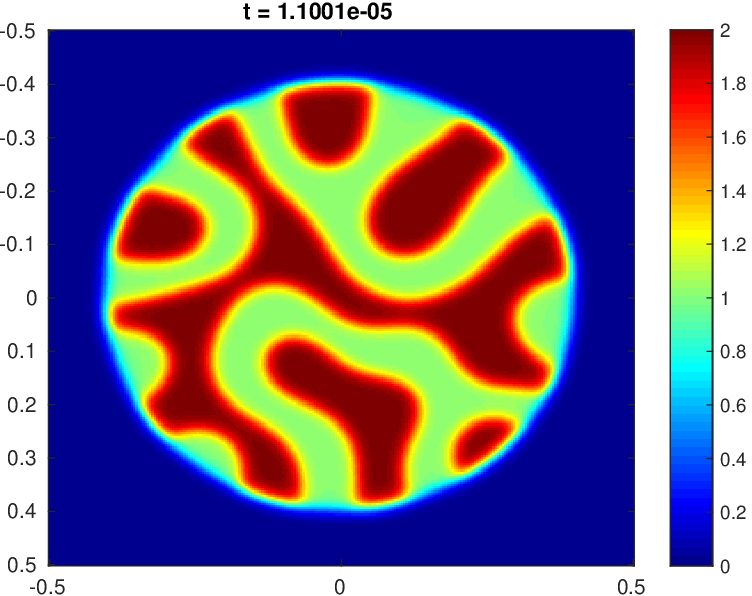}
	\includegraphics[width=3.5cm]{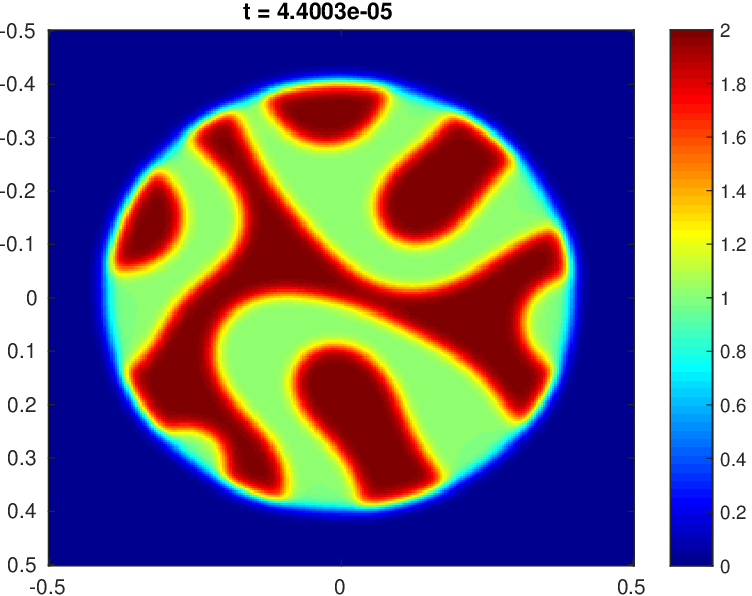}
	\includegraphics[width=3.5cm]{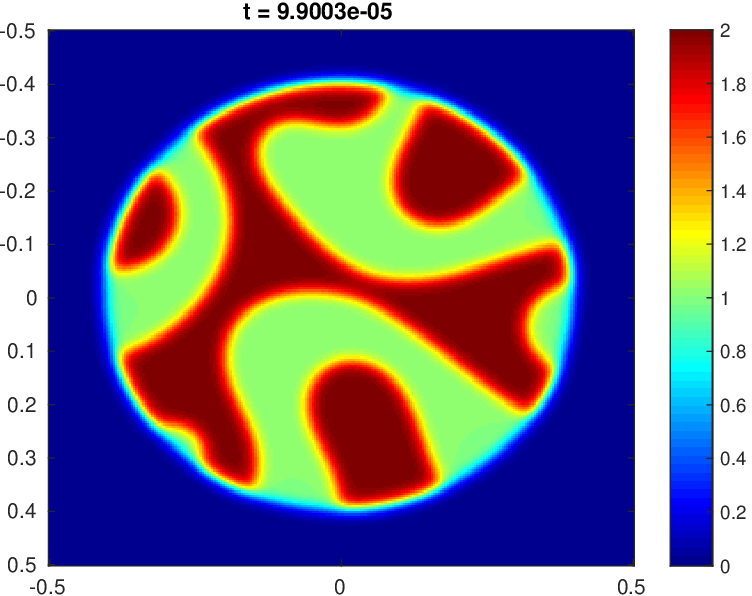} \\
\caption{Influence of the mobility coefficients using \NMNCH: evolution of ${\bf u}$ along the iterations;
First line using $\nu_1 = \nu_2 = \nu_3 = 1$; second line, using $\nu_1 = 0$ and $\nu_2 = \nu_3= 1$.}
\label{fig_test3}
\end{figure}

The second numerical experiment is intended to show the influence of surface mobilities $(\nu_{ij})$  on the velocity of each interface.
To illustrate this, we show in Figure~\ref{fig_test3} the evolution of ${\bf u}$ in two different cases: a first case where $\nu_i=1$
(see the first row on Figure~\ref{fig_test3}); a second case where $\nu_2=\nu_3=1$ and $\nu_1=0$ 
(see the second row). In both cases, the $(\sigma_i)$ coefficients associated with surface tensions
$(\sigma_{ij})$ are set to $\sigma_i=1$. As previously, we use the same numerical parameters
in each case: we set $N= 2^8$, $\varepsilon = 2/N$, $\delta_t = \varepsilon^4$,
$\alpha = 2$, $m=1$, and $\beta= 2/\varepsilon^2$. \\

As expected, we observe in the first row of Figure~\ref{fig_test3} that all phases are active along the
iterations since the mobility coefficients $\nu_i$ are all equal to $1$. 
On the contrary, in the second row, the first phase ($u_1$ in blue) is fixed along the iterations,
which is consistent with the fact that the coefficient mobility associated with the first phase $u_1$ is $\nu_1=0$.
Indeed, it is important to notice that mobilities play a role only in the gradient flow and therefore 
imposing a zero mobility $\nu_k=0$ forces the  $k$-th phase $u_k$ to be fixed. 
In particular, this allows us to deal easily and efficiently with the Cahn--Hilliard problem 
in irregular domains (see~\cite{shin_conservative_2011,li_conservative_2013,teigen_diffuse-interface_nodate,li_solving_nodate,YANG2020105276}) and the second row of Figure~\ref{fig_test3}  is a perfect illustration of it.
We insist that our model does not impose any boundary conditions on the complex domain,
nor the insertion of a surface energy. Another important remark is that the width of the diffuse interface depends only on $\varepsilon$ and does
not depend neither on surface tensions nor on mobilities.

\subsubsection{Influence of the surface tension coefficients using the  \NMNCH ~model}

\begin{figure}[htbp]
\centering
	\includegraphics[width=3.5cm]{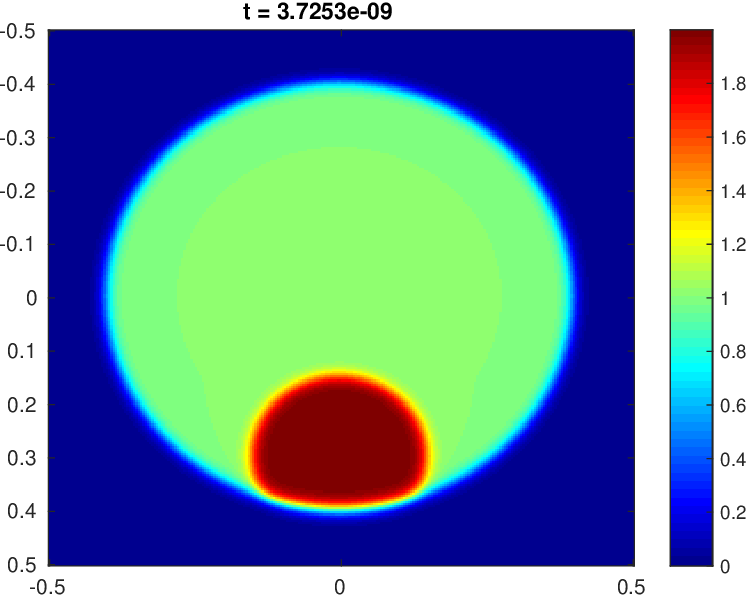}
	\includegraphics[width=3.5cm]{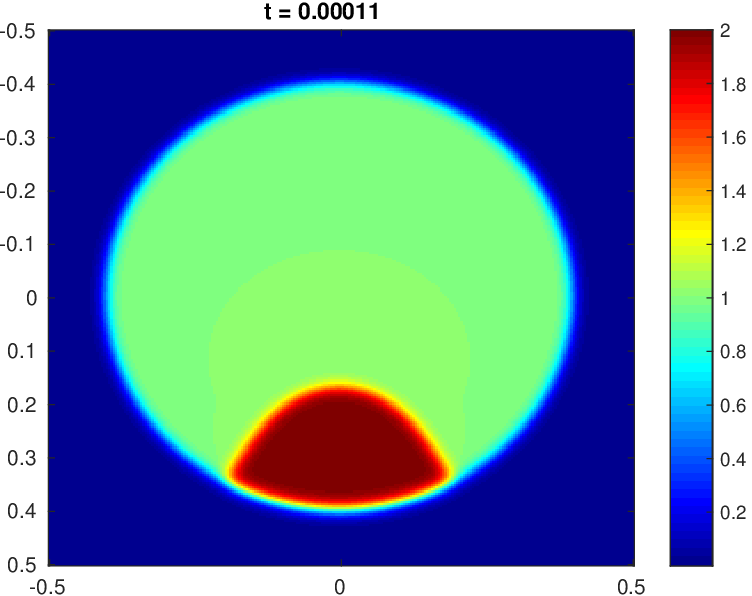}
	\includegraphics[width=3.5cm]{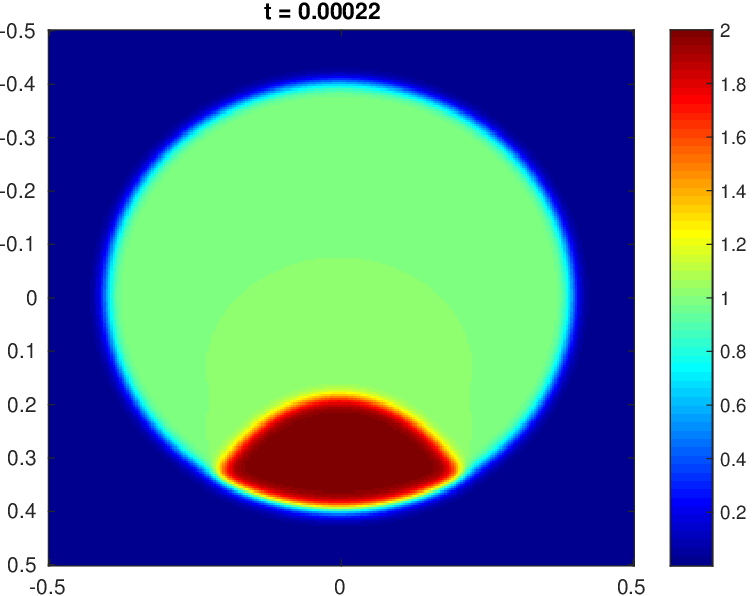}
	\includegraphics[width=3.5cm]{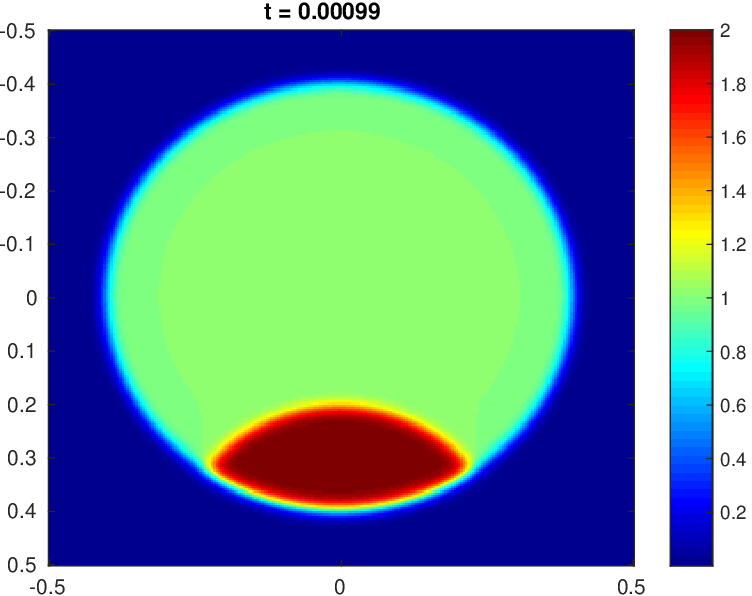} \\
	\includegraphics[width=3.5cm]{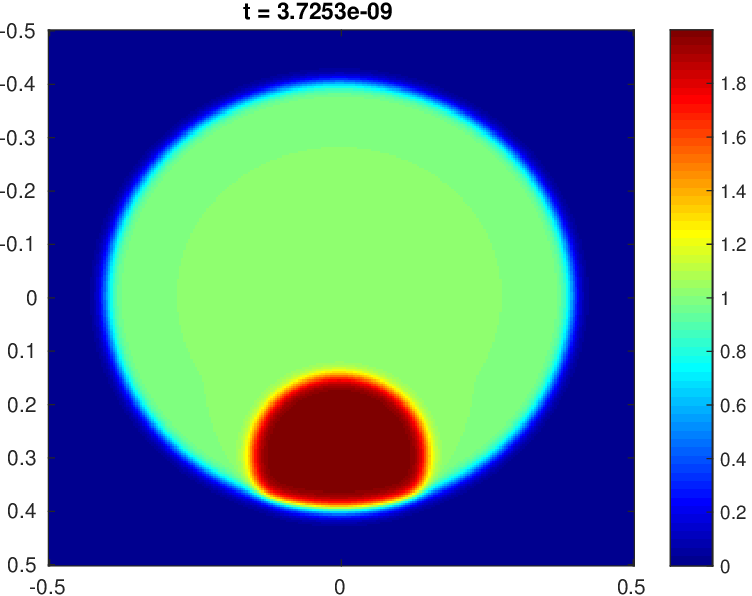}
	\includegraphics[width=3.5cm]{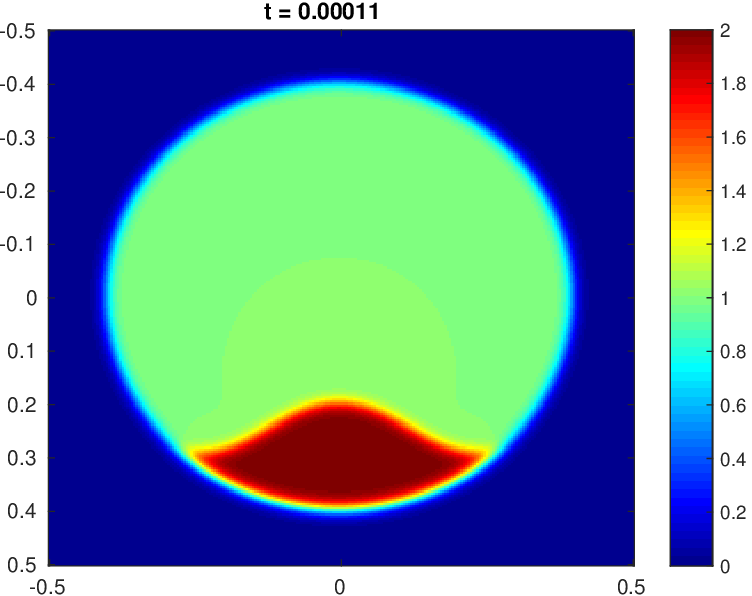}
	\includegraphics[width=3.5cm]{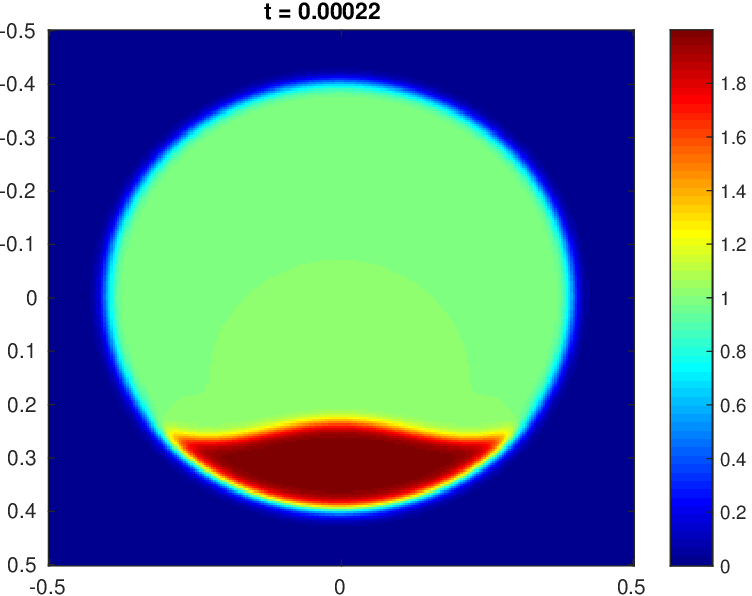}
	\includegraphics[width=3.5cm]{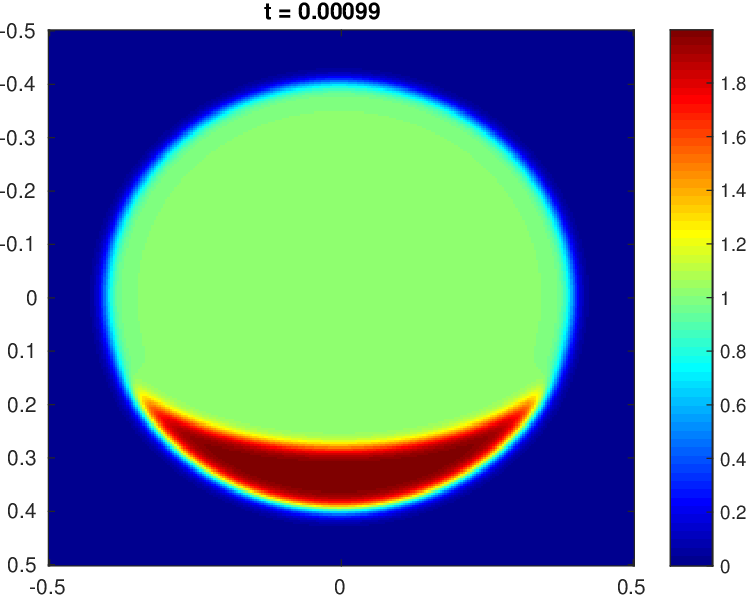} \\
		\includegraphics[width=3.5cm]{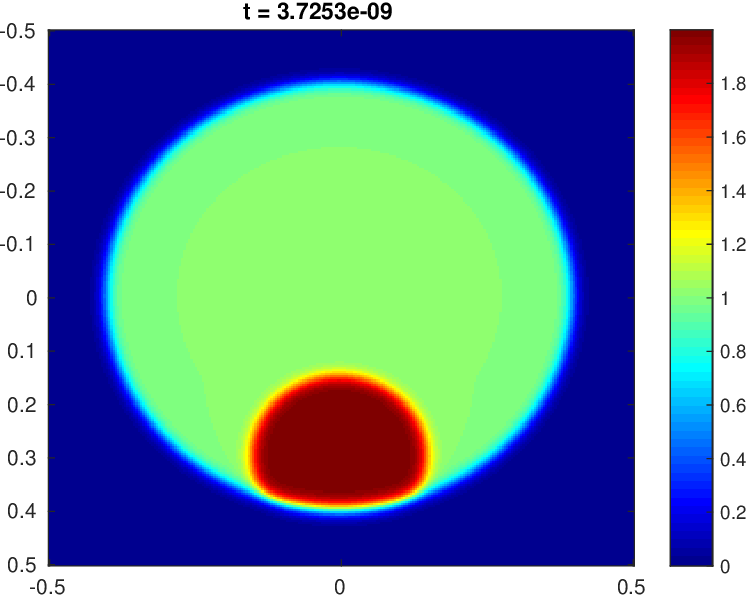}
	\includegraphics[width=3.5cm]{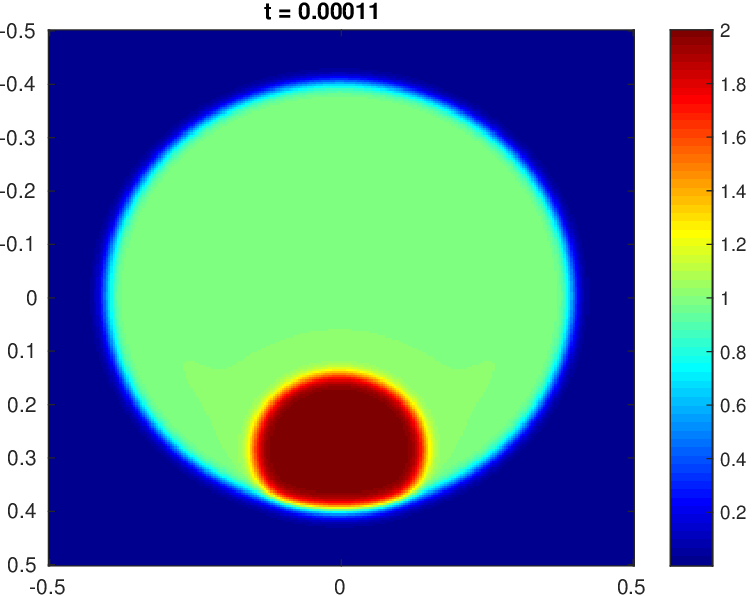}
	\includegraphics[width=3.5cm]{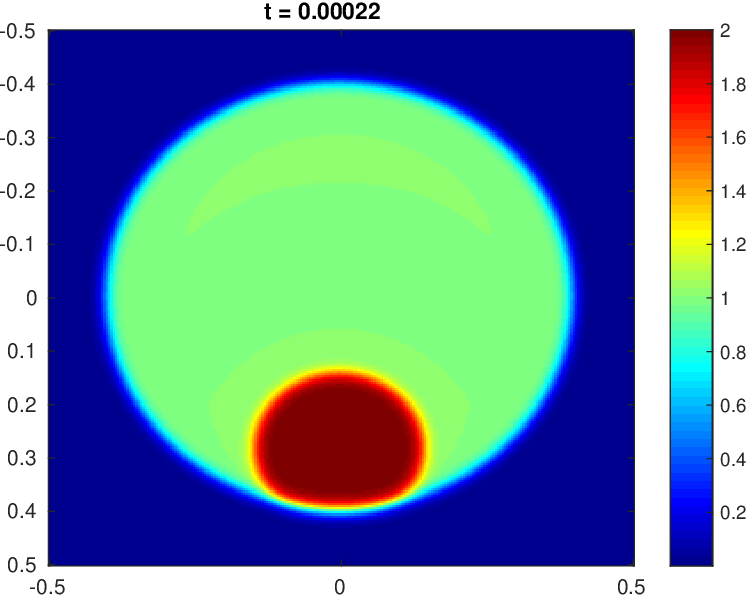}
	\includegraphics[width=3.5cm]{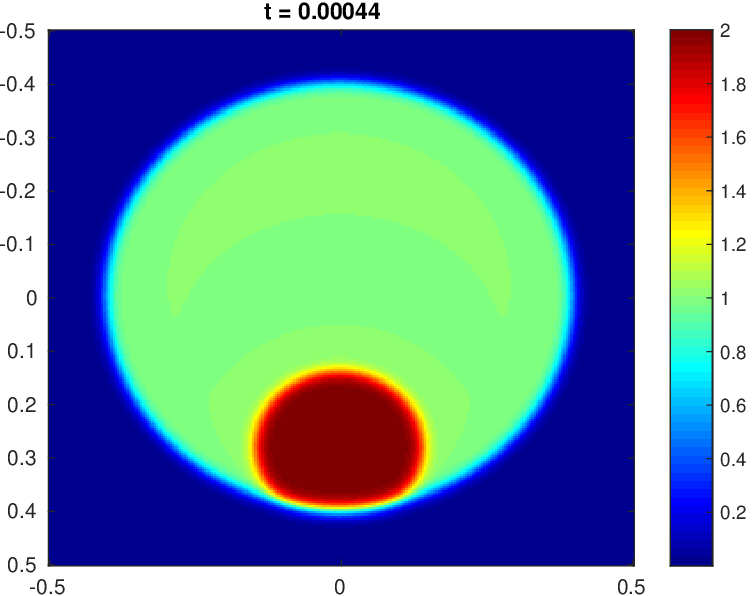} \\
\caption{Influence of the surface tension coefficients using \NMNCH: evolution of ${\bf u}$ along the iterations (the $\frac 1 2$ isolevel set of the liquid phase is plotted in gold, the $\frac 1 2$ isolevel set of the solid phase is in red);
First line using $\sigma_{12} = \sigma_{13} = \sigma_{23} = 1$; second line using $\sigma_{12} = 1.9$ and $\sigma_{13} = 
\sigma_{23}= 1$;  third line using $\sigma_{13} = 1.9$, and $\sigma_{12} = \sigma_{23}= 1$.}
\label{fig_test4}
\end{figure}

The \NMNCH ~model can also handle the case of the evolution of a liquid phase on a fixed solid surface by simply 
imposing a null mobility of the solid interface. Here we propose an application in space dimension 2. 
Figure~\ref{fig_test4} illustrates numerical results obtained with different sets of surface tension 
coefficients  $\mathbf{\sigma}=(\sigma_{12}, \sigma_{13}, \sigma_{23})$,  with mobilities $\nu_1=0,~\nu_2=\nu_3=1$ and
the same initial condition: $\sigma=(1,1,1)$, $\sigma=(1.9,1,1)$ and $\sigma=(1,1.9,1)$ for the first, the second 
and the third rows of Figure~\ref{fig_test4}, respectively. The solid $u_1$, liquid $u_2$ and vapor $u_3$ phases  are 
represented in blue, red, and green, respectively. Similarly to the previous computations, the numerical 
parameters are set to $N = 2^8$, $\varepsilon = 2 /N$, $\delta_t = \varepsilon^4$,
$\alpha = 2$, $m=1$, and $\beta= 2/\varepsilon^2$. As in the previous numerical experiment, 
we notice the ability of our model to handle the case of null mobilities 
(here to fix the exterior solid phase $u_1$ in blue).
In Figure~\ref{fig_test4}, we can also see the strong influence of the contact angle on the evolution of the liquid phase. 
We emphasize that our model does not prescribe the contact angle.  Rather, its value is an implicit consequence
of the multiphase interface energy considered in each simulation. 


\section{Application to the numerical simulation of dewetting}

{Various approaches have been proposed to approximate numerically the dewetting process as a constrained surface diffusion flow, see e.g. the papers \cite{wetting_physique,JIANG201969,Jiang2012PhaseFA,Dziwnik_2017,doi:10.1137/19M1251345,PhysRevB.91.045303} devoted to the simulation of solid-state dewetting based on either isotropic or anisotropic surface diffusion (anisotropic models are closer to physical reality but isotropic ones can be used as a first approximation). In particular, two phase field approaches are proposed in \cite{Jiang2012PhaseFA,Dziwnik_2017}. In these papers, an isotropic~\cite{Jiang2012PhaseFA} or anisotropic~\cite{Dziwnik_2017} Cahn-Hilliard model is coupled with surface energies that encode the contact angle on the support solid phase. As will be shown, the model proposed in the current paper also involves the isotropic Cahn-Hilliard equation and it can be used to simulate dewetting. However, our model makes no explicit reference to the contact line thanks to a multiphase formulation. In addition, we use degenerate mobilities which improve the order of approximation of the phase field model and bring more accurate numerical approximations.}\\

Other approaches have been developed to simulate either wetting without convection or dewetting. For example Cahn proposed in~\cite{cahn1977critical} a phase-field approach 
with an additional surface energy on the boundary
of the solid phase but the method is only applicable for a contact angle $\theta < \frac{\pi}{2}$.
An approach coupling the Allen-Cahn equation and smoothed boundary conditions to force the correct contact angle condition
is available in \cite{turco2009wetting}.
Other methods based on the Allen-Cahn equation and using this idea are proposed 
in \cite{ben2014phase,diewald2018investigating}. Alternative methods using wall boundary conditions with 
a third order polynomial to impose the contact angle 
are proposed in \cite{sibley2013contact,sibley2013moving,aymard2019linear}.  
A convexity splitting scheme using this idea with a sinusoidal boundary condition can be found in \cite{wei2020integral}.
In \cite{dong2012imposing,metzger2015numerical,carlson2011dissipation,shen2015efficient,boyer2017ddfv} 
the angle is imposed using wall boundary conditions again.
The dynamic case {with convection} can be treated via a coupled Cahn--Hilliard/Navier-Stokes system. In most cases, see for example
\cite{jacqmin1999calculation,liu2003phase,turco2009wetting,abels2009diffuse,boyer2010cahn},
the contact angle is set to the static contact angle $\frac{\pi}{2}$. 

In the convolution-thresholding framework, some recent approaches have been proposed to simulate the dewetting process or wetting without convection.
Expanding the original scheme of Bence, Merriman, and Osher \cite{merriman1992diffusion}, 
Esedoglu and Otto have proposed in \cite{esedog2015threshold} a multiphase convolution-thresholding method for arbitrary surface tensions satisfying the triangle inequality. Wang et al.~\cite{wang2019improved} then applied this generalization to the wetting case. A different approach proposed in \cite{xu2019adaptive} does not impose 
the contact angle in the formulation but requires the use of sophisticated techniques while solving the heat equation.

In \cite{bretin2017new,bretin2018multiphase}, two authors of the present paper proposed an Allen-Cahn equation coupled with a frozen solid phase to approximate droplet {dewetting} (or wetting without convection).
It was based on the use of zero surface mobilities  for the solid-vapor and solid-liquid interfaces. 
In this paper, we extend this idea to the Cahn--Hilliard equation and, coupled with a reformulation
of the problem, we introduce a new simple and effective method for simulating the dewetting phenomenon. 
It is important to emphasize that this method does not impose the contact angle, which is determined implicitly
by the surface tension coefficients $(\sigma_{SV},\sigma_{SL},\sigma_{LV})$.

\subsection{Rewriting of the model using the liquid phase only} 

We consider for the simplicity of presentation a liquid-solid-air dewetting situation but other situations could be considered equally. We will focus in particular on the rather difficult simulation of thin liquid tubes dewetting on arbitrary solid surfaces. Numerical simulation of dewetting in dimension $3$ with a complete model $(u_L,u_V,u_S)$ can be quite challenging numerically, it is therefore preferable to reduce the system. As the solid phase 
$\Omega_S$ is fixed and $\Omega_V$ can be obtained from $\Omega_L$, $\Omega_S$, it is possible to consider only one unknown, the liquid phase, and to use a reduced phase field model involving this phase only.\\

The Cahn-Hilliard energy reads as 
$$P_{\varepsilon}({\bf u}) = \sum_{k \in \{S,L,V\}} \frac{\sigma_k}{2} \int_{Q} \frac{\varepsilon}{2} |\nabla u_k|^2 + \frac{1}{\varepsilon} W(u_k),$$
where 
$$\sigma_L = \frac{\sigma_{LS} + \sigma_{LV} -  \sigma_{SV}}{2}, \,\sigma_S = \frac{\sigma_{LS} + \sigma_{SV} -  \sigma_{LV}}{2}\;\;
\text{ and }\;\; \sigma_V = \frac{\sigma_{LV} + \sigma_{SV} -  \sigma_{LS}}{2}.$$ 
Here, $u_S$ represents the phase field function associated with the solid set $\Omega_S$ and
the previous asymptotic developments show that $u_S$ should be of the form 
$u_S = q \left( \frac{\operatorname{dist}(x, \Omega_S)}{\varepsilon}\right)$. 
On the other hand, the vapor phase field function $u_V$ can be expressed from the partition constraint as 
$ u_V = 1 - (u_S + u_L)$.
Then the Cahn--Hilliard energy can be rewritten using only the variable $u_L$ as follows: 
\begin{eqnarray*}
\tilde{P}_{\varepsilon}(u_L) &=&  \frac{\sigma_{L}}{2} \int_{Q} (\frac{\varepsilon}{2} |\nabla u_L|^2 + \frac{1}{\varepsilon} W(u_L))dx \\
&+& \frac{\sigma_{V}}{2}\int_{Q} (\frac{\varepsilon}{2} |\nabla (1 - (u_S +  u_L)|^2 + \frac{1}{\varepsilon} W(1 - (u_L + u_S)))dx \\
&+&  \frac{\sigma_{S}}{2}\int_{Q} (\frac{\varepsilon}{2} |\nabla u_S|^2 + \frac{1}{\varepsilon} W(u_S))dx. \\
\end{eqnarray*}
Notice that its $L^{2}$-gradient satisfies 
\begin{eqnarray*}
 \nabla_{L^2} \tilde{P}_{\varepsilon}(u_L) &=&  \frac{\sigma_{SL}}{2} [- \varepsilon \Delta u_L  + \frac{1}{\varepsilon} W'(u_L)] + \frac{\sigma_{V}}{2} \varepsilon R_{u_S}(u_L) 
\end{eqnarray*}
where the first term
$$ \frac{\sigma_{SL}}{2} [- \varepsilon \Delta u_L  + \frac{1}{\varepsilon} W'(u_L)],$$
is a classical Allen-Cahn term  and the second  term 
$$ R_{u_S}(u_L) =    - \left[   \Delta u_S  +  \frac{1}{\varepsilon^2} (W'(u_L) + W'(1 - (u_L + u_S))  \right],$$
appears as a smooth penalization term which is active only on the boundary of $\Omega_S$. 

Finally, incorporating mobilities leads us to consider the following Cahn--Hilliard models:
\begin{itemize}
 \item {\bf \MCH ~model} \\
$$
 \begin{cases}
  \partial_t u_L &=  \div\left(M(u_k)\nabla( \sigma_{LV}/2  \mu_L + \sigma_V R_{u_S}(u_L) ) \right) \\
   \mu_L & = \frac{W'(u_L)}{\varepsilon^2} - \Delta u_L \\
 \end{cases}
 $$
 \item {\bf \NMNCH ~model} \\
$$
 \begin{cases}
    \partial_t u_L &=  N(u_L) \div\left(M(u_L)\nabla( N(u_L) \left(  \sigma_{LV}/2  \mu_L + \sigma_V R_{u_S}(u_L) \right) ) \right) \\
   \mu_L & = \frac{W'(u_L)}{\varepsilon^2} - \Delta u_L \\
 \end{cases}
 $$
  \end{itemize}
Note that, in practice, we used the \NMNCH ~model for all numerical simulations presented hereafter because the dewetting 
of a thin structure requires a model as accurate as possible.\\

The simulations are performed with the numerical scheme introduced earlier with an additional explicit treatment of the penalization term
$R_{u_S}(u_L)$. \\

Notice that the penalization term $R_{u_S}(u_L)$ is active on the whole boundary of $\Omega_S$. In particular, 
when $u_L = 0$ this term is still active and can be important as it corresponds to the Allen-Cahn term associated to $u_S$:  
$$R_{u_S}(u_L) = - \left(\Delta u_S + \frac{1}{\varepsilon^2}W'(1-u_S)\right) = - \Delta u_S - \frac{1}{\varepsilon^2} W'(u_S).$$
In practice, we propose to localize it only at the liquid phase boundary $u_L$, which can be done by considering 
the following variant
$$\tilde{R}_{u_S}(u_L) =  R_{u_S}(u_L) \frac{ \sqrt{2 W(u_L)} }{\sqrt{2 W(u_L) + \varepsilon }}.$$
This variant is interesting for it contributes to stabilizing the numerical scheme 
without disturbing the evolution of the liquid phase.
 
 \subsection{Influence of the surface tension coefficients}
 
\begin{figure}[htbp]
\centering
	\includegraphics[width=3.8cm]{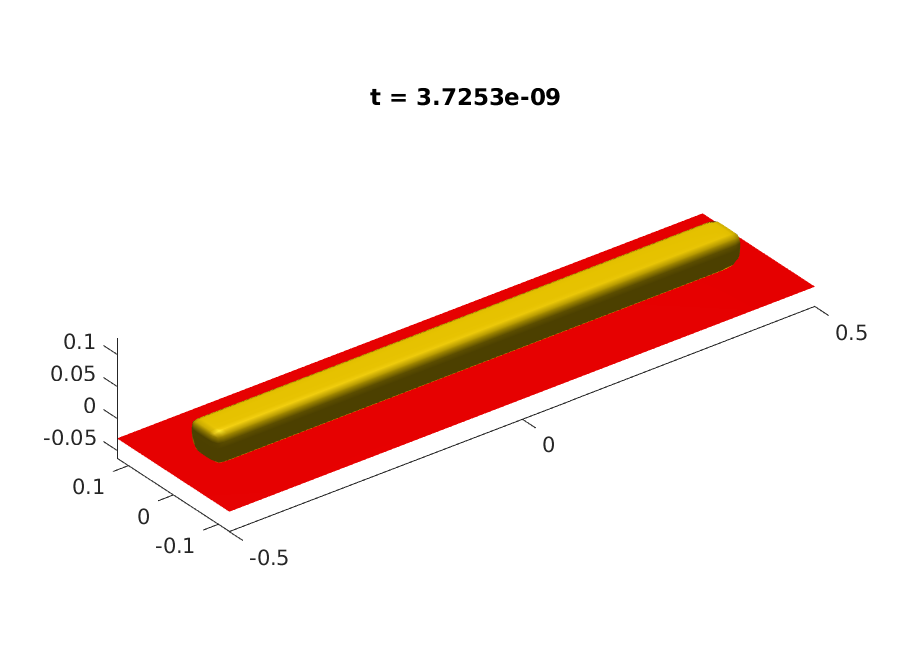}
	\includegraphics[width=3.8cm]{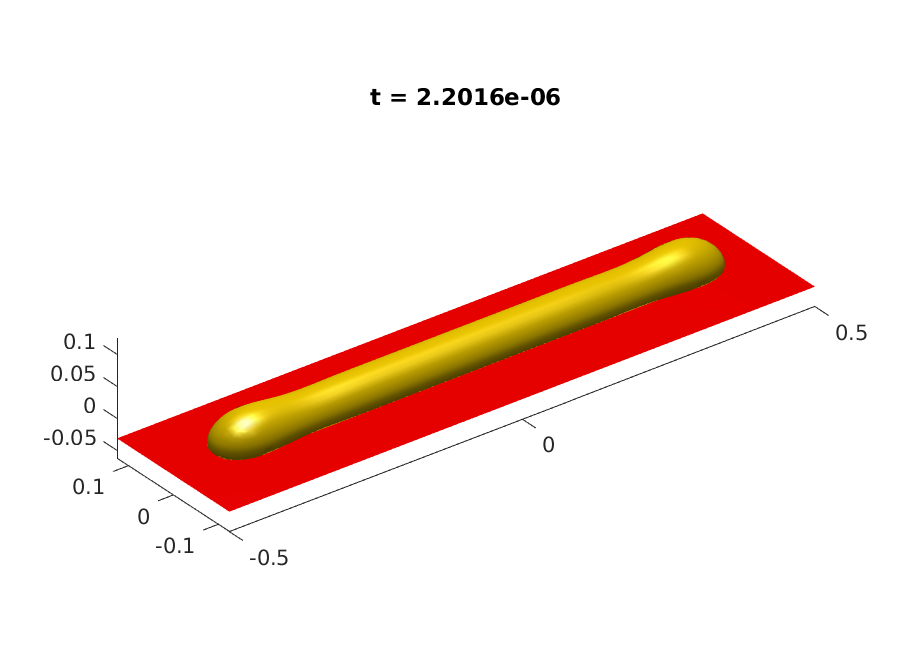}
	\includegraphics[width=3.8cm]{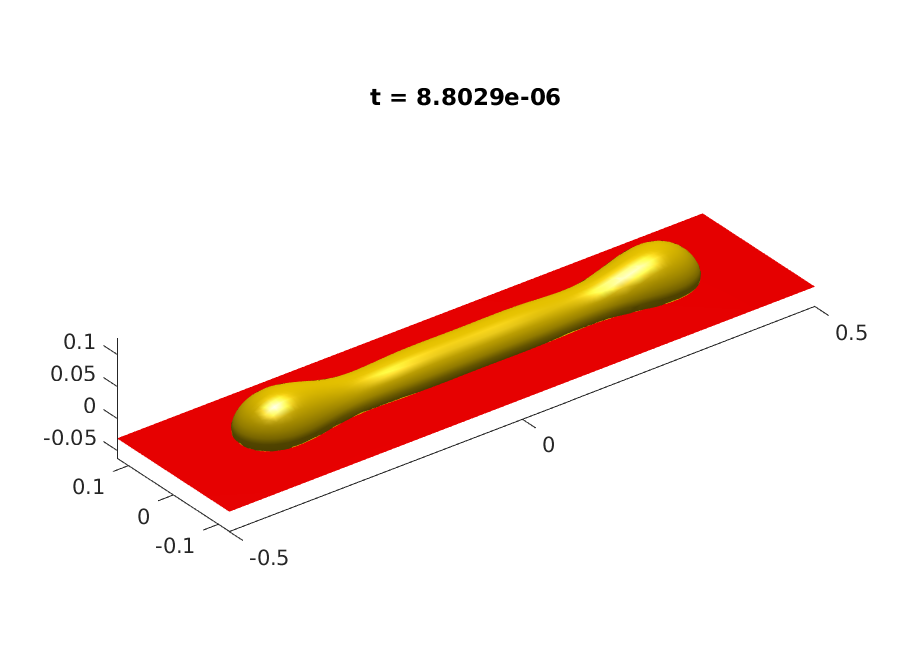}
	\includegraphics[width=3.8cm]{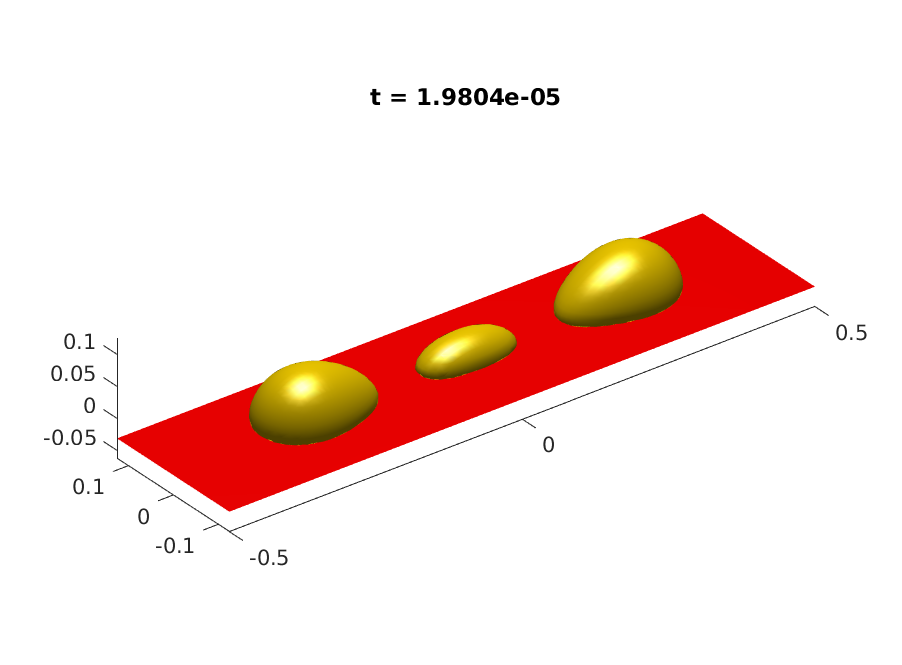} \\
	\includegraphics[width=3.8cm]{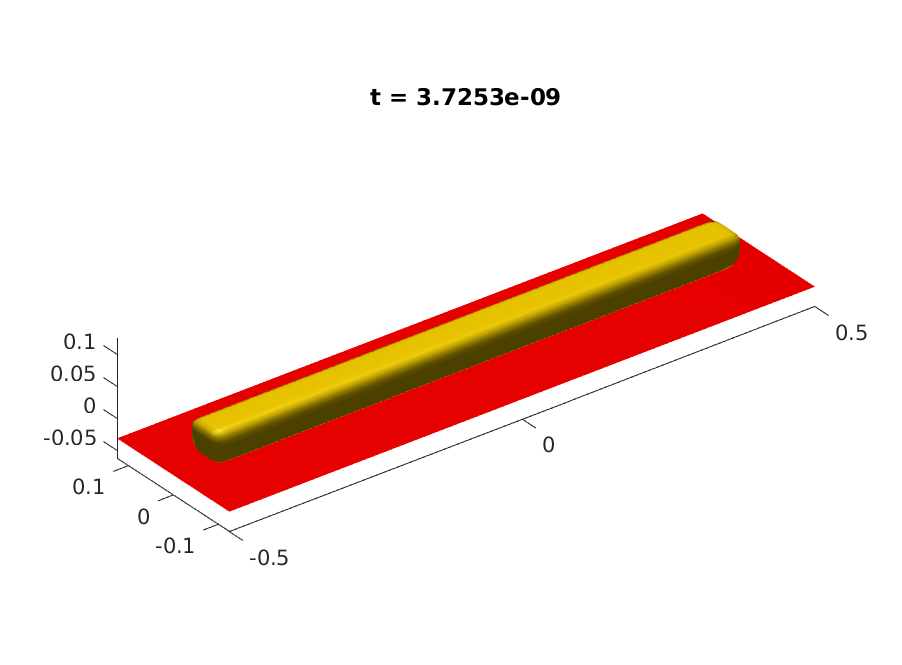}
	\includegraphics[width=3.8cm]{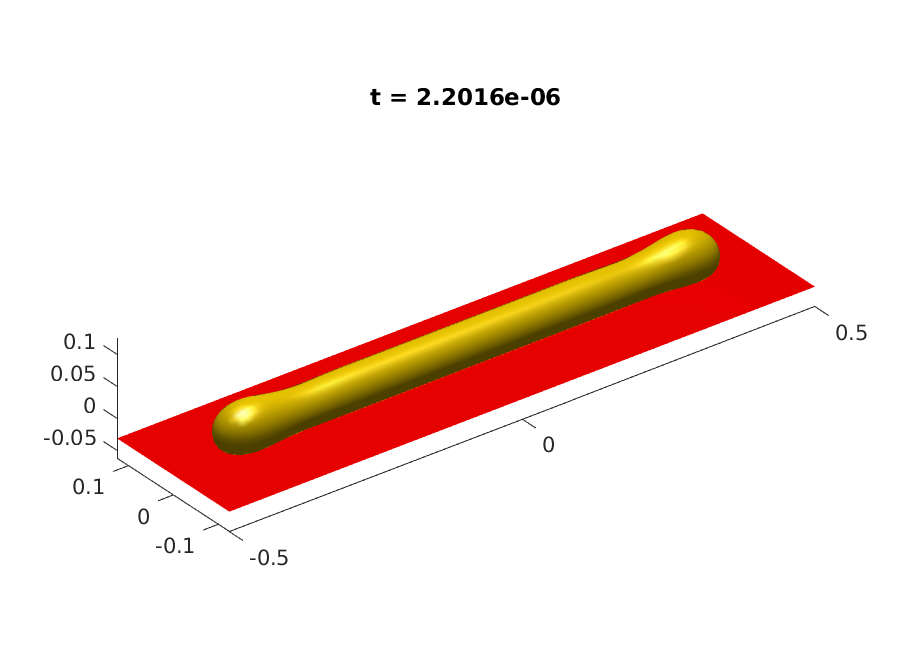}
	\includegraphics[width=3.8cm]{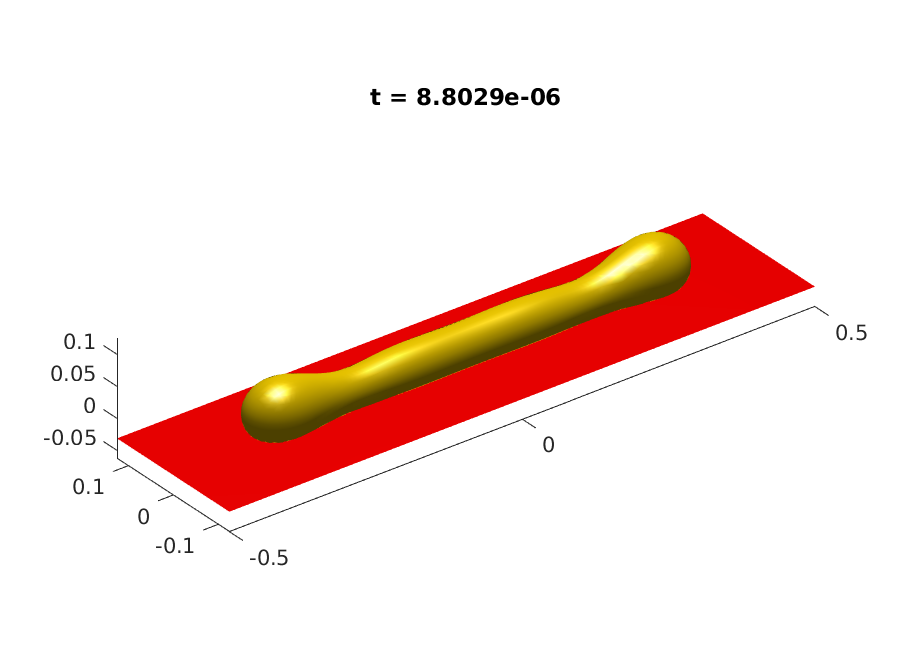}
	\includegraphics[width=3.8cm]{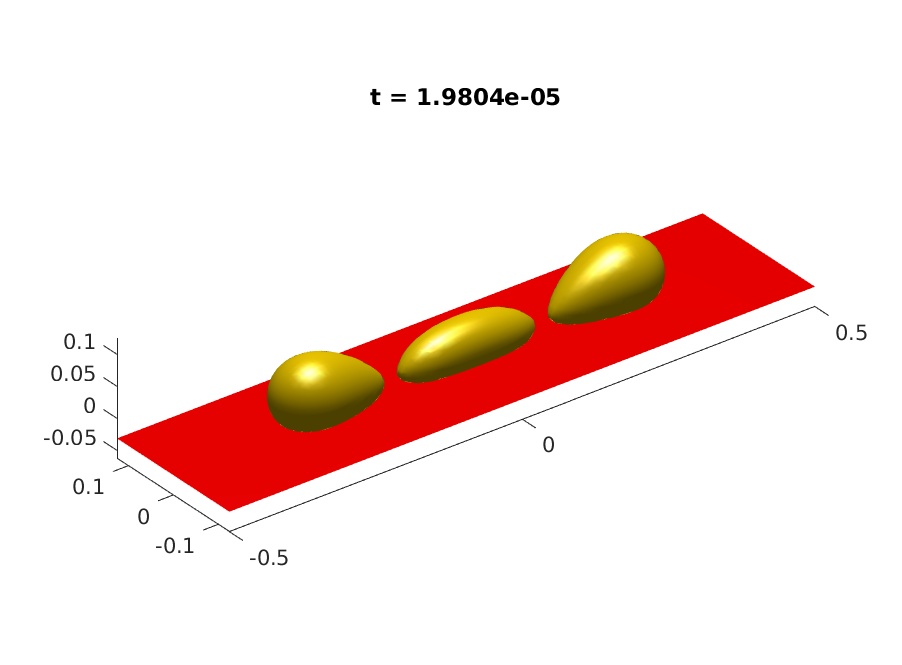} \\
	\includegraphics[width=3.8cm]{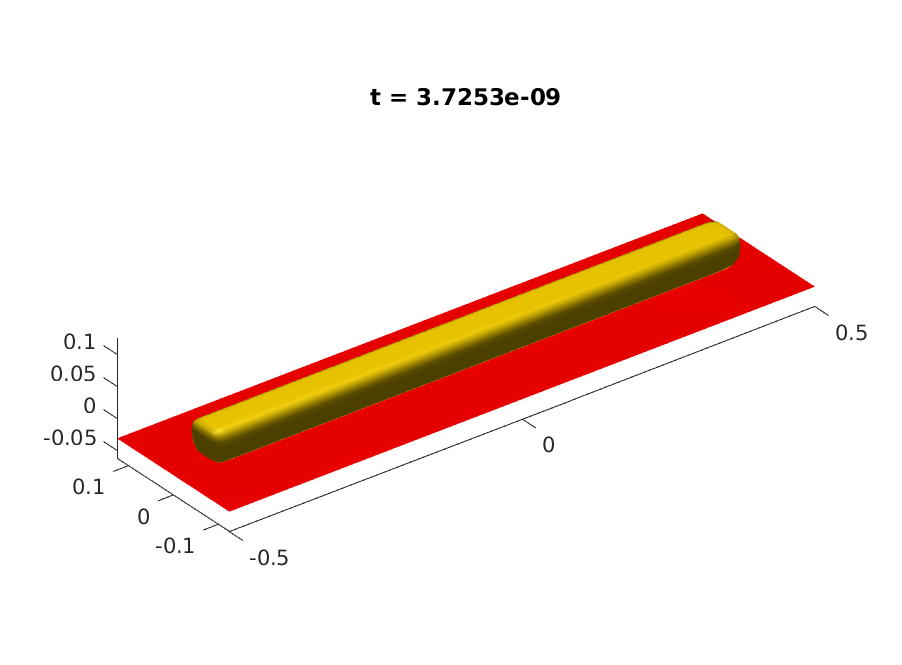}
	\includegraphics[width=3.8cm]{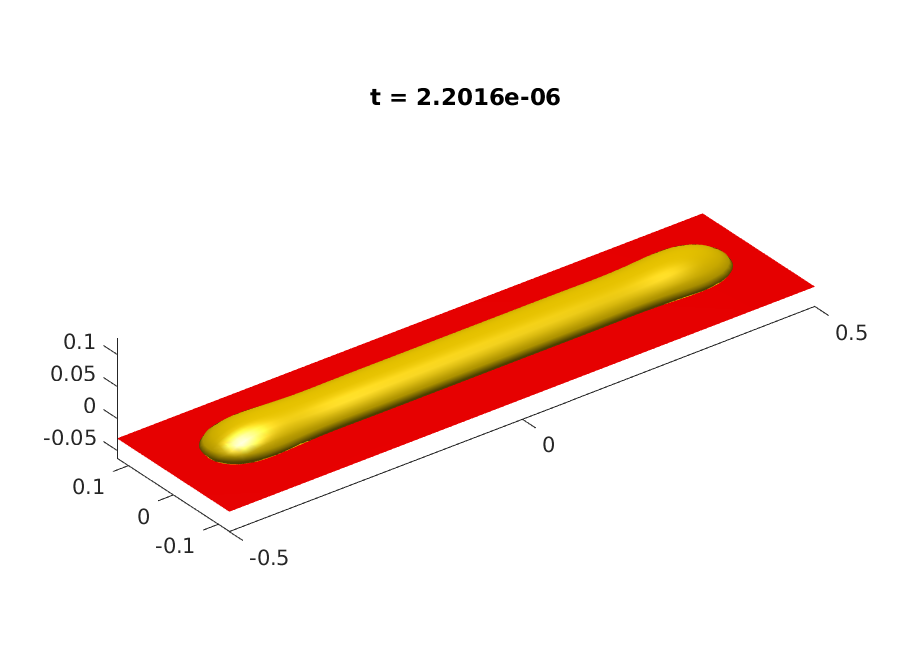}
	\includegraphics[width=3.8cm]{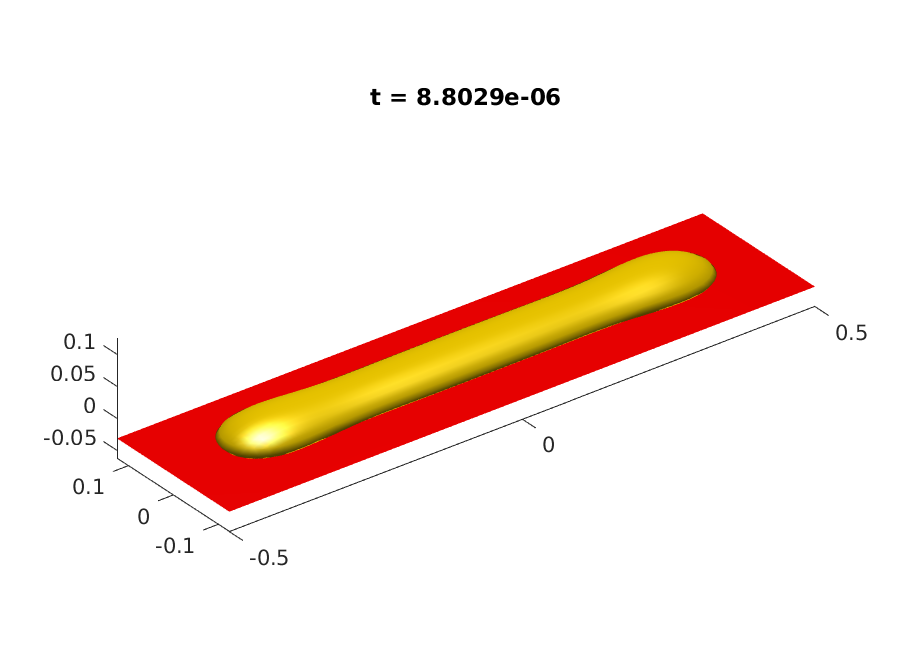}
	\includegraphics[width=3.8cm]{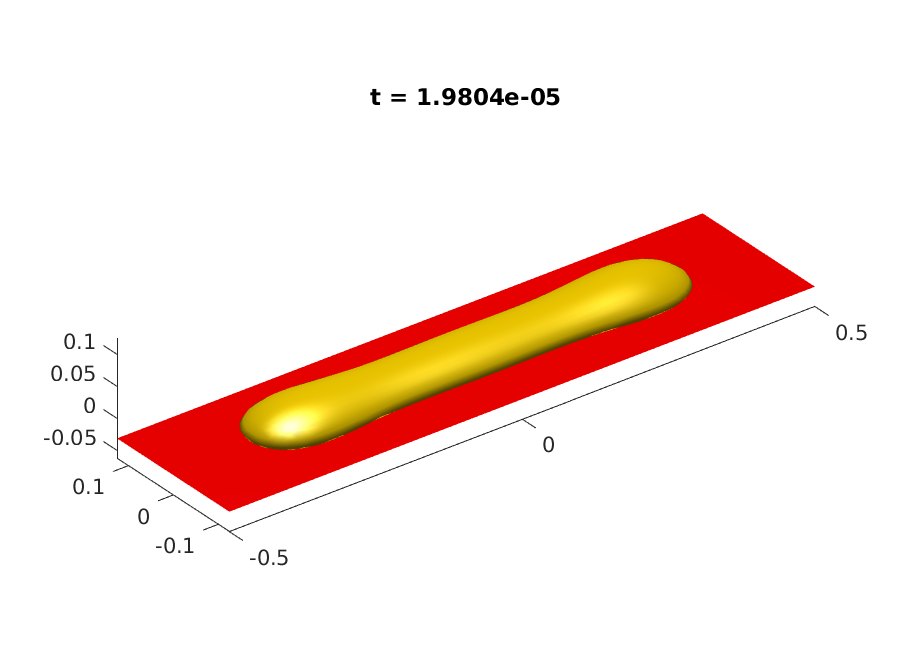} \\
\caption{Influence on dewetting of the surface tension coefficients: evolution along the iterations of the solution to the \NMNCH ~model for two different sets of coefficients; first line using $\sigma_{LV} = \sigma_{VS} = \sigma_{SL} = 1$; second line using $\sigma_{LS} = 1.7$ and $\sigma_{LV} = 
\sigma_{VS}= 1$;  third line, using $\sigma_{VS} = 1.7$, and $\sigma_{LV} = \sigma_{LS}= 1$.}
\label{fig_test5}
\end{figure}

We now propose a numerical experiment in dimension 3  where the initial set is a thin tube.
The numerical parameters are given by $N = 2^8$, $\varepsilon =  1/N$, $\delta_t = \varepsilon^4$,
$\alpha = 2$, $m=1$, and $\beta= 2/\varepsilon^2$. 
We plot on each image of Figure~\ref{fig_test5} the solution $\bf{u}$ calculated at different times $t$ with
the solid and liquid phase boundaries plotted in red and gold, respectively.
As in the 2D case, surface tension coefficients have a considerable 
influence on the evolution of the liquid phase.
They affect both the {dewetting} rate and the final shape of the liquid phase.

 \subsection{Influence of the roughness of the solid support} 
\begin{figure}[htbp]
\centering
		\includegraphics[width=3.8cm]{image/3D_1phase/Test_13D_sigma2_bord1_1}
	\includegraphics[width=3.8cm]{image/3D_1phase/Test_13D_sigma2_bord1_3}
	\includegraphics[width=3.8cm]{image/3D_1phase/Test_13D_sigma2_bord1_6}
	\includegraphics[width=3.8cm]{image/3D_1phase/Test_13D_sigma2_bord1_11} \\
		\includegraphics[width=3.8cm]{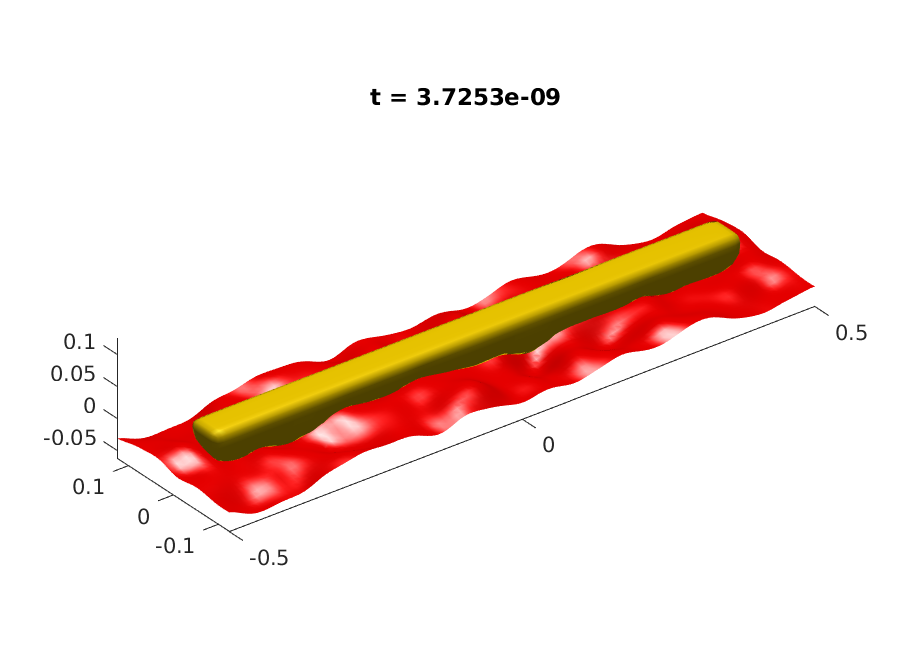}
	\includegraphics[width=3.8cm]{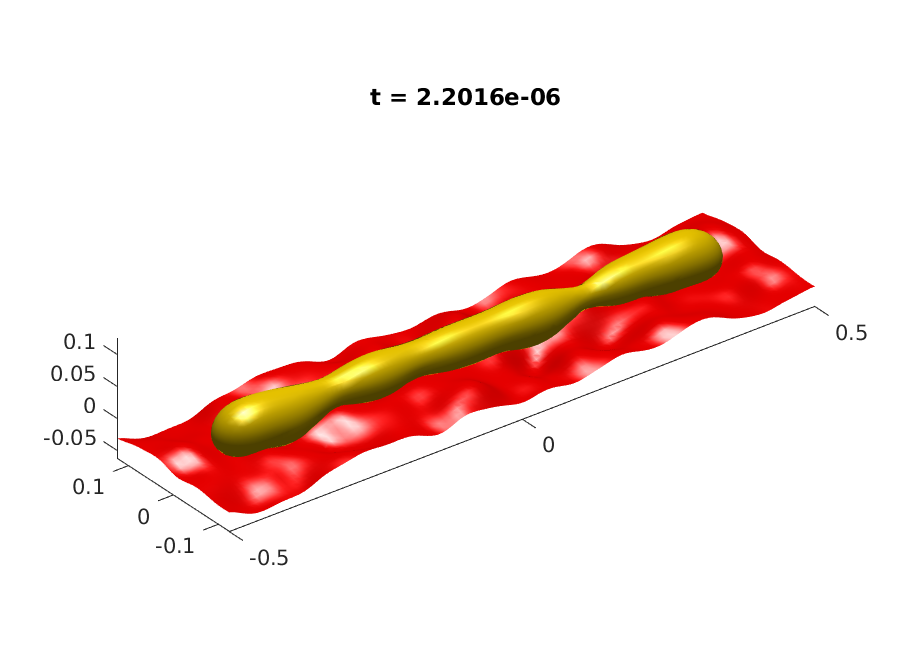}
	\includegraphics[width=3.8cm]{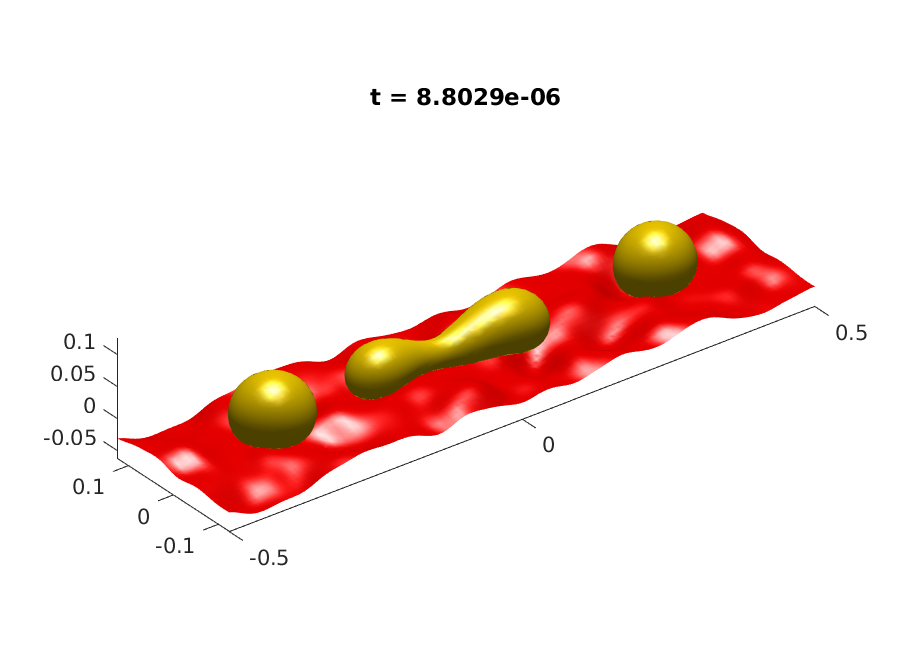}
	\includegraphics[width=3.8cm]{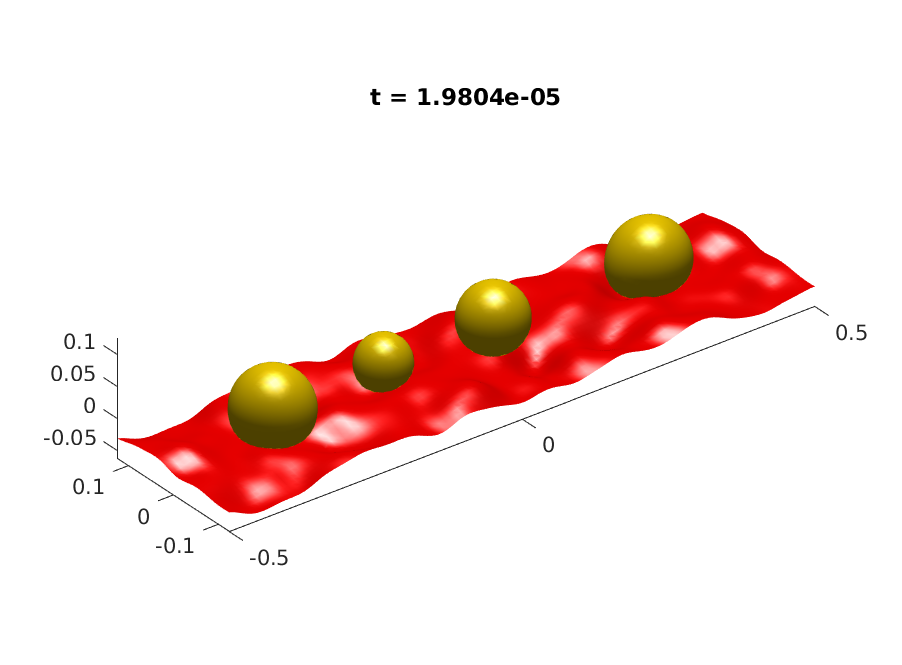} \\
		\includegraphics[width=3.8cm]{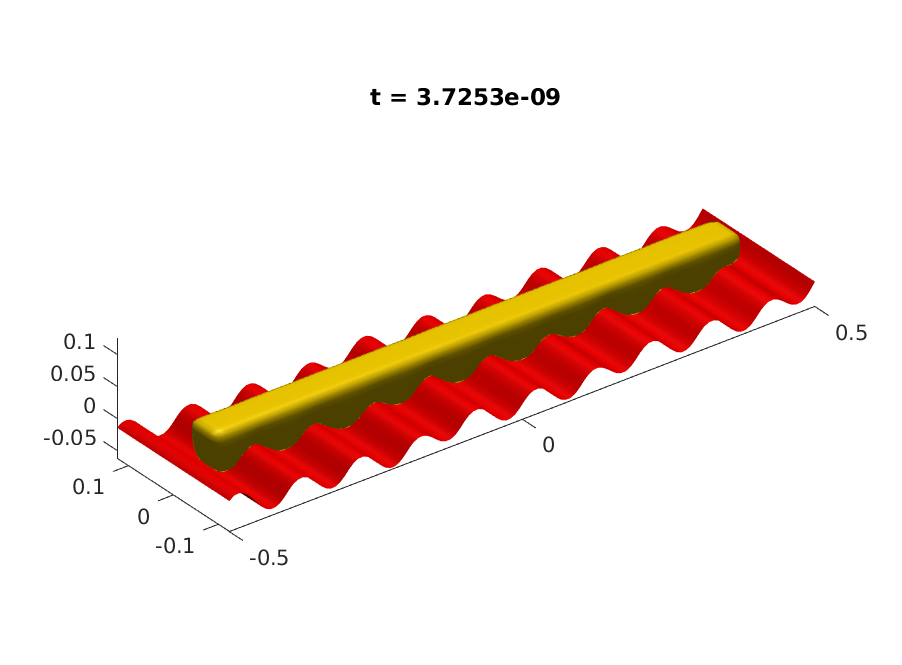}
	\includegraphics[width=3.8cm]{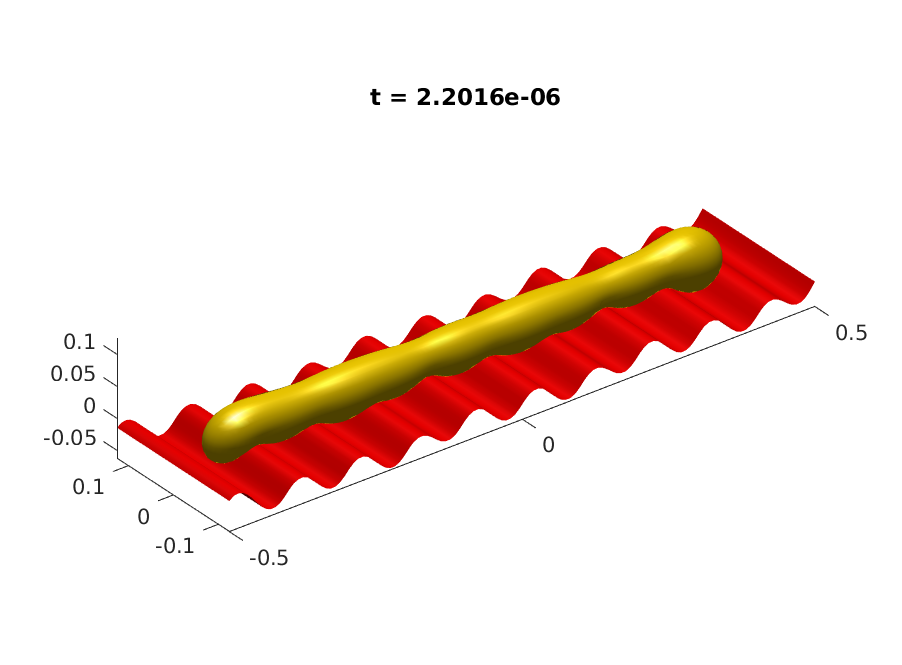}
	\includegraphics[width=3.8cm]{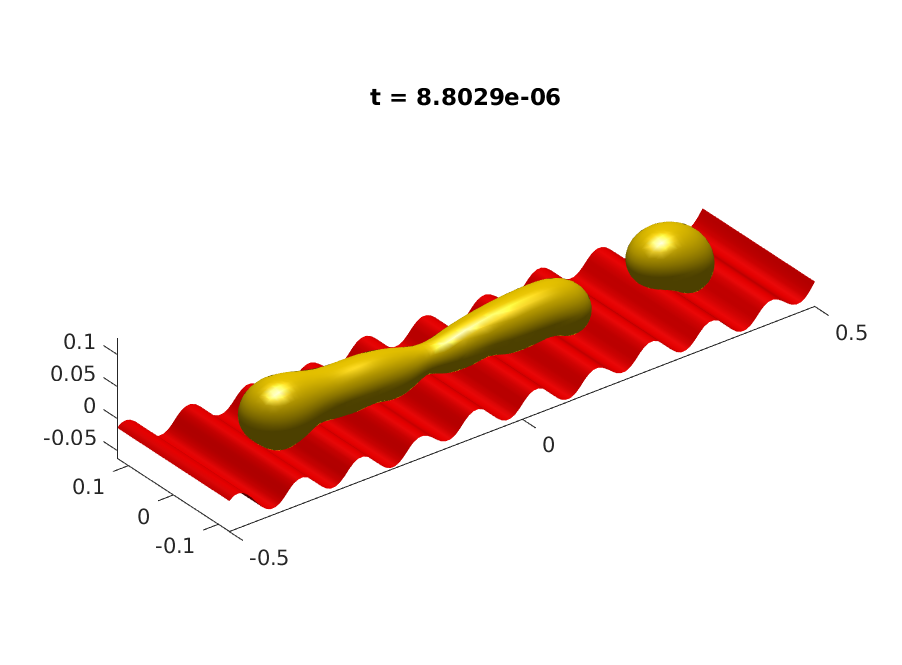}
	\includegraphics[width=3.8cm]{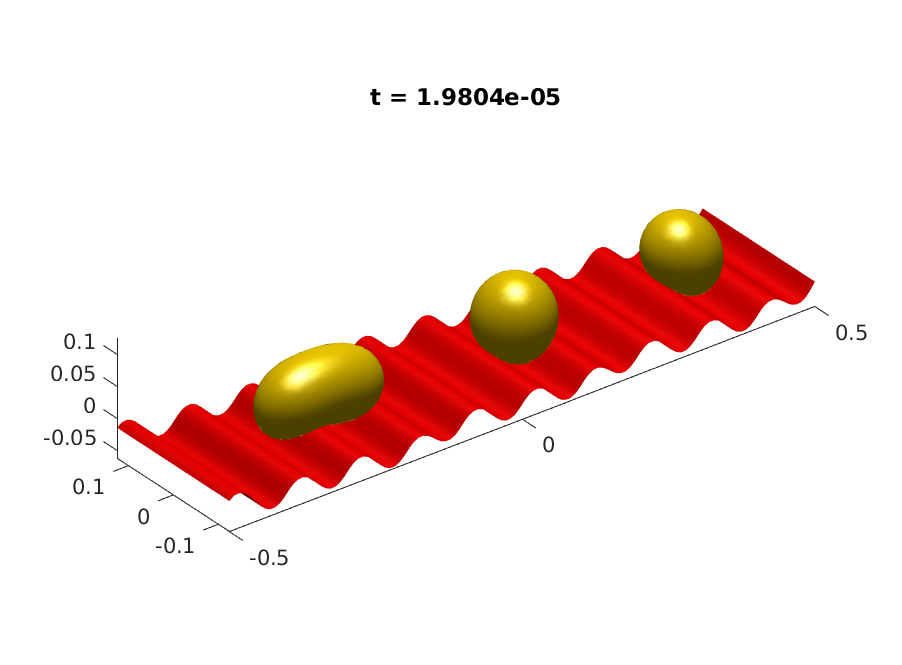} \\
\caption{Influence on dewetting of the roughness of the solid support: evolution along the iterations of the solution to the \NMNCH ~model for three different supports, with $\sigma_{LS} = 1.7$ and $\sigma_{LV} =\sigma_{VS}= 1$.}
\label{fig_test6}
\end{figure}

Our approach is also well suited for handling solid supports with roughness, i.e., notably difficult configurations for the simulation of {dewetting}. In Figure~\ref{fig_test6},  we test the case of a classical flat support, a support with randomly generated roughness, and an oscillating support.
 We observe a direct influence of the substrate roughness on the {dewetting} dynamics, each simulation being initialized in a similar way and using the same set of coefficients. 
 
\section*{Acknowledgment}
The authors acknowledge support from the French National Research Agency (ANR) under grants ANR-18-CE05-0017 (project BEEP) and ANR-19-CE01-0009-01 (project MIMESIS-3D). Part of this work was also supported by the LABEX MILYON (ANR-10-LABX-0070) of Universit\'e de Lyon, within the program "Investissements d'Avenir" (ANR-11-IDEX- 0007) operated by the French National Research Agency (ANR), and by the European Union Horizon 2020 research and innovation programme under the Marie Sklodowska-Curie grant agreement No 777826 (NoMADS).

\end{document}